%% file: main.tex
\documentclass[11pt]{amsart}
\usepackage{amssymb,amsmath,amsthm,amsfonts,mathrsfs}
\usepackage[hmargin=3cm,vmargin=3.5cm]{geometry}
\usepackage[dvipsnames,table,xcdraw]{xcolor}
\usepackage[colorlinks = true,
            linkcolor = WildStrawberry,
            urlcolor  = ForestGreen,
            citecolor = Purple,
            anchorcolor = blue]{hyperref}
\usepackage{graphicx,psfrag}
\usepackage{amscd}  
\usepackage[shellescape]{gmp}
\usepackage{stmaryrd}  
\usepackage[all,2cell,cmtip]{xy} \UseAllTwocells \SilentMatrices
\usepackage{tikz-cd}
\usepackage[dvips]{epsfig}
\usepackage{MnSymbol}
\usepackage{comment}
\usepackage{enumitem}
\usepackage{kbordermatrix}
\usepackage[colorinlistoftodos]{todonotes}

\renewcommand{\kbldelim}{(}
\renewcommand{\kbrdelim}{)}

\usetikzlibrary{arrows,automata}
\usetikzlibrary{decorations.markings}
\usetikzlibrary{decorations.pathmorphing}

\makeatletter
\def\@adminfootnotes{%
  \let\@makefnmark\relax  \let\@thefnmark\relax
  \ifx\@empty\@date\else \@footnotetext{\@setdate}\fi
  \ifx\@empty\@subjclass\else \@footnotetext{\@setsubjclass}\fi
  \ifx\@empty\@keywords\else \@footnotetext{\@setkeywords}\fi
  \ifx\@empty\thankses\else \@footnotetext{%
    \def\par{\let\par\@par}\@setthanks}%
  \fi
}
\makeatother

\newtheorem{theorem}{Theorem}[section]

\theoremstyle{definition}

\newtheorem{conjecture}[theorem]{Conjecture}

\newtheorem{example}[theorem]{Example}

\theoremstyle{remark}

\numberwithin{equation}{section}

\let\oldtocsection=\tocsection
\let\oldtocsubsection=\tocsubsection
\renewcommand{\tocsection}[2]{\hspace{0em}\oldtocsection{#1}{#2}}
\renewcommand{\tocsubsection}[2]{\hspace{1em}\oldtocsubsection{#1}{#2}}

\usepackage{chngcntr}
\counterwithin{figure}{subsection}

\makeatletter
\@namedef{subjclassname@2020}{%
  \textup{2020} Mathematics Subject Classification}

\makeatother
 
\title[From automata to tangle cobordisms]{From finite state automata to tangle cobordisms: a TQFT journey from 
one to four dimensions}
 
\author{Mee Seong Im}
\address{Department of Mathematics, United States Naval Academy, Annapolis, MD 21402, USA}
\email{\href{mailto:meeseongim@gmail.com}{meeseongim@gmail.com}}
\thanks{}

\author{Mikhail Khovanov} 
\address{Department of Mathematics, Johns Hopkins University, Baltimore, MD 21218, USA}
\email{\href{mailto:khovanov@jhu.edu}{khovanov@jhu.edu}}
\address{Department of Mathematics, Columbia University, New York, NY 10027, USA}
\email{\href{mailto:khovanov@math.columbia.edu}{khovanov@math.columbia.edu}}
\thanks{}
 
\subjclass[2020]{Primary: 
57K16, 
57K18, 
57K45, 
68Q45, 
18M30; 
Secondary: 
68Q70, 
18B20} 
\date{September 17, 2023}

\dedicatory{To Igor B. Frenkel, on his anniversary, with admiration.}

\providecommand{\keywords}[1]{\textbf{\textit{Key words and phrases.}} #1}

\keywords{Topological quantum field theory, TQFT, finite-state automata,  Reshetikhin--Turaev invariants, categorification,  link homology, Khovanov--Rozansky homology, foam evaluation. }

\begin{document}
\def\CMF{\mathsf{CMF}}
\def\Inv{\mathsf{Inv}}
\def\Tan{\mathsf{Tan}}
\def\mcF{\mathcal{F}}

\def\Aff{\mathsf{Aff}}
\def\AND{\mathsf{AND}}
\def\concatenate{\mathsf{concatenate}}
\def\Br{\mathsf{Br}}
\def\Gal{\mathsf{Gal}}
\def\gen{\mathsf{generators}}
\def\GL{\mathsf{GL}}
\def\SL{\mathsf{SL}}
\def\init{\mathsf{in}}
\def\t{\mathsf{t}}
\def\out{\mathsf{out}}
\def\I{\mathsf I}
\def\region{\mathsf{region}}
\def\plane{\mathsf{plane}}
\def\R{\mathbb R}
\def\Q{\mathbb Q}
\def\Z{\mathbb Z}
\def\mc{\mathcal{c}}
\def\finite{\mathsf{finite}}
\def\infinite{\mathsf{infinite}}
\def\N{\mathbb N} 
\def\C{\mathbb C}
\def\sep{\mathsf{sep}}
\def\S{\mathbb S}
\def\SS{\mathbb S} 
\def\CP{\mathbb P}
\def\Ob{\mathsf{Ob}}
\def\op{\mathsf{op}}
\def\new{\mathsf{new}}
\def\old{\mathsf{old}}
\def\OR{\mathsf{OR}}
\def\AND{\mathsf{AND}}
\def\rat{\mathsf{rat}}
\def\rec{\mathsf{rec}}
\def\tail{\mathsf{tail}}
\def\coev{\mathsf{coev}}
\def\ev{\mathsf{ev}}
\def\id{\mathsf{id}}
\def\s{\mathsf{s}}
\def\S{\mathsf{S}}
\def\t{\mathsf{t}}
\def\start{\textsf{starting}}
\def\Notation{\textsf{Notation}}
\def\circleft{\raisebox{-.18ex}{\scalebox{1}[2.25]{\rotatebox[origin=c]{180}{$\curvearrowright$}}}}
\renewcommand\SS{\ensuremath{\mathbb{S}}}
\newcommand{\kllS}{\kk\llangle  S \rrangle} 
\newcommand{\kllSS}[1]{\kk\llangle  #1 \rrangle}
\newcommand{\klS}{\kk\langle S\rangle}  
\newcommand{\aver}{\mathsf{av}}  
\newcommand{\ophana}{\overline{\phantom{a}}}
\newcommand{\Bool}{\mathbb{B}}
\newcommand{\dmod}{\mathsf{-mod}}
\newcommand{\lang}{\mathsf{lang}}
\newcommand{\pfmod}{\mathsf{-pfmod}}
\newcommand{\primitive}{\mathsf{irr}}
\newcommand{\Bmod}{\Bool\mathsf{-mod}}  
\newcommand{\Bmodo}[1]{\Bool_{#1}\mathsf{-mod}}  
\newcommand{\Bfmod}{\Bool\mathsf{-fmod}} 
\newcommand{\Bfpmod}{\Bool\mathsf{-fpmod}} 
\newcommand{\Bfsmod}{\Bool\mathsf{-}\underline{\mathsf{fmod}}}  
\newcommand{\undvar}{\underline{\varepsilon}} 
\newcommand{\RLang}{\mathsf{RLang}}
\newcommand{\undotimes}{\underline{\otimes}}
\newcommand{\sigmaacirc}{\Sigma^{\ast}_{\circ}} 
\newcommand{\cl}{\mathsf{cl}}
\newcommand{\PP}{\mathcal{P}} 
\newcommand{\wedgezero}{\{ \vee ,0\} } 
\newcommand{\whA}{\widehat{A}}
\newcommand{\whC}{\widehat{C}}
\newcommand{\whM}{\widehat{M}}
\newcommand{\Sigmalr}{\Sigma^{\Z}}
\newcommand{\Sigmal}{\Sigma^{-}}
\newcommand{\Sigmar}{\Sigma^{+}}
\newcommand{\Sigmaa}{\Sigma^{\ast}}
\newcommand{\SigmaZ}{\Sigma^{\Z}}  
\newcommand{\Sigmac}{\Sigma^{\circ}}

\newcommand{\alphai}{\alpha_I}  
\newcommand{\alphac}{\alpha_{\circ}}  
\newcommand{\alphap}{(\alphai,\alphac)} 
\newcommand{\alphalr}{\alpha_{\leftrightarrow}}
\newcommand{\alphaZ}{\alpha_{\Z}}
\newcommand{\mcCinfalpha}{\mcC^{\infty}_{\alpha}}
\newcommand{\mathT}{\mathsf{T}}
\newcommand{\mathF}{\mathsf{F}}
\newcommand{\mcS}{\mathcal{S}}
\newcommand{\mcN}{\mathcal{N}}
\newcommand{\wmcN}{\widetilde{\mcN}}
\newcommand{\Net}{\mathsf{Net}}

\let\oldemptyset\emptyset
\let\emptyset\varnothing

\newcommand{\undempty}{\underline{\emptyset}}
\newcommand{\undsigma}{\underline{\sigma}}
\newcommand{\undtau}{\underline{\tau}}
\def\basis{\mathsf{basis}}
\def\irr{\mathsf{irr}} 
\def\spanning{\mathsf{spanning}}
\def\elmt{\mathsf{elmt}}

\def\H{\mathsf{H}}
\def\I{\mathsf{I}}
\def\II{\mathsf{II}}
\def\l{\lbrace}
\def\r{\rbrace}
\def\o{\otimes}
\def\lra{\longrightarrow}
\def\Ext{\mathsf{Ext}}
\def\mf{\mathfrak} 
\def\mcC{\mathcal{C}}
\def\mcO{\mathcal{O}}
\def\Fr{\mathsf{Fr}}

\def\ovb{\overline{b}}
\def\tr{{\sf tr}} 
\def\str{{\sf str}} 
\def\det{{\sf det }} 
\def\tral{\tr_{\alpha}}
\def\one{\mathbf{1}}   

\def\lra{\longrightarrow}
\def\twoheadlra{\longrightarrow\hspace{-4.6mm}\longrightarrow}
\def\hooklra{\raisebox{.2ex}{$\subset$}\!\!\!\raisebox{-0.21ex}{$\longrightarrow$}}
\def\kk{\mathbf{k}}  
\def\gdim{\mathsf{gdim}}  
\def\rk{\mathsf{rk}}
\def\undep{\underline{\varepsilon}}
\def\mathM{\mathbf{M}}  

\def\CCC{\mathcal{C}} 
\def\wCCC{\widehat{\CCC}}  

\def\complement{\mathsf{comp}}
\def\Rec{\mathsf{Rec}} 

\def\Cob{\mathsf{Cob}} 
\def\Kar{\mathsf{Kar}}   

\def\dmod{\mathsf{-mod}}   
\def\pmod{\mathsf{-pmod}}    

\newcommand{\brak}[1]{\ensuremath{\left\langle #1\right\rangle}}
\newcommand{\brakspace}[1]{\ensuremath{\left\langle \:\: #1\right\rangle}}

\newcommand{\oplusop}[1]{{\mathop{\oplus}\limits_{#1}}}
\newcommand{\ang}[1]{\langle #1 \rangle } 
\newcommand{\ppartial}[1]{\frac{\partial}{\partial #1}} 
\newcommand{\checkr}{{\bf \color{red} CHECK IT}}
\newcommand{\checkb}{{\bf \color{blue} CHECK IT}}
\newcommand{\checkk}[1]{{\bf \color{red} #1}}

\newcommand{\mcA}{{\mathcal A}}
\newcommand{\cZ}{{\mathcal Z}}
\newcommand{\sq}{$\square$}
\newcommand{\bi}{\bar \imath}
\newcommand{\bj}{\bar \jmath}

\newcommand{\undn}{\underline{n}}
\newcommand{\undm}{\underline{m}}
\newcommand{\undzero}{\underline{0}}
\newcommand{\undone}{\underline{1}}
\newcommand{\undtwo}{\underline{2}}

\newcommand{\cob}{\mathsf{cob}} 
\newcommand{\comp}{\mathsf{comp}} 

\newcommand{\Aut}{\mathsf{Aut}}
\newcommand{\Hom}{\mathsf{Hom}}
\newcommand{\Idem}{\mathsf{Idem}}
\newcommand{\Ind}{\mbox{Ind}}
\newcommand{\Id}{\textsf{Id}}
\newcommand{\End}{\mathsf{End}}
\newcommand{\iHom}{\underline{\mathsf{Hom}}}
\newcommand{\Bools}{\Bool^{\mathfrak{s}}}
\newcommand{\mfs}{\mathfrak{s}}
\newcommand{\blueline}{line width = 0.45mm, blue}

\newcommand{\drawing}[1]{
\begin{center}{\psfig{figure=fig/#1}}\end{center}}

\def\endomCempt{\End_{\mcC}(\emptyset_{n-1})}

\def\MS#1{{\color{blue}[MS: #1]}}
\def\MK#1{{\color{red}[MK: #1]}}

\begin{abstract} 
This is a brief introduction to link homology theories that categorify Reshetikhin--Turaev $\SL(N)$-quantum link invariants. A recently discovered surprising connection between finite state automata and Boolean TQFTs in dimension one is explained as a warm-up.
\end{abstract}

\maketitle
\tableofcontents

%
%

\section{Introduction}
\label{section:intro}

This note is based on the talk that one of us (M.K.) gave at the First International Congress of Basic Science which was held at the Yanqi Lake Beijing Institute of Mathematical Sciences and Applications (BIMSA) in July 2023. 

We briefly discuss topological quantum field theories (TQFTs) and explain a recent surprising observation~\cite{GIKKL23} that Boolean-valued one-dimensional TQFTs with defects correspond to nondeterministic finite state automata. The precise statement is given in Theorem~\ref{thm_automata_TQFT} below and in~\cite{GIKKL23}. 

We then review various approaches to categorification of Reshetikhin--Turaev link invariants for the fundamental representations of quantum $\SL(N)$ emphasizing the original approach via matrix factorizations~\cite{KR04,KR05,Yon11,WuEquiv12} and Robert--Wagner foam evaluation approach~\cite{RW20}. 

Reshetikhin--Turaev link invariants \cite{RT90,Resh91} are part of the Chern--Simons three-dimensional TQFT discovered by Witten and Reshetikhin--Turaev~\cite{Witten89,RT91}. Categorification of Reshetikhin--Turaev link invariants can be viewed as a 4-dimensional TQFT restricted to links in $\R^3$ and to link cobordisms. Boolean one-dimensional TQFTs and finite state automata connection is explained first, as a warm-up to discussing these  sophisticated TQFTs in dimensions three and four. 

\vspace{0.05in} 

This is a relatively short paper which reviews Khovanov--Rozansky link homology in its second part. We provide key references but can not give a fully comprehensive coverage of the substantial amount of literature on this and closely related subjects due to the immense and rapid development over the last twenty-five years.

\vspace{0.1in} 

{\bf Acknowledgments:} The authors are grateful to Tsinghua University and Yau Mathematical Sciences Center for their hospitality. 
M.K. would like to acknowledge support from Beijing Institute of Mathematical Sciences and Applications (BIMSA) and the opportunity to give a plenary lecture at the first International Congress of Basic Science (ICBS).  The authors would like to thank Joshua Sussan for valuable feedback on the paper.   M.K. was partially supported by NSF grant DMS-2204033 and Simons Collaboration Award 994328 while working on this paper.

\section{TQFTs, Boolean TQFTs and automata}
\label{sec_boolean}

In this section, we survey the correspondence between automata and one-dimensional Boolean TQFTs explained in~\cite{GIKKL23}. 

\subsection{TQFTs} 
An $n$-dimensional TQFT (topological quantum field theory) over a field $\kk$ is a tensor functor from the category of $n$-dimensional cobordisms to the category of vector spaces 
\[
\mcF \ :  \ \mathsf{Cob}_n  \lra \mathbf{k}\mathsf{-vect}. 
\]
The category $\mathsf{Cob}_n$ has oriented $(n-1)$-manifolds as objects and $n$-dimensional cobordisms between them as morphisms. This category is symmetric tensor, and the functor $\mcF$ must respect this structure so that, in particular, 
\[
\mcF(K_1\sqcup K_2) \cong \mcF(K_1)\otimes \mcF(K_2)
\]
for $(n-1)$-manifolds $K_1$ and $K_2$. 


\subsection{One-dimensional TQFTs} 
A TQFT $\mcF$ in one dimension assigns a vector space $V$ to a point with positive orientation, a space $W$ to a point with negative orientation, which we can write as
\[ 
\mcF(+)= V, \hspace{1cm} \mcF(-)=W, 
\]
and maps $\kk\lra V\otimes W$ and $W\otimes V\lra \kk$  to the \emph{cup}, respectively, $\emph{cap}$, cobordisms shown in the top row of Figure~\ref{figure-0.1}. Transposition cobordisms, see Figure~\ref{figure-0.1} top right, are mapped to transposition maps between products of $V$ and $W$, such as $V\otimes W\lra W\otimes V$, $v\otimes w\mapsto w\otimes v$. 

Isotopy relations shown in the second row of that figure imply that vector space $V$ is finite-dimensional and vector space $W\cong V^{\ast}$ can be taken to be the dual of $V$, so that cup and cap cobordisms are sent by $\mcF$  to the \emph{evaluation} and \emph{coevaluation} maps: 
\[ 
\mathsf{coev}: \kk \lra V \otimes V^{\ast}, \ \ \ \ \mathsf{ev}: V^{\ast}\otimes V \lra \kk \]
given in a basis $\{v_1,\ldots, v_n\}$ of $V$ and the dual basis $\{v^1,\ldots, v^n\}$ of $V^{\ast}$ by 
\begin{eqnarray*}
\mathsf{coev}(1) & = & \sum_{i=1}^n v_i\otimes v^i, \\
\mathsf{ev}(v^i\otimes v_j) & = & \delta_{i,j}, \ \  \mathsf{ev}(v^{\ast}\otimes v)=v^{\ast}(v), \ \  v^{\ast}\in V^{\ast}, \ v \in V. 
\end{eqnarray*}
Composing with the transposition morphism in Figure~\ref{figure-0.1} in top right gives the cup and cap morphisms for the opposite orientation of the arc. 

\input{figure-0.1}

The TQFT $\mcF$ evaluated on a circle gives a linear map $\kk\lra \kk$ which is multiplication by $\dim V$, and that is the only invariant of a one-dimensional TQFT. Two such TQFTs 
\[
\mcF \ :  \ \mathsf{Cob}_1  \lra \mathbf{k}\mathsf{-vect}, \hspace{1cm} \mcF(+)=V, \hspace{1cm} \mcF(-)=V^{\ast} 
\]
are isomorphic if and only if the corresponding vector spaces are isomorphic, that is, if they have the same dimension (evaluate the same on the circle).

\vspace{0.07in} 

To spice up this simple classification, let us add $0$-dimensional defects, which are points on one-manifolds, labelled by elements of a finite set $\Sigma$, and also allow one-manifolds to end in the middle of a cobordism, see Figure~\ref{figure-A1}. In this way, the boundary of a one-manifold cobordism splits into the \emph{outer boundary} and \emph{inner boundary}. Connected components of a morphism in $\Cob_{\Sigma,1}$ can be classified into four types:
\begin{itemize}
    \item arcs with two outer endpoints, 
    \item \emph{half-intervals}, which are arcs with one outer and one inner endpoint, 
    \item \emph{floating intervals}, which are arcs with two inner endpoints, 
    \item circles. 
\end{itemize}
The morphism in Figure~\ref{figure-A1} consists of two arcs, three half-intervals, one floating interval and one circle. Each connected component of a morphism may be decorated by $\Sigma$-labelled dots. 

\vspace{0.07in} 

\input{figure-A1}

This bigger category $\mathsf{Cob}_{\Sigma,1}$ is also a  symmetric tensor category, containing $\mathsf{Cob}_1$ as a subcategory. How can one classify TQFTs 
(tensor functors) 
\[
\mcF \ :  \ \mathsf{Cob}_{\Sigma,1}  \lra \mathbf{k}\mathsf{-vect} \ ? 
\]
Necessarily $\mcF(+)\cong V$ and $\mcF(-)\cong V^{\ast}$ for some finite-dimensional vector space $V$, with cup and cap maps given by evaluation and coevaluation. 

Applying $\mcF$ to a dot labelled $a$ on an upward-oriented strand gives a morphism $\mcF(+)\lra \mcF(+)$, that is, a linear map $m_a:V\lra V$, see Figure~\ref{figure-3}. Applying $\mcF$ to a dot labelled $a$ on a downward-oriented strand then necessarily gives the dual linear map $m_a^{\ast}:V^{\ast}\lra V^{\ast}$ (to check the duality property, move a dot across a local maximum or minimum of a cobordism shown in Figure~\ref{figure-0.1} top left). A collection of labelled dots on a strand can be encoded by a word $\omega=a_1\cdots a_n$, with the corresponding map on $V$ given by the product, see Figure~\ref{figure-3}.  

\input{figure-3}

\input{figure-2}

A \emph{half-interval} (an interval with one outer and one inner endpoint) oriented upward with the outer endpoint at the top, upon applying $\mcF$,  gives a linear map $\kk\cong \mcF(\emptyset)\lra \mcF(+)\cong V$, that is, a vector $v_0\in V$, see Figure~\ref{figure-2}. For the other half-interval with a $+$ endpoint, also see Figure~\ref{figure-2}, applying $\mcF$ results in a morphism $V\lra \kk$,  described by a covector $v^{\ast}\in V^{\ast}$. 

No relations on $m_a$'s, $v_0$ and $v^{\ast}$ are imposed. Thus, a TQFT for the category $\mathsf{Cob}_{\Sigma,1}$ is given by a finite-dimensional vector space $V$, a collection of linear maps $m_a:V\lra V$, for $a\in \Sigma$, a vector $v_0$ and a covector $v^{\ast}$. Two such theories  are isomorphic if there is an isomorphism between the corresponding vector spaces that intertwines maps $m_a$, for all $a$, vectors $v_0$ and covectors $v^{\ast}$ for the two theories. 

Topological quantum field theory $\mcF$ has a large number of invariants, which are values of $\mcF$ on: 
\begin{itemize}
\item 
A circle carrying a word $\omega=a_1\cdots a_n$, up to cyclic permutation. The functor $\mcF$ evaluates this decorated circle to the trace $\tr_V(m_{a_n}\cdots m_{a_1})\in \kk$. 
\item A floating interval carrying a word $\omega =a_1\cdots a_n$. The functor $\mcF $ evaluates this interval to $v^{\ast}(m_{a_1}\cdots m_{a_n}v_0)\in \kk$. 
\end{itemize} 
These evaluations are depicted in Figures~\ref{figure-4} and~\ref{figure-5}. For instance, for an interval, start with $v_0\in V$ at the endpoint where the orientation looks into the interval, then apply the map $m_{a_n}$  to $v_0$ and so on until the opposite end of the interval is reached. Then we apply the covector $v^{\ast}$ to the product. 

\input{figure-4}

\input{figure-5}

\vspace{0.1in} 

\subsection{Automata and regular languages}
\label{subsection:automata_regular_lang}
By an \emph{alphabet} or \emph{set of letters} we mean a finite set $\Sigma$. A \emph{language} $L$ is any subset of the free monoid $\Sigma^{\ast}$ on $\Sigma$. 

A (nondeterministic) finite state automaton $(Q)$ over alphabet $\Sigma$ consists of a finite set of states $Q$, a \emph{transition function} 
\begin{equation}\label{eq_delta}
\delta: \Sigma\times Q\lra \mathcal{P}(Q),
\end{equation}
where $\mathcal{P}(Q)$ is the powerset of $Q$, a subset $Q_{\init}\subset Q$ of \emph{initial} states and a subset $Q_{\t}\subset Q$ of \emph{accepting} or \emph{terminal} states. 

To an automaton there is associated a graph $\Gamma_{(Q)}$ with $Q$ as the set of vertices, an oriented edge $q\lra q'$ labelled $a\in \Sigma$ if and only if $q'\in \delta(a,q)\subset Q$ and two subsets $Q_{\init}, Q_{\t}$ of distinguished vertices. Vice versa, such a data of a decorated oriented graph determines an automaton. 

A word $\omega=a_1\cdots a_n\in \Sigma^{\ast}$ is \emph{accepted} by the automaton $(Q)$ if there exists an initial state $q_{\init}\in Q_{\init}$, an accepting state $q_{\t}\in Q_{\t}$ and an oriented path from $q_{\init}$ to $q_{\t}$ where consecutive labels of oriented edges are $a_1,\ldots, a_n$, see Figure~\ref{fig_A1}. 

\vspace{0.1in} 

\input{fig_A1}

\vspace{0.1in} 

The set of words $\omega\in \Sigma^{\ast}$ accepted by $(Q)$ is called \emph{the language of automaton $(Q)$}. We denote this language by $L_{(Q)}$. A language $L\subset \Sigma^{\ast}$ is called \emph{regular} if it is the language of some automaton. One can check that a language is regular if and only if it can be described by a regular expression.  

\begin{example}\label{ex_language} Consider the language $L$ for the alphabet $\Sigma=\{a,b\}$ given by the regular expression $L=(a+b)^{\ast}b(a+b)$. This expression describes all words which have $b$ as the second from the last letter. An example of a nondeterministic automaton for $L$ is shown in Figure~\ref{figure-X1}. 

\input{figure-X1}

\end{example}

\subsection{Boolean TQFT from an automaton}
\label{subsection:boolean_TQFT}
In the definition of a TQFT, we can replace field $\kk$ by any commutative semiring. Commutativity is needed for the following reason: floating components of a cobordism are evaluated to elements of a ground ring or semiring. These components can change relative position, that is, float up or down past each other, which correspond to requiring that the ground (semi)ring is \emph{commutative}. 

\vspace{0.07in} 

By analogy with the category $\Cob_{\Sigma,1}$ and a TQFT $\mcF$ on it we can imagine encoding words $\omega$ in some language $L$ by placing their labels next to dots along a one-manifold. A floating interval that carries a word $\omega$ can evaluate to two values, depending on whether or not $\omega$ is in the language $L$. We make these values $0,1$, with $\omega$ evaluating to $1$ if it belongs to a language $L$ and to $0$ otherwise. It is then natural to replace a field $\kk$ by the Boolean semiring $\Bool=\{0,1 : 1+1=1\}$ which consists of these two elements, with the additional rule that $1+1=1$. In particular, $\Bool$ is not a ring, but a commutative semiring. 

We specialize to the commutative semiring $\Bool$ and fix an automaton $(Q)$ as above. 
Replace a $\kk$-vector space $V\cong \kk^n$ by a free module $\Bool Q$ over $\Bool$ with a set of generators $Q$. Elements of $\Bool Q$ are formal finite sums of distinct elements of $Q$ or, equivalently, subsets of $Q$, so we can identify $\Bool Q\cong \mathcal{P}(Q)$ with the set of subsets of $Q$. The zero element of $\Bool Q$ corresponds to the empty subset of $Q$. Note that $q+q=q$ for $q\in Q$ since $1+1=1$. We view $\Bool Q$ as a free $\Bool$-module with a basis $Q$. It has $2^{|Q|}$ elements. Unlike the case of a field, where the group $\mathsf{GL}(V)$ acts freely and transitively on the set of bases of $V$, the only basis of $\Bool Q$, up to changing the order of elements, is $Q$.

To the automaton $(Q)$ we associate a Boolean-valued TQFT 
\begin{equation}
    \mcF \ : \ \Cob_{\Sigma,1} \lra \Bool\mathsf{-fmod} 
\end{equation}
taking values in the category of free $\Bool$-modules by assigning $\Bool Q$ to the $+$ point and its dual $\Bool Q^*$ to the $-$ point: 
\begin{align*}
    &\mcF(+) = \Bool Q, \hspace{3.0cm} \mcF(-) = \Bool Q^{\ast}, \\
    &\mathsf{ev} : \ \Bool Q^{\ast} \otimes \Bool Q\lra \Bool, \hspace{1.35cm} \mathsf{ev}(q^{\ast}\otimes q') = \delta_{q,q'} \\
    &\mathsf{coev}: \ \Bool \lra \Bool Q \otimes \Bool Q^{\ast}, \hspace{1cm} \mathsf{coev}(1) = \sum_{q\in Q}q\otimes q^{\ast}. 
\end{align*}
Here $Q^{\ast}$ is a copy of the set $Q$ with elements labelled $q^{\ast}$, over $ q\in Q$, with the evaluation and coevaluation maps $\mathsf{ev}$ and $\mathsf{coev}$ given by the above formulas. 

For each $a\in \Sigma$ the transition map $\delta_a : Q\lra \mathcal{P}(Q)$ extends to a $\Bool$-linear endomorphism $\delta_a$ of the free module $\Bool Q$ that takes $q$ to $\delta(a,q)$, using the identification $\mathcal{P}(Q)\cong \Bool Q$, see above. Note that $\delta_a$ denotes the transition map while $\delta_{q,q'}$ stands for the delta function on a pair of states. 

For example, in the basis $\{q_1,q_2,q_3\}$ of $\Bool Q$ in Example~\ref{ex_language} the two maps for $a,b\in \Sigma$ are given by $\Bool$-valued matrices 
\[
\delta_a = \begin{pmatrix}
1 & 0 & 1\\
0 & 0 & 0\\
0 & 1 & 0 
\end{pmatrix}, 
\hspace{1cm}
\delta_b = \begin{pmatrix}
0 & 0 & 0\\
1 & 1 & 1\\
0 & 1 & 0 
\end{pmatrix},
\]
see Figure~\ref{figure-X1}. 

To an upward-oriented arc cobordism from $+$ to $+$ with an $a$-dot, $a\in \Sigma$, we associate the endomorphism $\delta_a$ of $\Bool Q$, see Figure~\ref{figure-X2}. 

\input{figure-X2}

To a downward-oriented arc, which is a cobordism from $-$ to $-$, we associate the dual map $\delta_a^{\ast}:\Bool Q^{\ast}\lra \Bool Q^{\ast}$ given by the transposed matrix to that of $\delta_a$, see Figure~\ref{figure-X2}.

To upward-oriented half-intervals (cobordisms from the empty 0-manifold $\emptyset_0$ to $+$ and from $+$ to $\emptyset_0$) associate maps between free $\Bool$-modules $\Bool$ and $\Bool Q$ given by 
\begin{eqnarray*}
 & \Bool \lra \Bool Q, \ \ 1\longmapsto \displaystyle{\sum_{q\in Q_{\init}}} q, \\
 & \Bool Q\lra \Bool, \ \ q\longmapsto 
 \begin{cases}
     1 & \mathsf{if }\ q\in Q_{\t}, \\
     0 & \mathsf{if }\  q\notin Q_{\t}, 
 \end{cases}
\end{eqnarray*} 
see Figure~\ref{figure-1-1}, and likewise for downward-oriented half-intervals.
These maps are determined by subsets $Q_{\init}$ of initial and $Q_{\t}$ of accepting states. Alternatively the image of $1$ under the first map can be written as $Q_{\init}$, using identification $\Bool Q\cong \mathcal{P}(Q)$, and the second map can be denoted $Q_{\t}^{\ast}$.

\input{figure-1-1}

Given a floating interval with the word $\omega=a_1\cdots a_n$ written on it so that $a_1$ is at the tail and $a_n$ at the head, the interval evaluates to $1\in \Bool$ if and only if $\omega$ is in the language $L_{(Q)}$ defined by the automaton $(Q)$, see Figure~\ref{figure-5a}. If $\omega$ is not in $L_{(Q)}$ the interval evaluates to $0$. The evaluation can also be written as $Q_{\t}^{\ast}(\delta_{a_n}\cdots \delta_{a_1}Q_{\init})$, that is, evaluating the product $\delta_{a_n}\cdots \delta_{a_1}Q_{\init}\in \Bool Q$ on the functional $Q_{\t}^{\ast}$.  

\input{figure-5a}

A circle with a circular word $\omega=a_1\cdots a_n$ placed on it evaluates to the trace of operator $\delta_{\omega}:=\delta_{a_n}\cdots \delta_{a_1}$ on $\Bool Q$, see Figure~\ref{figure-4a}.  Equivalently, it evaluates to $1$ if and only if there is a state $q\in Q$ and a sequence of arrows labelled $a_1,a_2,\ldots, a_n$ terminating back in $q$ (\textit{i.e.}, and oriented cycle with word $\omega$ in the graph of the automaton). 

\input{figure-4a}

We state a main result in \cite{GIKKL23}:

\begin{theorem}[Gustafson--Im--Kaldawy--Khovanov--Lihn]
\label{thm_automata_TQFT} A nondeterministic finite state automaton $(Q)$ on alphabet $\Sigma$ defines a one-dimensional Boolean-valued TQFT 
\begin{equation}\label{eq_thm_TQFT}
    \mcF \ : \ \Cob_{\Sigma,1} \lra \Bool\mathsf{-fmod} 
\end{equation}
with $\mcF(+)=\Bool Q$, transition function of $(Q)$ encoding TQFT maps for $\Sigma$-labelled defects on strands, and sets of initial and accepting states encoding the maps for undecorated half-intervals. This correspondence gives a bijection between isomorphism classes of Boolean-valued TQFTs for $\Cob_{\Sigma,1}$ with $\mcF(+)$ a free $\Bool$-module and isomorphism classes of nondeterministic finite state automata on alphabet $\Sigma$. 
\end{theorem} 
$\Bool\mathsf{-fmod}$ in \eqref{eq_thm_TQFT} denotes the category of free $\Bool$-semimodules and semimodule maps. 

Regular language $L_{(Q)}$ of the automaton $(Q)$ describes evaluation of floating intervals decorated by words in $\Sigma^{\ast}$ (with $\omega$-decorated interval evaluating to $1$ if and only if $\omega \in L_{(Q)}$), while evaluation of circles with defects in the TQFT $\mcF$ is determined by oriented cycles in the graph of $(Q)$.

The table in Figure~\ref{mfig_022} summarizes the correspondence between generating morphisms in $\Cob_{\Sigma,1}$ and structural parts of an automaton (the set of states, transition maps, and sets of initial and accepting states). 

\input{mfig_022}

Consider a one-dimensional TQFT 
\begin{equation}\label{eq_thm_TQFT_2}
    \mcF \ : \ \Cob_{\Sigma,1} \lra \Bool\mathsf{-mod} 
\end{equation}
valued, more generally, in the category of all $\Bool$-modules rather than free ones. Then necessarily $P=\mcF(+)$ is a finitely-generated semimodule which is \emph{projective} in the sense of being a retract of a free module~\cite{GIKKL23,IK-top-automata}. Namely, there are semimodule maps 
\begin{equation}\label{eq_retract}
  P \stackrel{\iota}{\lra} \Bool^n \stackrel{p}{\lra} P, \hspace{1cm} 
  p\circ\iota=\id_P
\end{equation}
for some $n$. Note that $P$ is a direct summand of $\Bool^n$ only if $P$ is free, otherwise it is just a retract. Then $P^{\ast}\cong \mcF(-)$ is finitely-generated projective as well, with the retract maps $p^{\ast},\iota^{\ast}$ given by dualizing those for $P$. 

Finitely-generated projective $\Bool$-modules $P$ are described by finite distributive lattices, \textit{i.e.}, see \cite{IK-top-automata}: it is a theorem that goes back at least to Birkhoff that any such $P$ is isomorphic to the $\Bool$-semimodule $\mathcal{U}(X)$ of open sets in a finite topological space $X$, with the empty set $\emptyset$ the $0$ element of the semimodule and addition $U_1+U_2:=U_1 \cup U_2$ given by the union of sets. 

The structure of a TQFT $\mcF$ in this case is given by 
\begin{itemize}
\item 
a collection of endomorphisms $\delta_a:\mathcal{U}(X)\lra \mathcal{U}(X)$ for $a\in \Sigma$, taking open sets to open sets, $\emptyset$ to $\emptyset$, and  preserving unions of sets, 
\item
initial element $Q_{\init}\in \mathcal{U}(X)$ and a terminal map $Q_{\t}:\mathcal{U}(X)\lra \Bool$ taking $\emptyset$ to $0$ and intertwining union of sets with addition in $\Bool$. 
\end{itemize}
Such structures are called \emph{quasi-automata} in~\cite{GIKKL23}. It is an interesting question whether they can be of use in computer science.

\subsection{Extending to arbitrary commutative semirings}
It is straightforward to extend the above TQFT construction from vector spaces over a field $\kk$ and $\Bool$-semimodules to semimodules over a commutative semiring $R$. 
 A commutative semiring $R$ is an abelian group under addition, a commutative monoid under multiplication, and distributivity property holds, $a(b+c)=ab+ac$. Semiring $R$ has the zero element $0$ and the unit element $1$. Subtraction operation $a-b$ is usually not available in semirings. It is straightforward to define the notion of a module $M$ over $R$ (alternatively called a \emph{semimodule}) and introduce the category $R\mathsf{-dmod}$ of $R$-modules.  

In the definition of a TQFT over $R$ the 
tensor product of vector spaces is replaced by the tensor product of $R$-modules, and vector space $V=\mcF(+)$ is replaced by an $R$-module $P=\mcF(+)$. To have cup and cap morphisms subject to the isotopy relations above requires $P$ to be a projective $R$-module of finite rank, see~\cite{GIKKL23}, for instance.  

This observation quickly leads to the following in \cite{GIKKL23}: 
\begin{theorem}[Gustafson--Im--Kaldawy--Khovanov--Lihn]
\label{thm_R}
    One-dimensional TQFTs 
    \[ 
    \mcF:\Cob_1\lra R\mathsf{-mod}
    \] 
    over a commutative semiring $R$ correspond to finitely-generated projective $R$-modules. 
    One-dimensional TQFTs with defects, that is, tensor functors 
    \begin{equation} \label{TQFT_S_R}
    \mcF \ : \ \mathsf{Cob}_{\Sigma,1}\lra R\mathsf{-mod} 
    \end{equation}
    are in a correspondence with finitely-generated projective $R$-modules $P$ equipped with endomorphisms $m_a:P\lra P$ for $a\in \Sigma$, an element $v_0\in P$ and a covector $v^{\ast}\in \Hom_R(P,R)$. This correspondence is a bijection between isomorphism classes of TQFTs and isomorphism classes of data $(P,\{m_a\}_{a\in\Sigma},v_0,v^{\ast})$. 
\end{theorem}
Here, an $R$-semimodule $P$ is defined to be finitely-generated projective if it is a retract of free semimodule $R^n$ for some $n$: 
\begin{equation}\label{eq_retract_2}
  P \stackrel{\iota}{\lra} R^n \stackrel{p}{\lra} P, \hspace{1cm} 
  p\circ\iota=\id_P . 
\end{equation}

It may be interesting to look at examples of such TQFTs when $R$ is, for instance, the tropical semiring, see~\cite{JK_tropical14,IJK_tropical_poly16} where projective modules over the tropical semiring are studied.  

Another question is whether the notion of a \emph{finite state machine}, which extends the notion of a finite state automaton, has the TQFT counterpart. 

The authors are not aware of any studies or results on  Boolean TQFTs (and, generally, TQFTs over commutative semirings that are not rings) in dimension two and higher. The above correspondence between finite state automata and one-dimensional Boolean TQFTs with defects, observed in~\cite{GIKKL23}, see also~\cite{IK-top-automata}, remains a curiosity, for now. 

Two-dimensional TQFTs for oriented surfaces, without defects and over a field $\kk$, are classified by commutative Frobenius $\kk$-algebras, and one open problem is to find a supply of commutative Frobenius $\kk$-semialgebras when $\kk$ is the Boolean semiring $\Bool$ or the tropical semiring.  

%
%
 
\section{Reshetikhin--Turaev \texorpdfstring{$\mathsf{SL}(N)$}{SLN}-invariants and their categorification}

Quantum link polynomials were discovered by Vaughan Jones~\cite{Jones85,Jones97} (the Jones polynomial), 
Louis Kauffman~\cite{Kauffman90} (the Kauffman polynomial and bracket), J.~Hoste, A.~Ocneanu, K.~Millet, P.~Freyd, W.B.R.~Lickorish, D.~Yetter~\cite{FYHLMO85}, J.~Przytycki and P.~Traczyk~\cite{Przytycki88} (the HOMFLYPT polynomial) 
 and others (Alexander polynomial, discovered several decades prior, was an outlier). N.~Reshetikhin and V.~Turaev~\cite{RT90} put these polynomials into the framework of quantum groups and their representations (V.~Drinfeld~\cite{Drinfeld86,Drinfeld87}, M.~Jimbo~\cite{Jimbo85}), see also earlier preprints of N.~Reshetikhin~\cite{Resh87,Resh87_II} and many other references. 

Furthermore, while at generic values of the parameter $q$ of the quantum group $U_q(\mathfrak{g})$ of a simple Lie algebra $\mathfrak{g}$, its representation theory produces link invariants, root of unity values give rise to Witten--Reshetikhin--Turaev~\cite{Witten89,RT91,Resh91} and related invariants of 3-manifolds, see~\cite{MR1191386,Kuperberg96,KL01}, which can also be thought of as three-dimensional TQFTs. 

Quantum group $U_q(\mathfrak{g})$ is a Hopf algebra deformation of the universal enveloping algebra $U(\mathfrak{g})$ of a simple Lie algebra $\mathfrak{g}$.
For generic $q$, representation theory of $U_q(\mathfrak{g})$ gives rise to Reshetikhin--Turaev invariants $P_{\mathfrak{g}}(L)\in q^{\ell}\Z[q,q^{-1}]$ of knots and links $L$ in $\mathbb{R}^3$, where $\ell$ is a rational number that depends on $\mathfrak{g}$ and the linking number of $L$, see~\cite{Le00}. A key property of $U_q(\mathfrak{g})$ is that it is \emph{quasitriangular}.

\vspace{0.05in}

To define the Reshetikhin--Turaev invariant of a link, its connected components need to be labelled by irreducible representations of $U_q(\mathfrak{g})$. The latter are parametrized by positive integral weights $\lambda\in \Lambda^+$ of $\mathfrak{g}$.
    Quantum invariants of 3-manifolds (Witten--Reshetikhin--Turaev invariants) are given by an appropriate
    sum of these invariants when $q$ is a root of unity \cite{RT91,KM91}. 
    A 3-manifold is given by surgery on a link.

%
%

\subsection{Tangles and Reshetikhin--Turaev invariants}

 Reshetikhin--Turaev invariants~\cite{RT90} are defined for tangles, which are links with boundary (more carefully,  the invariant is usually defined for framed links and tangles).
  Tangles constitute a braided monoidal category $\mathsf{Tan}$. Composition of tangles is given by concatenation, while the tensor product on morphisms is given by placing tangles in parallel. Objects in the category of tangles are finite sequences of signs, which are orientations of a tangle at its endpoints. 

\input{mfig_001}

\input{mfig_002}

  \vspace{0.05in}

  The Reshetikhin--Turaev invariant $f(T)$ of a tangle $T$ is an intertwiner
  (homomorphism of representations) between
  tensor products of representations, read off from the endpoints of a tangle. If one picks a representation $V$ of the quantum group $U_q(\mathfrak{g})$, the invariant is an intertwiner between tensor products of $V$ and its dual $V^{\ast}$, according to the orientations of the endpoints, see an example in Figure~\ref{mfig_002}.

 To define the Reshetikhin--Turaev invariant in full generality, one first modifies the category $\mathsf{Tan}$ by labeling the components of a tangle by positive integral weights $\lambda\in \Lambda^+$. The Reshetikhin--Turaev invariant for a labelled tangle is built from intertwiners between tensor products of irreducible representations $V_{\lambda}$  of $U_q(\mathfrak{g})$.

\vspace{0.07in}

  A link $L$ is a tangle with the empty boundary. Reshetikhin--Turaev invariant
  $f(L): \C \lra \C$ is then a scalar, depending on $q$, see Figure~\ref{mfig_003}, and $f(L)\in \Z[q^{1/D},q^{-1/D}]$, where $D$ divides the determinant of the Cartan matrix of $\mathfrak{g}$, see~\cite{Le00}. 
  This integrality is a special property of Reshetikhin--Turaev invariants.

\input{mfig_003}

\subsection{Crane--Frenkel conjecture}

  Around 1994, Igor B. Frenkel (the graduate advisor of the second author) and Louis Crane proposed in \cite{CF94}:
 \begin{conjecture}[Crane--Frenkel]
 \label{conjecture:Crane-Frenkel}
 There exists a categorification of the quantum group $U_q(\mathfrak{sl}_2)$ at roots of unity giving rise to a 4D TQFT.
 \end{conjecture}

One motivation for the conjecture was that  
Floer homology was already known at the time. It has the Casson invariant as its Euler characteristic. Floer homology can be viewed as a 4-dimensional TQFT, defined for at least some 3-manifolds and 4-cobordisms. It was natural to wonder whether other quantum invariants of 3-manifolds can be realized as Euler characteristics.    

Another motivation came from geometric representation theory, with the discovery of the Kazhdan--Lusztig basis in the Hecke algebra in \cite{KL79} and its geometric interpretation via sheaves on flag varieties. This was followed by the Beilinson--Lusztig--MacPherson's~\cite{BLM90} geometric
    interpretation of $V^{\otimes k}$, for $V$ a fundamental representation
    of $U_q(\mathfrak{sl}_n)$, via sheaves on partial flag varieties, with generators $E_i,F_i$ of the quantum group acting by correspondences. Lusztig's geometric realization~\cite{Lusztig91} of Kashiwara--Lusztig bases~\cite{Kash94} of $U_q(\mathfrak{g})$ and bases of its irreducible representations and Lusztig's discovery of his bases of idempotented quantum groups~\cite{Lusztig10}. 

    I.~Frenkel's insight, beyond the above conjecture,
    was that positivity and integrality structure of these bases should be
    used to systematically lift Hecke algebra and quantum group elements to functors acting
    on categories that replace representations of quantum groups. On the TQFT level this should
    correspond to lifting quantum invariants one dimension up, from 3D to 4D (categorification).

\vspace{0.05in}

    While Conjecture~\ref{conjecture:Crane-Frenkel} is still open, significant progress in the past thirty 
    years has been made on:
  \begin{itemize}
    \item
    Link homology theories, which are four-dimensional counterparts of Reshetikhin--Turaev quantum link invariants,
    \item
    Categorification of quantum groups at generic $q$
    (A.~Lauda, A.~Lauda--Khovanov, R.~Rouquier and further foundational work by many researchers),
    \item
    Categorification at prime roots of unity (Y.~Qi, J.~Sussan, B.~Elias, Khovanov).
  \end{itemize}

\subsection{Semisimple versus triangulated}

Reshetikhin--Turaev link invariants are governed by \emph{semisimple} categories of representations of quantum groups $U_q(\mathfrak{g})$. In one dimension up (4D), categories cannot be semisimple. They are also unlikely to be abelian. 

In four dimensions, an \emph{extended} TQFT for link cobordisms requires one to assign:
\begin{itemize}
\item a category $\mcC_n$ to a plane with $n$ points (ignoring orientations, for simplicity),
\item a functor $F(T):\mcC_n\lra \mcC_m$ between these categories to a tangle $T$ with $n$ bottom and $m$ top boundary points,  and 
\item 
a natural transformation $F(T_1)\Rightarrow F(T_2)$ between these functors to a tangle cobordism between $T_1$ and $T_2$.
\end{itemize}
This assignment should assemble into a 2-functor 
\[
\mcF \ : \ \mathsf{Tan}_2 \lra \mathcal{NT}
\]
from the 2-category of tangle cobordisms to the 2-category $\mathcal{NT}$ of natural transformations (between functors between appropriate categories). Usually one wants to convert topological structures into something algebraic, so functors and natural transformations must be additive and defined over a field or a commutative ring. 

In particular, the braid group $\mathsf{B}_n$ on $n$ strands (the mapping class group of a plane with $n$ points) needs to act on $\mcC_n$, the category assigned to the plane with $n$ points. An action of a group $G$ on a semisimple category essentially just permutes its simple objects, thus reduces to a homomorphism $G\lra S_n$ into the symmetric group. Homomorphisms from braid groups $\mathsf{B}_n$ to symmetric groups $S_n$ are unlikely to be part of a sophisticated structure that carries interesting information about four-dimensional topology. The same argument applies to full four-dimensional TQFTs, replacing braid groups by mapping class groups of closed surfaces. In some cases one can expect that endomorphism rings $\End(L_i)$ of simple objects are not the ground field $\kk$ but a field or a division ring $D$ over $\kk$, in which case homomorphisms of $\mathsf{B}_n$ into $\Aut_{\kk}(D)$ or into cross-products of its direct products with the symmetric group may be available, but the crux of this informal argument (even just a gut feeling) that such homomorphisms cannot be upgraded to a sophisticated 4D TQFT remains. 

 This informal argument about unsuitability of semisimple categories and the need for triangulated categories in four-dimensional TQFTs can be found in the old paper of one of us~\cite[Section 6.5]{Kho02} 
    which further argues that interesting four-dimensional TQFTs are unlikely to assign abelian categories to surfaces. 
 For a much more recent and precise work we refer to Reutter~\cite{Reu20} and~\cite{RSP22}. 

\vspace{0.07in}

Replacing semisimple or abelian categories by triangulated categories removes the obstacle which is the lack of interesting categorical actions. Consider a ring $A$, take its category of modules $A\mathsf{-mod}$ and form the category $\mathcal{H}(A\mathsf{-mod})$ of finite length complexes of $A$-modules up to chain homotopies:  
\[
\xymatrix{
\ldots \ar[r]  & M^i 
\ar[r]^{d^i}  & M^{i+1} 
\ar[r]^{d^{i+1}} & M^{i+2} 
\ar[r]^{d^{i+2}} & \ldots, & \hspace{-0.5cm} 
\mbox{where } M^i\in A\mathsf{-mod}, \ \ 
d\circ d=0.
}
\]
One can come up with a specific finite-dimensional ring $A$ (interestingly, graded $A$-modules are essentially the same as double complexes) and a \emph{faithful} action of the braid group on the  category $\mathcal{H}(A\mathsf{-mod})$, see~\cite{KS02}, also~\cite{ST01} in the context of algebraic geometry.
Passing to the Grothendieck group of $\mathcal{H}(A\mathsf{-mod})$  recovers the Burau representation of the braid group (or the permutation representation, ignoring the extra grading. 
An additional grading on the category of modules turns the Grothendieck group into a $\Z[q,q^{-1}]$-module). Faithfulness holds in either case. 

There are many ways to construct algebras with a faithful braid group action on their homotopy categories but the above example seems minimal, in a sense. One can, for instance, ask for a finite-dimensional algebra $A_n$ over a field $\kk$ with a faithful action of the braid group $\mathsf{B}_n$ on $n$ strands on its homotopy category. In the example in~\cite{KS02} $\dim_{\kk}(A_n)=4n-6$, and we do not know examples of algebras of dimension less than $4n-6$ for $n\ge 3$ with a faithful action of $\mathsf{B}_n$ on their homotopy categories of complexes. 

Graded $A_n$ and $A$-modules correspond to bicomplexes, and one is working in the homotopy category of complexes over them. The latter category can be replaced by the stable category of tricomplexes, with a faithful braid group action on it~\cite{KQ20}. It is an open question whether this step into stable categories and tricomplexes can be extended beyond categorified Burau representation, to categorification of other braid group and quantum group representations. 

\vspace{0.1in} 

\subsection{The HOMFLYPT polynomial}
The Reshetikhin--Turaev invariant for the quantum group $U_q(\mathfrak{sl}(N))$ of $\mathfrak{sl}(N)$ and its fundamental $N$-dimensional representation $V$ is determined by the skein relation

\vspace{0.05in}

\begin{center}
\input{mfig_004}
\end{center}

  and normalization on the unknot:
  \begin{center}
  \begin{tikzpicture}[scale=0.6]
  \node at (-1,2) {$P_N\Big($};
  \draw[thick,->] (0.90,2) arc (0:360:0.50);
  \node at (7.55,2) {$\Big)=[N]=\dfrac{q^N-q^{-N}}{q-q^{-1}}=q^{N-1}+q^{N-3}+\ldots + q^{1-N}.$};
  \end{tikzpicture}
\end{center}  
Sometimes the normalization on the unknot is taken to be $1$, since otherwise the invariant of any nonempty link is divisible by $[N]$. The disjoint union with the unknot multiplies either invariant by $[N]$ and the above normalization is natural from the categorification viewpoint. 

  The skein relation above is due to the space of intertwiners (homomorphisms of quantum group representations) $V^{\otimes 2}\lra V^{\otimes 2}$
  being two-dimensional, so that any three maps are related by a linear equation. 
\vspace{0.05in}

  The HOMFLYPT polynomial~\cite{FYHLMO85,Przytycki88} is a 2-variable invariant $P(L)\in \Z[a^{\pm 1},b^{\pm 1}]$ given by replacing the coefficients in the above skein relation by  $a=q^N,b = q-q^{-1}$. One-variable specializations $P_N(L)$ of $P(L)$ have representation-theoretic interpretation via the quantum group $U_q(\mathfrak{sl}(N))$, as briefly explained earlier. Replacing $N$ by $-N$ and $q$ by $q^{-1}$ preserves the invariant, so one can restrict to $N\ge 0$. We record special cases: 
  \begin{itemize}
      \item $N=0$: $P_0(L)$ is the Alexander polynomial. 
      \item $N=1$: $P_1(L)$ is a trivial invariant. 
      \item $N=2$: $P_2(L)$ is the Jones polynomial. 
      \item $N=3$: $P_3(L)$ is the Kuperberg $\mathfrak{sl}(3)$ quantum invariant~\cite{Kup96}. 
  \end{itemize}

\subsection{MOY graphs and their invariants}
The above specializations $P_N(L)$ take values in the ring $\Z[q,q^{-1}]$ of Laurent polynomials with \emph{integer} coefficients. It is possible to reduce links to linear combinations of planar objects (webs or graphs) on which the invariant takes values in $\Z_+[q,q^{-1}]$, the semiring of Laurent polynomials with \emph{non-negative integer} coefficients. 
 Decompose $V^{\otimes 2}\cong \Lambda^2_q(V)\oplus S^2_q(V)$ into the sum of two irreducible representations---the second quantum exterior and symmetric powers of $V$---and
  consider projection and inclusion operators 
  \[ 
  \xymatrix{
  V^{\otimes 2}  \ar[r]^p & 
  \Lambda^2_q(V) \ar[r]^{\iota} & 
  V^{\otimes 2}
  }
  \]
  scaled so that $p\circ\iota = (q+q^{-1})
  \id$ is the identity map of $\Lambda^2_q(V)$ times $q+q^{-1}$, see Figure~\ref{mfig_005}. 
\begin{center}
\input{mfig_005}
\end{center}
From these two basic pieces one can assemble oriented planar graphs with thin and thin edges and all vertices of valency three, as in Figure~\ref{mfig_005} left and center, with one thick and two thin edges at each vertex. Vertices correspond to the intertwiners $p$ and $\iota$. 
These are the simplest instances of Murakami--Ohtsuki--Yamada (MOY) graphs. We think of thin edges as carrying label one (for $V$) and thick edges as carrying label two (for $\Lambda_q^2 V$). Such a graph with boundary defines an intertwiner between tensor products of $V$, $\Lambda_q^2 V$ and their duals. A closed MOY graph $\Gamma$ defines an endomorphism of the trivial representation of $U_q(\mathfrak{sl}(N))$, thus a function $P_N(\Gamma)\in \C(q)$. 

\input{mfig_006}

The invariant $P_N(\Gamma)$ is both integral and positive, see \cite{MOY98}: 
\begin{theorem}[Murakami--Ohtsuki--Yamada] 
\label{thm_MOY}
For any planar MOY graph $\Gamma$ as above, its
  invariant $P_N(\Gamma)\in \Z_+[q,q^{-1}]$.
\end{theorem}
Integrality and positivity properties of $P_N(\Gamma)$ are key to its categorification and categorification of the corresponding link invariants $P_N(L)$. 

Let us first observe that the link invariant $P_N(L)$ reduces to the invariants $P_N(\Gamma)$ of MOY planar graphs via skein relations: 

\begin{center}
\input{mfig_008}
\end{center}
implying the HOMFLYPT relation 
\begin{center}
\input{mfig_004}
\end{center}

  More generally, one can allow lines of thickness
  $a$, $1\le a\le N-1$, which correspond to fundamental representations $\Lambda^a_q V$ of the quantum group $U_q(\mathfrak{sl}(N))$ and $(a,b,a+b)$-vertices, see the top row of Figure~\ref{mfig_007}, which
  denote the scaled projection and inclusion between the corresponding representations 
  \[
  \xymatrix{
  \Lambda^a_q V \otimes \Lambda^b_q V 
  \ar[r]^{\hspace{0.35cm} p_{a,b}} &  
  \Lambda^{a+b}_q V, 
  }
  \qquad 
  \quad 
  \xymatrix{
  \Lambda^{a+b}_q V \ar[r]^{\iota_{a,b}\hspace{0.50cm}} & 
  \Lambda^a_q V \otimes \Lambda^b_q V.
  }
  \]
  normalized so that the composition 
  \[
  p_{a,b}\circ \iota_{a,b} = 
  \begin{bmatrix}
      a+b \\
      b \\ 
  \end{bmatrix}
\,\id = \frac{[a+b]!}{[a]![b]!}\cdot\id , 
\qquad 
[a]! :=\prod_{c=1}^a[c], 
\qquad 
[c]=\frac{q^c-q^{-c}}{q-q^{-1}}
  \]
  is the identity endomorphism of $\Lambda_q^{a+b}V$ scaled by the $q$-binomial coefficient. It is also natural to allow thickness $N$ as well and replace $\mathfrak{sl}(N)$ by the Lie algebra $\mathfrak{gl}(N)$. The exterior power $\Lambda^N V$ is a trivial representation of $\mathfrak{sl}(N)$ but not of $\mathfrak{gl}(N)$, and likewise for their quantum deformations, and vertices $(a,N-a)$ are now allowed, with lines of thickness $a$ and $N-a$ merging into a line of thickness $N$.

Murakami--Ohtsuki--Yamada planar graphs or webs $\Gamma$ with edges of arbitrary thickness from $1$ to $N$ have vertices $(a,b)$ as above ($a+b\le N$). The quantum invariant $P_N(\Gamma)$ of an MOY graph with edges of arbitrary thickness takes values in $\Z_+[q,q^{-1}]$, and Theorem~\ref{thm_MOY} holds for such graphs as well. 

\begin{center}
\input{mfig_007}
\end{center}

Consider a link $L$ with components labelled by numbers in $\{1,2,\dots, N\}$ and pick a planar projection. The analogue of the skein relations above $(I)$ and $(II)$ are the relations in Figures~\ref{mfig_023},~\ref{mfig_024}, see for instance~\cite{Wu12}. 
These skein relations allow to define $P_N(L)\in \Z[q,q^{-1}]$ for links with components labelled by numbers from $1$ to $N$.  

\input{mfig_023}

\input{mfig_024}


\subsection{Matrix factorizations and link homology}

\quad 

{\bf Idea:}
To realize $P_N(L)$ as the Euler characteristic of a bigraded homology theory $H_N(L)$  of links,
first build homology $H_N(\Gamma)$ for planar graphs $\Gamma$ as \emph{singly-graded} vector spaces. There must be an equality of Laurent polynomials
\[ P_N(\Gamma) \ = \ \mathsf{gdim} H_N(\Gamma),\]
where for a $\Z$-graded vector space $V=\oplus_i V_i$, the graded dimension $\mathsf{gdim} V := \sum_i \dim(V_i)q^i.$
Then lift skein relations above to the long exact sequence in Figure~\ref{mfig_009}. 

Commutative ground ring $R$ of the theory may be different from a field, in which case one expects $H_N(\Gamma)$ to be a free graded $R$-module of graded rank $P_N(\Gamma)$.

\input{mfig_009}

 This idea was successfully realized by L.~Rozansky and one of us back in 2004 in~\cite{KR04}. To define homology groups (or state spaces) 
  $H_N(\Gamma)$ of planar graphs $\Gamma$ \emph{matrix factorizations} were used. Start with a polynomial ring $S=\kk[x_1,\ldots, x_k]$ in several variables over a field $\kk$, where $\mathsf{char}(\kk)=0$. A polynomial $\omega\in S$ is called a \emph{potential} if the ideal $(\partial \omega/\partial x_1,\dots, \partial \omega/\partial x_k)$ has finite codimension in $S$. Informally, this means that $\omega$ is sufficiently generic. In this case the quotient algebra $J_{\omega}:=S/(\partial \omega/\partial x_1,\dots, \partial \omega/\partial x_k)$ is known to be Frobenius, via the Grothendieck residue construction. The quotient algebra is called \emph{the Milnor algebra} of $\omega$. 
  
   For a potential $\omega\in S$, consider 2-periodic generalized complexes $M$ of free $S$-modules

\begin{center}
\input{mfig_011}    
\end{center}

\noindent 
    and maps between them such that $d^2(m) = \omega m$. Modulo homotopies, these constitute a triangulated
    category of matrix factorizations $MF_{\omega}$. When matrix factorizations $M,M'$ have finite ranks (meaning $M^0,M'^0$ are finite rank free $S$-modules) hom spaces $\Hom(M,M')$ in the category $MF_{\omega}$ are finite-dimensional $\kk$-vector spaces, due to multiplications by $\partial \omega/\partial x_i$ being homotopic to $0$. The Milnor algebra $J_{\omega}$ then acts on $\Hom(M,M')$, implying finite-dimensionality of hom spaces. 

    \vspace{0.05in}

 To a potential $\omega$ there is assigned a two-dimensional TQFT with corners, built out of categories of matrix factorizations for potentials that are signed sums of copies of $\omega$ over several sets of variables. This construction goes back to Kapustin--Li \cite{KL03}, see~\cite{CM16} for a more recent treatment, while the one-variable case is thoroughly explained in~\cite{KR04}.  Without going into full details and specializing to $\kk=\Q$ and $\omega=x^{N+1}$, 
 to an arc (viewed as a one-manifold with boundary) one assigns a factorization $L$ for $\omega = x_1^{N+1}-x_2^{N+1}$, with $S=\Q[x_1,x_2]$:

\input{complex_001}

Term $x_1-x_2$ ensures that $x_1,x_2$ act the same up to homotopy, when viewed as endomorphisms of factorization $L$ in the homotopy category. 
One thinks of $L$ as the identity factorization, implementing the identity functor on the category of matrix factorizations. Namely, tensoring $L$, say over the subring $\Q[x_1]$ with a factorization $M$, in variables $x_1$ and some other variables (not $x_2$, say $x_3,x_4$) for the potential $-x_1^{N+1}+\omega'(x_3,x_4)$ results in
the factorization $M':=L\otimes_{\Q[x_1]}M$ over $\Q[x_2,x_3,x_4]$, and one can check that upon substitution $x_2\mapsto x_1$ factorizations $M,M'$ are naturally isomorphic in the homotopy category. Factorization $M'$ is of infinite rank over $\Q[x_2,x_3,x_4]$, but upon removing contractible summands it becomes finite rank and isomorphic to $M$. 

\begin{center}
\input{mfig_012}
\end{center}

 Closing up the arc into a circle and equating $x_1=x_2$, see Figure~\ref{mfig_012} on the right,  gives a complex 
 \[
\xymatrix{
\Q[x_1] \ar[rr]^{(N+1)x^N} & & \Q[x_1] \ar[rr]^0 & & \Q[x_1],
}  
 \]
 where the ring is reduced to $\Q[x_1,x_2]/(x_1-x_2)\cong \Q[x_1]$. This a general feature of building topological theories from matrix factorizations: closing up a factorization with boundary points and equating the variables, as long as the potentials at the endpoints match, results in a two-periodic complex, with $d^2=0$, rather than just a factorization, due to the cancellation of terms in $\omega$ upon identifying the variables. 

 For the circle as above, the homology group of the complex is $\Q[x]/(x^N)\cong H^{\ast}(\mathbb{CP}^{N-1})$ (the other homology group is $0$). Algebra $\Q[x]/(x^N)$ is commutative Frobenius and gives rise to a 2D TQFT once a nondegenerate trace on it is chosen. Building an extended TQFT from matrix factorizations is discussed in~\cite{KR04}, also see~\cite{CM16} and references there.

\vspace{0.03in}

In the above example one can work with graded polynomial rings with $\deg(x_j)=2$ and factorizations $\sum_j \pm x_j^{N+1}$. Then the graded degree of the vector space of the circle is 
 \[\mathsf{gdeg} (H_N(\mathsf{circle}))=
1+q^2+\ldots + q^{2(N-1)}=q^{N-1}[N]=q^{N-1}P_N(\mathsf{unknot}).
\] 
This equality points to a possible match between the 2D TQFT for the one-variable potential $x^{N+1}$ and homology of links and planar graphs, for the simplest possible link and MOY graph, which is the circle in the plane (of thickness     1). Up to a power of $q$, the quantum invariant of the unknot equals the graded dimension of the commutative Frobenius algebra that the matrix factorization TQFT for the potential $\omega=x^{N+1}$ assigns to the circle. That algebra is also isomorphic to the cohomology ring of 
$\mathbb{CP}^{N-1}$. 

\begin{center}
\input{mfig_013}
\end{center}

To move beyond circles, recall that, for link diagrams with components colored by $V$, the corresponding MOY graphs have edges of thickness one and two only, with any thickness two edge having two thin ``in'' edges and two thin ``out'' edges, see Figure~\ref{mfig_013} on the left.

Matrix  factorization associated to a neighbourhood of a double edge should have the potential  
\begin{equation}\label{eq_omega}\omega = x_1^{N+1}+x_2^{N+1}-x_3^{N+1}-x_4^{N+1},
\end{equation} 
see Figure~\ref{mfig_013}  on the left. The four variables $x_1,x_2,x_3,x_4$ are assigned to the endpoints of the diagram. The term $x_i^{N+1}$ enters $\omega$ with the plus sign, respectively the minus sign, if the orientation at that point is out, respectively into, the diagram. Potential $\omega$ is the difference of two terms, and the first term $x_1^{N+1}+x_2^{N+1}$ can be written as a polynomial in the elementary symmetric functions: $x_1^{N+1}+x_2^{N+1}= g(x_1+x_2,x_1x_2).$  

To write the identity factorization for the 2-variable polynomial $g$ we decompose the difference of two $g$'s for two different sets of variables:  
\begin{eqnarray*} 
g(y_1,y_2)-g(z_1,z_2)& = & (g(y_1,y_2)-g(z_1,y_2))+(g(z_1,y_2)-g(z_1,z_2)) \\
& = & (y_1-z_1)u_1(\underline{y},\underline{z}) + (y_2-z_2)u_2(\underline{y},\underline{z})
\end{eqnarray*}
for some polynomials $u_1,u_2$ in the four variables $y_1,y_2,z_1,z_2$. 
Tensor the two factorizations 
\[
\xymatrix{
(S  \ar[r]^{u_1} & 
S \ar[r]^{y_1 - z_1} & 
S) 
}
\quad 
\mbox{ and }
\quad 
\xymatrix{
(S  \ar[r]^{u_2} & 
S \ar[r]^{y_2 - z_2} & 
S) 
}
\]
to get a factorization associated to a double edge $j$ of $\Gamma$:
\[
\xymatrix{
M_j \ := \  (S \ar[r]^{\hspace{0.6cm} u_1} & 
S \ar[r]^{y_1 - z_1 \hspace{0.60cm}} & 
S) \otimes_S 
(S \ar[r]^{\hspace{0.50cm} u_2} & 
S \ar[r]^{y_2-z_2 \hspace{0.1cm}} & 
S). 
}
\]
Notice that $M_j$ has the potential given by \eqref{eq_omega}. 

For a general planar MOY graph $\Gamma$ with edges of thickness one and two, place marks on thin edges, with at least one mark on each edge. Denote by $I$ the set of marks and consider the ring $\Q[x_i]_{i\in I}$. To a thick edge there is assigned a factorization with the potential $x_{i_1}^{N+1}+x_{i_2}^{N+1}-x_{i_3}^{N+1}-x_{i_4}^{N+1}$, to a thin arc -- an identity factorization as described earlier, see Figure~\ref{mfig_014} on the left. In the Figure~\ref{mfig_014} example, the factorization for the thin arc shown has potential $x_{i_5}^{N+1}-x_{i_1}^{N+1}$. 

\vspace{0.05in} 

  Given a graph $\Gamma$ and a set of markings $I$ of thin edges as described, tensor together factorizations $M_j$ over all double edges $j$ 
   and factorizations for thin arcs (if some thin edge carries more than one mark) to get the product factorization $ M_{\Gamma} $ . Then the differential $D= \displaystyle{\sum_j} d_j$ in $M_{\Gamma}$ satisfies
  $D^2 = \displaystyle{\sum_j} d_j^2 
  = \displaystyle{\sum_i} (x_i^{N+1}-x_i^{N+1}) = 0$.
  
\begin{center}
\input{mfig_014}
\end{center}

  Define $H_N(\Gamma) := H_N(M_{\Gamma},D)$. This is our homology (or state space) of the graph $\Gamma$. Placing extra marks on any thin edges of $\Gamma$ does not change $H_N(\Gamma)$. 

\vspace{0.07in}

From \cite{KR04}, we cite the following: 
\begin{theorem}[Khovanov--Rozansky] $H_N(\Gamma)$ lives is a single $\Z/2$ degree and 
  $\mathsf{gdim}(H_N(\Gamma))= P_N(\Gamma).$
  \end{theorem}

  To extend the homology to links, consider resolutions $\Gamma_0,\Gamma_1$ of a crossing in Figure~\ref{mfig_015}. There are homomorphisms of factorizations \[
  \chi_0: M_{\Gamma_0}\lra M_{\Gamma_1}, 
  \qquad 
  \qquad 
  \chi_1: M_{\Gamma_1}\lra M_{\Gamma_0}, 
  \] 
  and to positive and negative crossings one assigns two-step complexes of factorizations with the differential given by $\chi_0,\chi_1$, respectively, see Figure~\ref{mfig_015} on the right.  
\begin{center}
\input{mfig_015}
\end{center}
Homology of a link diagram $D$ is given by tensoring these two-step complexes over all crossings of $D$, computing the homology of all terms (for the inner differential in factorizations), and then taking the homology again for the differential build from maps $\chi_0,\chi_1$ over all crossings. 

\vspace{0.07in} 

From \cite{KR04}, we also recall: 

\begin{theorem}[Khovanov--Rozansky]
The resulting homology does not depend on the choice of a link diagram $D$ of an oriented link $L$ and can be denoted $H_N(L)$.  
\begin{enumerate}
    \item 
$H_N(L)$ is bigraded and its Euler characteristic 
\[\chi_N(L) = \displaystyle{\sum_{i,j}} (-1)^i q^j \dim H^{i,j}_N(L)=P_N(L)
\]
is the Reshetikhin--Turaev link invariant for the fundamental $\mathsf{SL}(N)$ representation. 
    \item 
$H_N(L)$ is functorial under link cobordisms. 
\end{enumerate}
\end{theorem}

Homology theory $H_N$ of links and link cobordisms can be viewed as a categorification of quantum invariant $P_N$, the latter a one-variable specialization of the HOMFLYPT polynomial. The functoriality under link cobordisms, shown in~\cite{KR04}, is up to overall scaling by elements of $\Q^{\ast}$. The theory extends to tangles and tangle cobordisms as well. 

\vspace{0.05in} 
 
Y.~Yonezawa~\cite{Yon11} and H.~Wu~\cite{Wu14} extended the homology from \emph{thin} MOY graphs to arbitrary MOY graphs, with general $(a,b,a+b)$ trivalent vertices, and from a categorification of the link invariant $P_N$ to categorification of Reshetikhin--Turaev invariants where components of a link are colored by arbitrary quantum exterior powers $\Lambda^a_q V$ of the fundamental representation $V\cong \C^N$, $1\le a\le N$.  

In their construction the homology of the unknot colored $a$ is the cohomology of the Grassmannian of $a$-planes in $\C^N$: 
 \begin{center}
  \begin{tikzpicture}[scale=0.6]
  \node at (-1.2,2) {$H_N\Big($};
  \node at (0.8,2.8) {$a$};
  \draw[thick,->] (0.70,2) arc (0:360:0.50);
  \node at (4.00,2) {$\Big) =\mathsf{H}^{\ast}(\mathsf{Gr}(a,N),\Q)$.} ; 
  \end{tikzpicture}
\end{center} 
This cohomology ring is a commutative Frobenius algebra and gives rise to a 2D TQFT. 
Standard merge, split, birth and death cobordisms between unnested circles in the plane and relations on compositions of these cobordisms show that in a  functorial link homology theory homology of the unknot is a commutative Frobenius algebra, and for the theory above that algebra is the cohomology ring of the complex Grassmannian. 

\vspace{0.1in}

This theory admits a generalization where the potential $x^{N+1}$ is replaced by the potential $\omega = x^{N+1} + a_1 X^N + \ldots + a_N x$, where $a_i$ is a formal variable of degree $2i$. The potential is then homogeneous and defined over the ground ring $R_N=\Q[a_1,\dots, a_N]$. The 2D TQFT for matrix factorizations and link homology construction can be generalized to the  ground ring $R_N$. In the resulting homology theory $\widetilde{H}_N$ the ring $R_N$, which is the homology of the empty link, can be interpreted as the $U(N)$-equivariant cohomology of a point $p$: 
\[
 \widetilde{H}_N(\emptyset) \ = \ R_N \  \cong \ \mathsf{H}^{\ast}_{U(N)}(p,\Q). 
\]
Homology of an $a$-labelled unknot is then the $U(N)$-equivariant cohomology of the complex Grassmannian: 
\begin{center}
  \begin{tikzpicture}[scale=0.6]
  \node at (-1.2,2) {$\widetilde{H}_N\Big($};
  \node at (0.8,2.8) {$a$};
  \draw[thick,->] (0.70,2) arc (0:360:0.50);
  \node at (4.45,2) {$\Big) =\mathsf{H}_{U(N)}^{\ast}(\mathsf{Gr}(a,N),\Q)$.} ; 
  \end{tikzpicture}
\end{center} 
We refer to papers~\cite{Gornik04,Kra10,WuEquiv12,Wu_equivar15} for the construction of equivariant link homology, see also~\cite{LL16,Ras15,MY19} for some applications. Deforming the potential proved to be a prolific idea, also leading to a categorification of the bigraded HOMFLYPT polynomial~\cite{KR05} rather than of its singly-graded specializations $P_N$. This triply-graded link homology theory is closely related to the category of Soergel bimodules~\cite{Kho_triply07}, which is a fundamental object at the center of modern geometric representation theory~\cite{EMTW20}. We refer to~\cite{GSV05,Witt18,Aga22} and many other papers for physical interpretations of $\mathsf{SL}(N)$ link homology. 

Bigraded and triply-graded categorifications of the HOMFLYPT polynomial and its specializations, including for torus and algebraic links, relate to deep structures in representation theory and geometry \cite{GORS14,GHM21}. 

A refinement of bigraded and triply-graded homology groups, now known as \emph{$y$-ified homology}, was discovered by J.~Batson and C.~Seed~\cite{BS_spectral15} in the $N=2$ case, S.~Cautis and J.~Kamnitzer \cite{CK_coherent17} in the geometric language for any $N$ and E.~Gorsky and M.~Hogancamp for the HOMFLYPT homology~\cite{GH_yification22}, also see~\cite{CLS_Rickard20} for a generalization and an approach via the KLR categories and a related paper~\cite{HRW21}. This homology exhibits additional symmetries and relates link homology to Hilbert schemes. Sophisticated relations between flavors of link homology and equivariant matrix factorizations have been investigated by A.~Oblomkov and L.~Rozansky in~\cite{OR20} and follow-up papers, see~\cite{Oblomkov19} for a review. 

\vspace{0.1in}


\subsection{Model scenario for categorification of Reshetikhin--Turaev link invariants}
\label{subsubsection:model_scenario_RT_inv}

\quad 

Reshetikhin--Turaev link invariant functor~\cite{RT90} assigns the tensor product $V_{\underline{\lambda}}=V_{\lambda_1}\otimes \dots \otimes V_{\lambda_n}$ of irreducible $U_q(\mathfrak{g})$ representations to a plane with $n$ points marked by positive integral weights $\lambda_1,\dots, \lambda_n$, where $\underline{\lambda}:=(\lambda_1,\dots, \lambda_n)$, and taking orientations into account. The invariant of a decorated tangle $T$ with $\underline{\lambda}$ and $\underline{\mu}$ endpoint label sequences at the bottom and top planes is an intertwiner 
\begin{equation}
    f(T) \ : \ V_{\underline{\lambda}}\lra V_{\underline{\mu}}.
\end{equation}

Let us explain the best possible scenario for categorification of Reshetikhin--Turaev link invariants. 
Upon a categorification, the tensor product $V_{\underline{\lambda}}$ of representations should be replaced by a triangulated category $\mathcal{C}_{\underline{\lambda}}$. 
The Grothendieck group $K_0(\mathcal{C}_{\underline{\lambda}})$ of that category
 should be related to $V_{\underline{\lambda}}$, for instance, be a $\Z[q,q^{-1}]$-lattice in $V_{\underline{\lambda}}$, so that 
 \begin{equation}
     K_0(\mathcal{C}_{\underline{\lambda}}) \otimes_{\Z[q,q^{-1}]} \C(q) \ \cong \ V_{\underline{\lambda}} 
 \end{equation}
 (assuming $U_q(\mathfrak{g})$ is defined over the field of rational functions $\Q(q)$). 

\vspace{0.05in} 

To a collection of $n$ points on the plane labelled by $\underline{\lambda}$ assign the category $\mcC_{\underline{\lambda}}$ as above. 
 To a tangle $T$ there should be assigned an exact functor $\mathcal{F}(T):\mcC_{\underline{\lambda}}\lra \mcC_{\underline{\mu}}$, see Figure~\ref{mfig_020} on the left. On the Grothendieck group level the functor $\mathcal{F}(T)$ should give the Reshetikhin--Turaev tangle invariant, $[\mathcal{F}(T)]=f(T)$. This equality can be written as a commutative diagram

\[
\xymatrix{
K_0(\mathcal{C}_{\underline{\mu}})  \ar@{^{(}->}[rr]^{}&   & 
V_{\underline{\mu}} \\ 
& & \\ 
K_0(\mathcal{C}_{\underline{\lambda}}) \ar@{^{(}->}[rr]^{} \ar[uu]^{[\mcF(T)]} & & V_{\underline{\lambda}} \ar[uu]_{f(T).} \\
}
\]
Horizontal inclusions are those of $\Z[q,q^{-1}]$-modules and become isomorphisms upon tensoring with $\Q(q)$. Alternatively,  $K_0(\mathcal{C}_{\underline{\lambda}}) \otimes_{\Z[q,q^{-1}]} \C(q)$ could be a subspace of $V_{\underline{\lambda}}$, such as the subspace $\Inv_{U_q(\mathfrak{g})}(V_{\underline{\lambda}})$ of quantum group invariants in the tensor product.  

For composable tangles $T',T$ there should be fixed an isomorphism of functors $F(T'\circ T)\cong F(T')\circ F(T)$. 

\begin{center}
\input{mfig_020}    
\end{center}

A tangle is an oriented decorated one-manifold properly embedded in $\R^2\times [0,1]$. 
A tangle cobordism between oriented tangles $T_0,T_1$ is a surface $S$ with boundary and corners properly embedded in $\R^2\times [0,1]^2$ such that its boundary is the union of $T_0,T_1$ on the intersections $S\cap (\R^2\times[0,1]\times \{i\})$, $i=0,1$, and the union of product tangles (finite sets of points times $[0,1]$) on the intersections $\R^2\times \{i\}\times [0,1]$, $i=0,1$.  A tangle cobordism is depicted rather schematically in Figure~\ref{mfig_020}, in the middle. Since tangles are additionally decorated by positive integral weights $\lambda$, connected components of $S$ must be decorated by integral weights as well, with boundary decorations induced by those of $S$.  
 
 To a tangle cobordism $S$ between
 tangles $T_0,T_1$ there should be assigned a natural transformation of functors
 $F(S):\mathcal{F}(T_0)\longrightarrow \mathcal{F}(T_1)$.

\vspace{0.05in} 

As one varies over all collections of labelled points in the plane, functors $\mcF(T)$ and natural transformations $\mcF(S)$, these should fit into 
 a \emph{2-functor} $\mcF$ from the \emph{2-category of tangle cobordisms}
 to the 2-category of natural transformations between exact functors. (Tangle cobordisms are oriented and decorated by positive integral weights of $\mathfrak{g}$.) Having such a 2-functor is the most functorial scenario for a link homology theory. 

\begin{center}
\input{mfig_010}
\end{center}

In the matrix factorization approach to link homology, described earlier, the link homology does extend to a 2-functor from the category of tangle cobordisms to the category of bigraded vector spaces (or bigraded $R_N$-modules, for the equivariant theory).

The Lie algebra $\mathfrak{g}=\mathfrak{sl}(N)$ and  all components are colored by $V$. For $n$ points on the plane one forms the category $MF_{\omega}$ with the potential $\omega = \displaystyle{\sum_{i=1}^n} \pm x_i^{N+1}$, where the signs depend on orientations of points. Then one forgets that $MF_{\omega}$ is a triangulated category, views it only as an additive category, and forms the homotopy category of complexes $\mathcal{C}_n=\mathcal{H}(MF_{\omega})$ over it. To a tangle $T$ there is assigned a complex of matrix factorizations and to a tangle cobordism -- a homomorphism of complexes. The entire construction results in a 2-functor as above: 

\begin{center}
\begin{tikzpicture}[scale=0.6]
\node at (0,2.5) {2-category $2\Tan$ of};
\node at (0,1.5) {Tangle cobordisms};

\node at (6,2.65) {$\mcF$};

\draw[thick] (4,2.1) -- (8,2.1);
\draw[thick] (4,1.9) -- (8,1.9);

\draw[thick] (7.75,2.5) -- (8.25,2);
\draw[thick] (7.75,1.5) -- (8.25,2);

\node at (14,2.5) {2-category $\CMF$ of (complexes of)};
\node at (14,1.5) {matrix factorizations.};
\end{tikzpicture}
\end{center}

Objects of the category $\mathsf{CMF}$ are potentials $\omega$ in finite sets of variables $\underline{x}$, morphisms from $\omega_1(\underline{x})$ to $\omega_2(\underline{y})$ are complexes of matrix factorizations with potential $\omega_2(\underline{y})-\omega_1(\underline{x})$, and two-morphisms are homomorphisms of complexes of factorizations modulo chain homotopies. 

The 2-functor $\mcF$ above is $\Q^{\ast}$-projective, with the map associated to a link cobordism well-defined up to scaling by elements of $\Q^{\ast}$. Extending T.~Sano's approach to strict cobordism invariance~\cite{Sano20} from $N=2$ to any $N$ should be one way to resolve the $\Q^\ast$-indeterminacy. 

To get a homomorphism from the Grothendieck group of the category $\mcC_n$ to the tensor product $V_{\underline{\lambda}}$, more precisely, to its subspace of invariants $\mathsf{Inv}(V_{\underline{\lambda}})$, one should restrict to the subcategory generated by matrix factorizations associated to planar graphs with a given boundary.

Versions of $\SL(N)$ homology and the 2-functor $\mathcal{F}$ can also be recovered from
\begin{itemize}
\item Parabolic-singular blocks of highest weight categories for $\mathfrak{sl}(k)$, where $k$
depends on $\underline{\lambda}$, see~\cite{Sussan07,MS09},
and, in the $\mathfrak{sl}(2)$ case,  \cite{BFK99,Str_Duke_05,FKS_quantum06,Str_perverse09,BS11,CK14}. 
\item  Fukaya--Floer categories of  quiver varieties  for $\SL(k)$, $k$ a function of weights $\underline{\lambda}$ ($N=2$ case~\cite{SS06,Reza09,AS19}) and arbitrary $N$ case restricted to braid and braid closures \cite{Man07},  see also~\cite{WW15}. 

\item Derived categories of coherent sheaves on the convolution varieties of affine Grassmannians; $N=2$ case:~\cite{CK_Duke08} and arbitrary $N$: \cite{CK_slm08}, and also on quiver varieties~\cite{AN23}. For a recent physics-motivated approach, see M.~Aganagic~\cite{Aga22,Agan21}. 
 
\item Various categories of importance in representation theory, see~\cite{Web17} and related papers.
\end{itemize}
The principle here is that in each of these cases there is (or expected to exist) a categorical action
of the 2-category of tangle cobordisms $\mathsf{Tan}_2$ on the corresponding derived or homotopy categories, giving a  2-functor similar to $(I)$ and $(II)$ in Figure~\ref{mfig_010}. Tangles act by functors on the categories and tangle cobordisms act by natural transformations.
Often, the corresponding categories for various examples are equivalent~\cite{MW_Howe18}, sometimes with minor modifications. One expects the equivalences to respect the 2-functor actions.

Motivation to relate quiver varieties (Kronheimer--Nakajima varieties) and link homology came from observing that (a) homology groups of $\mathsf{SL}(k)$ quiver varieties carry an action of that Lie algebra and can be identified with weight spaces of $\mathsf{SL}(k)$ representations~\cite{Naka94}, (b) due to the level-rank duality, those weight spaces can be identified with the invariants of tensor products of exterior powers $\Lambda^a V$ of the fundamental $\mathfrak{sl}(N)$ representation $V$, where $k$ depends on the sequence $(a_1,\dots, a_n)$ of these parameters in the tensor product. One then expects that replacing homology by the derived category of coherent sheaves will result in a categorical action of the Lie algebra and a commuting categorical action of tangles and tangle cobordisms. This program was realized in~\cite{CK_slm08} restricting to tensor products of $V$ and $V^{\ast}$ and using  related varieties which are suitable convolution varieties of the affine Grassmannian and should contain quiver varieties as open subvarities. 

    \vspace{0.05in} 

{\it Extended TQFTs and biadjoint functors:}
    When one has an extended TQFT in dimension $n$, to an $(n-2)$-manifold $M$ assign a category $\mcF(M)$, to an $(n-1)$-cobordism $K$ a functor $\mcF(K):\mcF(\partial_0 K)\lra \mcF(\partial_1 K)$, to an $n$-cobordism $L$ with corners a natural transformation $\mcF(\partial_0 L)\lra \mcF(\partial_1 L)$.

    \vspace{0.05in} 

    For each $(n-1)$-cobordism $K$ there is a dual cobordism $K^{\ast}$ from $\partial_1 K$ to $\partial_0 K$ (reflect $K$). There are four canonical $n$-cobordisms between compositions $K K^{\ast}$, $K^{\ast}K$ and identity cobordisms $\id_{\partial_1 K}$ and   $\id_{\partial_0 K}$. Applying $\mcF$ to them tells us that functors $\mcF(K)$ and $\mcF(K^{\ast})$ are biadjoint (both left and right adjoint). 

    \vspace{0.05in}

    Examples of biadjoint functors appear in algebraic geometry (Fourier--Mukai kernels between Calabi--Yau manifolds) and in symplectic topology (convolutions with Lagrangians in Fukaya--Floer categories). In particular, quiver varieties are Calabi--Yau (around their compact part) and their derived categories of coherent sheaves admit plenty of biadjoint functors.     
In three dimensions, the model example of a TQFT is the Witten--Reshetikhin--Turaev theory, for which the categories assigned to $3-2=1$-dimensional manifolds are semisimple. Any linear functor between these categories has a biadjoint functor. In higher dimensions, for an extended TQFT, one expects non-semisimple, likely triangulated, categories associated to $n-2$-dimensional manifolds. For an exact functor between triangulated categories having a biadjoint is a strong condition. It is then natural to pay special attention to Calabi--Yau varieties (and their derived categories of coherent sheaves), Fukaya--Floer categories, and suitable categories of representations built out of symmetric Frobenius algebras, for these types of categories provide large supply of biadjoint functors that can be constructed naturally~\cite{Kho02}. 

\vspace{0.07in} 

Reshetikhin--Turaev $\SL(N)$ (or $\GL(N)$) invariants for fundamental representations $\Lambda^a_q V$ are  distinguished by the existence of  {\emph{positive integral}} diagrammatical calculus of intertwiners (\emph{MOY calculus}). It guides the categorification of Reshetikhin--Turaev invariants. Positivity and integrality property  is missing already for other representations of quantum $\SL(2)$, since $q$-spin networks are not positive and have denominators~\cite{KL94}. This creates serious problems trying to  extend the above approaches  beyond $(\SL(N),\Lambda_q^{\ast}V)$ case. 

\vspace{0.10in}

In the remarkable work~\cite{Web17}, Webster categorified 
Reshetikhin--Turaev invariants for any simple Lie algebra $\mathfrak{g}$ and any labelling of a link's components by irreducible $\mathfrak{g}$-modules. Categorification of quantum groups~\cite{KL01,rouquier2008} was one of the motivations for Webster's construction. It is likely that for $\mathfrak{sl}(N)$ and coloring by minuscule representations his construction gives homology isomorphic to those coming from matrix factorizations (and those from Robert--Wagner foam evaluation). 

Beyond the minuscule representation case, Webster homology and a number of other known link homology theories  do not fully extend or not known to extend to link cobordisms. 

For a subset of these theories~\cite{FSS_fractional12,Rozan_Jones_Wenzl14,CK_Jones_Wenzl12,Web17,GH_yification22,SS_cat22,OR19}, the reason is the following. Take a functorial link homology theory (and subject to TQFT assumption) and consider standardly embedded cobordisms in $\R^3\times [0,1]$ between unlinks (disjoint unions of unknots). Restricting a link homology theory to these cobordisms results in a 2D TQFT. Consequently, 
homology $A$ of the unknot must be a commutative Frobenius algebra over the homology $R$ of the empty link. In particular, if $R$ is a field, $A$ must be finite dimensional over $R$, and this property fails in the above examples. 

Other theories~\cite{KR05,Caut17,QRS18,RW_sym20,OR_glmk22} have finite-dimensional homology groups over $R$ (upon normalization, if needed), but their construction is only available for braid closures, making it hard to extend them to link cobordism. 
For an early categorification of the colored Jones polynomial~\cite{Kho_colored_Jones05} 
functoriality is not known either, see also~\cite{BW_Rasmussen08}. 

It is an important open problem to (in the first case) modify these link homology theories, including the Webster homology, to make them functorial under link cobordisms and (in the second) find whether they can redefined in a more functorial way, for all link diagrams and extending to link cobordisms or even to tangle cobordisms. 

\begin{center}
\input{mfig_016}
\end{center}

\subsection{Foam evaluation and link homology}
    L.-H.~Robert and E.~Wagner~\cite{RW20} found a purely combinatorial approach to $\SL(N)$ link homology (again, for minuscule representations $\Lambda^a_q V$). Their approach is based on foam evaluations. Foams are cobordisms between planar graphs $\Gamma$ as above and implicitly appear in most approaches to $\SL(N)$ link homology. Maps in the long exact sequences in Figure~\ref{mfig_009} between the two planar graphs should be induced by foam cobordisms 
  shown in Figure~\ref{mfig_010}   between graphs $\Gamma_0,\Gamma_1$, also see Figure~\ref{mfig_016}. 
    To build a combinatorial theory, Robert and Wagner construct a subtle evaluation of closed foams to symmetric functions $R=\Z[x_1,\ldots, x_N]^{S_N}$ in $N$ variables.

Let us briefly review foams and Robert--Wagner foam evaluation. 
  $\GL(N)$ foams can have facets of thickness $a\in \{1,\ldots, N\}$, seams where facets of thickness $a,b,a+b$ come together and vertices with 6 adjoint facets of thickness $a,b,c,a+b,b+c,a+b+c$, see Figure~\ref{mfig_021}. At first, consider \emph{closed foams} embedded in $\R^3$, that is, foams without boundary, while later we will need foam with boundary, viewed as cobordisms in $\R^2\times [0,1]$ between MOY graphs. 

\begin{center}
\input{mfig_021}
\end{center}

 A \emph{coloring} of foam $F$ is a map 
$c: \mathsf{facets} \lra \mathcal{P}(N)$,  where $|c(f)|$ is the thickness of a facet $f$,  with the flow condition at each seam, $ c(f_3)=c(f_1)\sqcup c(f_2)$, see Figures~\ref{mfig_000},~\ref{mfig_017}. In particular, $c(f_1),c(f_2)$ are disjoint sets. Here $\mathcal{P}(N)$ is the set of subsets of $\{1,\ldots, N\}$, so that a coloring maps facets to subsets of the set of colors from $1$ to $N$. A facet of thickness $a$ is mapped to a subset of cardinality $a$; subsets for $a$, $b$ thickness facets that meet along a seam and become an $a+b$ facet are disjoint and their union is the subset for the $a+b$ facet.

\begin{center}
\input{mfig_000}
\end{center}

\begin{center} 
\input{mfig_017}
\end{center}

Define the \emph{bicolored surface} $F_{ij}(c)$, $1\le i<j\le N$ as the union of facets that contain exactly one color from $\{i,j\}$. For a closed foam $F$, its bicolored surface is a closed compact surface without boundary embedded in $\R^3$. Figures~\ref{mfig_000},~\ref{mfig_017} on the right explain why $F_{ij}(c)$ has no singularities along seams of the foam, always containing either none or two facets along a seam. A similar computation shows that $F_{ij}(c)$ has no singularities at vertices of the foam.  

Thus, $F_{ij}(c)$ is an orientable surface and the union of its connected components, each of Euler characteristic $2-2g$, for a component of genus $g\ge 0$. Robert--Wagner evaluation $\langle F,c\rangle$ of a closed foam $F$ and its coloring $c$ has the form 
\[\langle F,c\rangle = \pm \prod_{i<j} (x_i-x_j)^{-\chi(F_{ij}(c))/2}\prod (\mathsf{facet\ decoration\  contributions}).
\]
We refer to~\cite{RW20} and a review in~\cite{KK_deformation_RW20} for the subtle formula for the minus sign, facet decorations and their contributions. Facets are decorated by dots (observables) labelled by homogeneous symmetric functions $f$ in $a$ variables, for a facet of thickness $a$. A coloring tells one which $a$ variables out of $N$ to select from $x_1,\dots, x_N$ to turn $f$ into a function in $x$'s. 

In general $\langle F,c\rangle$ has denominators $x_i-x_j$ but the sum 
\[
\langle F \rangle \ = \ \sum_c \langle F,c\rangle \in R
\]
is a symmetric polynomial in $x_1,\ldots, x_N$. 

With that sophisticated yet beautiful foam evaluation $\langle F \rangle$ at hand, Robert and Wagner build the state space (or homology) for a planar MOY graph $\Gamma$ using the universal construction. 

Universal construction of topological theories~\cite{BHMV95,Kh20_univ_const_two,KQR,IK_22_linear,IK-top-automata} begins with an evaluation of closed objects (closed foams, in our case) and builds state spaces for generic cross-sections of these objects. A generic cross-section $\Gamma$ of a foam in $\R^3$ by a plane $\R^2$ is a planar MOY graph $\Gamma$. Fix $\Gamma$ and consider the free $R$-module $\mathsf{Fr}(\Gamma)$ with a basis $[F]$ of all foams $F$ with $\partial F = \Gamma$, see Figure~\ref{mfig_X1} on the left. Notice that we no longer restrict to closed foams, instead looking at foams $F$ in $\R^2\times (-\infty,0]$ with $\Gamma$ as the boundary. 

\begin{center}
\input{mfig_X1}
\end{center}

Define an $R$-bilinear form \[(\:\:,\:\:):\mathsf{Fr}(\Gamma)\times \mathsf{Fr}(\Gamma) \lra R\] 
by $([F],[F_1])= \brak{\overline{F}F_1}$ and extending via $R$-bilinearity. Here we glue $F$ and $F_1$ along the common boundary $\Gamma$ into a closed foam $\overline{F}F_1$, where $\overline{F}$ denotes the reflection of $F$ in a horizontal plane. Define the state space 
\[ H(\Gamma) \ := \ \mathsf{Fr}(\Gamma)/\ker((\:\:,\:\:)_\Gamma)\]
of an MOY graph $\Gamma$ as the quotient of the (large) free module $\mathsf{Fr}(\Gamma)$ by the kernel of the bilinear form $(\:\:,\:\:)_\Gamma$. This means that an $R$-linear combination of foams $\sum_i \lambda_i F_i$, with $\lambda_i\in R$ and $\partial F_i = \Gamma$ is $0$ in $H(\Gamma)$ if and only if for any foam $F$ with $\partial F=\Gamma$, we have 
\[
\sum_i \lambda_i \langle \overline{F}F_i \rangle =0 \in R.
\]
One thinks of $\sum_i \lambda_i F_i=0$ as a linear skein relation on foams with boundary $\Gamma$. 

It is easy to see that state spaces $H(\Gamma)$, over all $\Gamma$, form a functorial topological theory. 
Namely, given a foam $F\in \R^2\times [0,1]$ with top and bottom boundary, so that $\partial F= \partial_1 F \sqcup (- \partial_0 F),$ composition with $F$ induces a map from $\mathsf{Fr}(F_0)$ to $\mathsf{Fr}(F_1)$. The map takes the kernel of the bilinear form $(\:\:,\:\:)_{\Gamma_0}$ to the kernel of the bilinear form $(\:\:,\:\:)_{\Gamma_1}$.  
Consequently, there is an induced map on state spaces 
\[ H(F) \ : \ H(\partial_0 F)\lra H(\partial_1 F) .\]
Varying over all foams $F$ with boundary, this results in a functor $H$ from the category of foams to the category of $R$-modules. 

For a general foam evaluation $\langle F\rangle$ such a functor is not interesting. In particular, one wants state spaces $H(\Gamma)$ to be sufficiently small, for instance, be finitely-generated $R$-modules. Even under such an assumption, functor $H$ will not be a TQFT in general, with the natural map  $H(\Gamma_1)\otimes_R H(\Gamma_2) \lra H(\Gamma_1\sqcup \Gamma_2)$ not an isomorphism of $R$-modules. 

For their evaluation sketched above, Robert and Wagner in \cite{RW20} proved: 
    \begin{theorem}[Robert--Wagner] $H$ is 
   a functorial TQFT from \emph{Foams} to free graded $R$-modules, with graded ranks the Murakami--Ohtsuki--Yamada planar graph invariants:  
   \[ 
   \mathsf{grank}_R(\Gamma) = P_N(\Gamma). 
   \]
  
   \end{theorem}

Methods of Yonezawa--Wu~\cite{Yon11,Wu12} 
for constructing homology from matrix factorizations for MOY graphs with edges of arbitrary thickness apply to Robert--Wagner state spaces. Robert--Wagner foam TQFT gives rise to a link homology theory which categorifies Reshetikhin--Turaev invariant for $\SL(N)$ and link components labelled by exterior powers of the fundamental representation~\cite{ETW18}. 

Furthermore, in~\cite{ETW18} this result is extended to:
    \begin{theorem}[M.~Ehrig--D.~Tubbenhauer--P.~Wedrich] Robert--Wagner link homology theory is functorial for link cobordisms. 
    \end{theorem}

 Robert--Wagner's approach~\cite{RW20,ETW18} to categorification of Reshetikhin--Turaev $\SL(N)$ invariants is \emph{complementary} to all the others, which require \emph{specific categories} (highest weight categories, coherent sheaves or Fukaya--Floer categories on quiver varieties, particular representation theory categories). In the Robert--Wagner construction, categories appear at the last step only (when extending to graphs with boundary, tangles, and their cobordisms). 
 It is a \emph{categorification of the state sum approach to quantum invariants} and should also be relevant to some models of 3-dimensional statistical mechanics. 

\vspace{0.07in} 

 Prior to Robert--Wagner's work, foams were heavily used in link homology and categorified quantum groups, see for example~\cite{MSV_foam09,Mack09,QR_comb16,Cap13,RW_deformation16,Wed19}. Foam evaluation should help to streamline and clarify a significant amount of prior work on the subject, including replacing the ground ring $\Q$ by $\Z$ or by the ring of symmetric functions $\Z[x_1,\ldots, x_N]^{S_N}$.

%
%

\section{Interactions and applications}

\subsection{Some interactions}

\quad 


{\it Soergel category.}
A version of Robert--Wagner evaluation~\cite{RW_sym20} allows one to describe the category of Soergel bimodules~\cite{EMTW20} for the symmetric group $S_N$, with homs given by foams in $\R^3$ between braid-like graphs in the plane modulo universal construction relations~\cite{RW_sym20,KRW_inprogress}. The Soergel category is central in geometric representation theory and categorifies the Hecke algebra of the symmetric group. Originally, the relation between foams and the Soergel category was established by D.~Rose and P.~Wedrich~\cite{RW_deformation16} and P.~Wedrich~\cite{Wed19}. An earlier diagrammatic approach to the Soergel category by Elias--Khovanov~\cite{EK10} is via a two-dimensional graphical calculus rather than foams in $\R^3$, with the missing dimension encoded by labels in the ordered set $\{1,\dots, N-1\}$. 

 \vspace{0.1in}
{\it Kronheimer--Mrowka theory.}
Inspired by the Robert--Wagner foam evaluation, Robert and one of us~\cite{KR21} related \emph{unoriented} $\SL(3)$ foams and Kronheimer--Mrowka gauge $SO(3)$ theory for 3-orbifolds~\cite{KM_Tait19}, proving one of Kronheimer--Mrowka conjectures. Kronheimer--Mrowka theory, in a rather special case, assigns homology groups to planar trivalent graphs $\Gamma$, and similar groups can be recovered from an unoriented version of Robert--Wagner foam  evaluation. (Kronheimer--Mrowka theory is much more general, and, in particular, assigns homology groups to spacial trivalent graphs, for which a combinatorial counterpart is unknown.) Kronheimer--Mrowka theory and its combinatorial counterpart in~\cite{KR21} for planar trivalent graphs relate  to the 4-Color Theorem, aiming to prove and rethink the  4-Color Theorem in a more conceptual way and via its relations to TQFTs, gauge theory, and low-dimensional topology. These two theories have been investigated by D.~Boozer~\cite{Booz19,Booz23}, who have shown, in particular, that they give rise to two non-isomorphic functors on foams in $\R^3$.  

\vspace{0.1in} 

{\it APS homology.}
Foam evaluation allows for a natural extension of the Asaeda--Przytycki--Sikora \cite{APS_Kauffman04} annular $\SL(2)$ homology to the 
equivariant setting of a larger ground ring as well as its extension to annular $\SL(N)$ homology~\cite{AkKh21,Akh23}.  

\subsection{\texorpdfstring{$\SL(2)$}{SL2} and \texorpdfstring{$\SL(3)$}{SL3} homology theories}
\label{subsection:SL2_SL3}

$\SL(2)$ homology (aka Khovanov homology) is noticeably simpler than $\SL(N)$ homology and was discovered first~\cite{Kho_Duke00}. Foams are replaced by surfaces,  and homology of links is built from a 2-dimensional TQFT where homology of the unknot has rank two over homology of the empty link. Khovanov homology categorifies the Jones polynomial. Odd Khovanov homology, discovered by P.~Ozsvath, J.~Rasmussen and Z.~Szabo~\cite{ORS13}, is another bigraded categorification of the Jones polynomial. 

\vspace{0.10in}

$\SL(3)$ homology, which categorifies the Kuperberg quantum $\mathfrak{sl}(3)$ invariant~\cite{Kup96}, is in-between $N=2$ and $N\ge 4$ cases complexity-wise~\cite{Kho04,MV_universalsl307,MN_su308,Clark_functoriality09}. $\SL(3)$ foams do not have vertices, and closed $\SL(3)$ foams can be evaluated via ``localization along singular circles'' and manipulation of foams with a single such circle. Single-circle foam evaluation is encoded in the cohomology groups of the flag variety of $\C^3$ (commutative algebra structure plus the trace map).


\subsection{Two applications to 4D topology}

\quad 

{\bf I.} {\it Rasmussen invariant.} 
The Rasmussen invariant $s(K)$ comes from an equivariant version of Khovanov homology ($N=2$) given by replacing the ground ring $\Z$ by $\Q[t]$ and the Frobenius algebra
$A=\Z[X]/(X^2)$  of homology of the unknot by $A_t=\Q[t,X]/(X^2-t)$ over $\Q[t]$. Then $\SL(2)$-homology $H_t(K)$ of a knot $K$ becomes a $\Q[t]$-module, and it can be written as a sum of its $t$-torsion and a free summand, where the free summand is a free rank one $A_t$-module in homological degree $0$ and $q$-degree $s(K)-1$: 
\[
H_t(K) \ \cong \ \mathsf{Tor}(H_t(K))\oplus A_t\{s(K)-1\}.
\]
\hspace{0.03in}

From \cite{Ras_slicegenus10}, one obtains the following: 
\begin{theorem}[J.~Rasmussen]
  The value $s(K)$ is an invariant of knot concordance and  gives a lower bound on the slice genus of $K$.
  This bound is explicitly computable and sharp on positive knots.
\end{theorem}
The Rasmussen invariant provides a combinatorial proof of Kronheimer--Mrowka--Milnor theorem (formerly
Milnor's conjecture) that the slice genus of $(p,q)$-torus knot is $\frac{(p-1)(q-1)}{2}$.

\vspace{0.10in}

Replacing $\Q$ by $\mathbb{F}_p$ or extending to
$\SL(N)$ homology, $N\ge 3$ and its deformations leads to families of  Rasmussen-like concordance invariants~\cite{MTV07,MTV_rasmusseninv13,Lobb_perturbation12,Lewark_spectral14}. Some linear independence results on these invariants are known and there is even a postcard on the topic \cite{postcard_21}, but the general theory of such concordance invariants remains a mystery. 

\vspace{0.1in} 

{\bf II.} {\it Exotic surfaces in $\R^4$.}
K.~Haiden and I.~Sundberg in~\cite{HS_exotic21} give examples of surfaces $S_1,S_2\subset \mathbb{D}^4$ that bound a knot $K\subset \mathbb{S}^3$
such that $S_1,S_2$ are homeomorphic but not diffeomorphic rel boundary. Non-diffeomorphic property follows by checking that maps induced by $S_1,S_2$ on Khovanov homology 
\[ H(S_i) \ : \ H(K) \lra H(\mathsf{empty} \ \mathsf{link})\cong\Z \]
are not equal. 

\vspace{0.07in} 

S.~Akbulut showed in \cite{Akb91} that $S_1,S_2$ have this property by using Donaldson theory.
Thirty years later, there is a combinatorial proof of his result. 

\vspace{0.10in} 

{\it Link homology relates smooth 4D topology and geometric representation theory in a highly intricate way. Braid groups and tangle categories act by exact functors on key categories in geometric representation theory. Tangle cobordisms act by natural transformations between these functors. 
These tangle cobordism invariants are subtle enough to distinguish between distinct smooth structures on properly embedded surfaces in the four-ball $\mathbb{D}^4$ with the same underlying topological structure. 
}

\vspace{0.10in} 

Going beyond mere homology groups, we should mention the spectrification of Khovanov homology by R.~Lipshitz and S.~Sarkar~\cite{LS_stable14,LLS_Burnside20} and~\cite{HKK16}; see N.~Kitchloo \cite{Kitch_I19} for spectrifications of $\SL(N)$ homology.  

\vspace{0.05in} 

Heegaard--Floer theory~\cite{OS06_lectures,OS06_intro} 
and link Floer homology of Ozsvath--Szabo and Rasmussen~\cite{Manol16,OS_Floer18} play a fundamental role in modern low-dimensional topology and should relate to yet-unknown categorified quantum group $\GL(1|1)$ and its categorified representations, as studied by Y.~Tian, A.~Ellis--I.~Petkova--V.~Vertesi, A.~Manion--R.~Rouquier and others~\cite{Tian_UT14,EPV_tangle19,MR_Hegaard20}.

\vspace{0.05in} 

This informal write-up discusses only a fraction of the rich developments in the past twenty-five years that were inspired by the work~\cite{CF94} of Igor Frenkel and Louis Crane and other related ideas of Igor Frenkel. We are delighted to celebrate Igor Frenkel's anniversary and wish him many more years of inspiring and exciting mathematical discoveries.


\bibliographystyle{amsalpha} 
\bibliography{z_brauer-group}

\end{document}

%% file: figure-0.1.tex
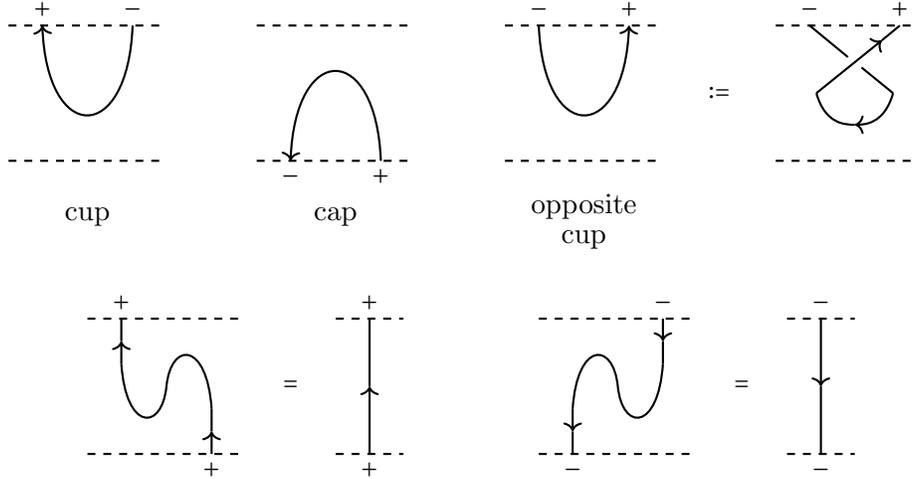
\begin{figure}
    \centering
    \begin{tikzpicture}[scale=0.6]

    \begin{scope}[shift={(0,0)}]
    \draw[thick,dashed] (-0.25,3) -- (3.25,3);
    \draw[thick,dashed] (-0.25,0) -- (3.25,0);
    
    \node at (0.5,3.35) {$+$};
    \node at (2.5,3.35) {$-$};
    \draw[thick,<-] (0.5,3) .. controls (0.6,0.35) and (2.4,0.35) .. (2.5,3);
    \node at (1.5,-1.25) {cup};
    \end{scope}  

    \begin{scope}[shift={(5.5,0)}]
    \draw[thick,dashed] (-0.25,3) -- (3.25,3);
    \draw[thick,dashed] (-0.25,0) -- (3.25,0);
    
    \node at (0.5,-0.35) {$-$};
    \node at (2.5,-0.35) {$+$};
    
    \draw[thick,<-] (0.5,0) .. controls (0.6,2.65) and (2.4,2.65) .. (2.5,0);
    \node at (1.5,-1.25) {cap};
    \end{scope}

    \begin{scope}[shift={(11,0)}]
    \draw[thick,dashed] (-0.25,3) -- (3.25,3);
    \draw[thick,dashed] (-0.25,0) -- (3.25,0);
    
    \node at (0.5,3.35) {$-$};
    \node at (2.5,3.35) {$+$};
    \draw[thick,->] (0.5,3) .. controls (0.6,0.35) and (2.4,0.35) .. (2.5,3);
    \node at (1.5,-1) {opposite};
    \node at (1.5,-1.75) {cup};
    
    \node at (4.5,1.5) {$:=$};
    \end{scope}

    \begin{scope}[shift={(17,0)}]
    \draw[thick,dashed] (-0.25,3) -- (3.25,3);
    \draw[thick,dashed] (-0.25,0) -- (3.25,0);
    
    \node at (0.5,3.35) {$-$};
    \node at (2.5,3.35) {$+$};
    \draw[thick] (0.65,1.5) .. controls (0.85,0.8) and (1.3,0.8) .. (1.5,0.8);
    \draw[thick,<-] (1.5,0.8) .. controls (1.7,0.8) and (2.15,0.8) .. (2.35,1.5); 
    \draw[thick,->] (0.65,1.5) -- (2.1,2.67);
    \draw[thick] (2.1,2.67) -- (2.5,3);
    

    \draw[thick] (0.5,3) -- (1.35,2.3);
    \draw[thick] (1.65,2.07) -- (2.35,1.5);
    \end{scope}
    \begin{scope}[shift={(1.75,-6.5)}]
    \draw[thick,dashed] (-0.25,3) -- (3.25,3);
    \draw[thick,dashed] (-0.25,0) -- (3.25,0);
    
    \node at (0.5, 3.35) {$+$};
    \node at (2.5,-0.35) {$+$};
    
    \draw[thick] (0.5,3) -- (0.5,2.5);
    \draw[thick,<-] (0.5,2.5) -- (0.5,2.0);    
    \draw[thick,<-] (2.5,0.5) -- (2.5,0);
    \draw[thick] (2.5,1.0) -- (2.5,0.5); 
    
    \draw[thick] (0.5,2) .. controls (0.6,0.5) and (1.4,0.5) .. (1.5,1.5);
    \draw[thick] (1.5,1.5) .. controls (1.6,2.5) and (2.4,2.5) .. (2.5,1);    
    \node at (4.25,1.5) {$=$};
    \end{scope}
    
    \begin{scope}[shift={(7.25,-6.5)}]
    \draw[thick,dashed] (-0.25,3) -- (1.25,3);
    \draw[thick,dashed] (-0.25,0) -- (1.25,0);
    \node at (0.5, 3.35) {$+$};
    \node at (0.5,-0.35) {$+$};
    \draw[thick,<-] (0.5,1.5) -- (0.5,0);
    \draw[thick] (0.5,3) -- (0.5,1.5);
    \end{scope}
    
    \begin{scope}[shift={(11.75,-6.5)}]
    \draw[thick,dashed] (-0.25,3) -- (3.25,3);
    \draw[thick,dashed] (-0.25,0) -- (3.25,0);
    
    \node at (2.5, 3.35) {$-$};
    \node at (0.5,-0.35) {$-$};
    
    \draw[thick,->] (2.5,3) -- (2.5,2.5);
    \draw[thick] (2.5,2.5) -- (2.5,2.0);    
    \draw[thick] (0.5,0.5) -- (0.5,0);
    \draw[thick,->] (0.5,1.0) -- (0.5,0.5); 

    \draw[thick] (0.5,1) .. controls (0.6,2.5) and (1.4,2.5) .. (1.5,1.5);
    \draw[thick] (1.5,1.5) .. controls (1.6,0.5) and (2.4,0.5) .. (2.5,2);   
    \node at (4.25,1.5) {$=$};
    \end{scope}
    
    \begin{scope}[shift={(17.25,-6.5)}]
    \draw[thick,dashed] (-0.25,3) -- (1.25,3);
    \draw[thick,dashed] (-0.25,0) -- (1.25,0);
    \node at (0.5, 3.35) {$-$};
    \node at (0.5,-0.35) {$-$};
    \draw[thick] (0.5,1.5) -- (0.5,0);
    \draw[thick,->] (0.5,3) -- (0.5,1.5);
    \end{scope}

    \end{tikzpicture}
    \caption{Top left: the cup and cap cobordisms. The cup and cap cobordisms for the opposite orientation are obtained by composing these two with the permutation cobordisms, see top right for the cup cobordism with the opposite orientation. Bottom row shows isotopy relations on the cup and cap cobordisms.}
    \label{figure-0.1}
\end{figure}

%% file: figure-A1.tex
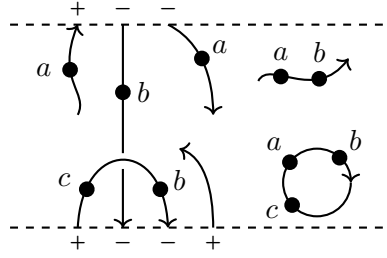
\begin{figure}
    \centering
\begin{tikzpicture}[scale=0.6]
\begin{scope}[shift={(0,0)}]
\draw[thick,dashed] (-1,4.5) -- (7.5,4.5);
\draw[thick,dashed] (-1,0) -- (7.5,0);

\node at (0.5,-0.35) {$+$};
\node at (1.5,-0.35) {$-$};
\node at (2.5,-0.35) {$-$};
\node at (3.5,-0.35) {$+$};

\node at (0.5,4.85) {$+$};
\node at (1.5,4.85) {$-$};
\node at (2.5,4.85) {$-$};

\draw[thick,->] (0.5,0) .. controls (0.6,2) and (2.4,2) .. (2.5,0);
\draw[thick,fill] (0.85,0.85) arc (0:360:1.5mm);
\node at (0.20,1.05) {$c$};
\draw[thick,fill] (2.48,0.85) arc (0:360:1.5mm);
\node at (2.75,1.05) {$b$};

\draw[thick,<-] (0.5,4.5) .. controls (0,3.25) and (0.75,3.0) .. (0.5,2.5);
\draw[thick,fill] (0.48,3.5) arc (0:360:1.5mm); 
\node at (-0.25,3.5) {$a$};

\draw[thick] (1.5,4.5) -- (1.5,1.65);
\draw[thick,->] (1.5,1.3) -- (1.5,0);
\draw[thick,fill] (1.65,3.0) arc (0:360:1.5mm); 
\node at (1.95,3.0) {$b$};

\draw[thick,->] (2.5,4.5) .. controls (3.5,4) and (3.5,2.75) .. (3.5,2.5);
\draw[thick,fill] (3.40,3.75) arc (0:360:1.5mm); 
\node at (3.65,4.0) {$a$};

\draw[thick,<-] (2.75,1.75) .. controls (3.5,1.50) and (3.5,0.25) .. (3.5,0);

\draw[thick,->] (4.5,3.25) .. controls (4.75,3.75) and (5.75,2.75) .. (6.50,3.75);
\draw[thick,fill] (5.15,3.35) arc (0:360:1.5mm);
\node at (5,3.85) {$a$};
\draw[thick,fill] (6.00,3.30) arc (0:360:1.5mm);
\node at (5.85,3.90) {$b$};

\draw[thick,<-] (6.55,1) arc (0:360:0.75); 
\draw[thick,fill] (6.45,1.55) arc (0:360:1.5mm);
\node at (6.65,1.95) {$b$};
\draw[thick,fill] (5.35,1.45) arc (0:360:1.5mm);
\node at (4.85,1.85) {$a$};
\draw[thick,fill] (5.40,0.50) arc (0:360:1.5mm);
\node at (4.80,0.35) {$c$};

\end{scope}

\end{tikzpicture}
    \caption{A morphism from oriented 1-manifold $(+--+)$ to $(+--)$ in category $\Cob_{\Sigma,1}$. This morphism has seven outer and five inner boundary points. It has two arcs, three half-intervals, one floating interval and one circle.}
    \label{figure-A1}
\end{figure}

%% file: figure-3.tex
\begin{figure}
    \centering
\begin{tikzpicture}[scale=0.6]
\begin{scope}[shift={(0,0)}]
\node at (0.5, 3.35) {$+$};
\node at (0.5,-0.35) {$+$};
\draw[thick,dashed] (0,3) -- (1,3);
\draw[thick,dashed] (0,0) -- (1,0);
\draw[thick,->] (0.5,0) -- (0.5,3); 
\draw[thick,fill] (0.65,1.5) arc (0:360:1.5mm);
\node at (1.1,1.5) {$a$};

\node at (2.25,3) {$V$};
\draw[thick,<-] (2.25,2.5) -- (2.25,0.5);
\node at (3,1.5) {$m_a$};
\node at (2.25,0) {$V$};
\end{scope}

\begin{scope}[shift={(6.5,0)}]
\node at (0.5, 3.35) {$-$};
\node at (0.5,-0.35) {$-$};
\draw[thick,dashed] (0,3) -- (1,3);
\draw[thick,dashed] (0,0) -- (1,0);
\draw[thick,<-] (0.5,0) -- (0.5,3); 
\draw[thick,fill] (0.65,1.5) arc (0:360:1.5mm);
\node at (1.1,1.5) {$a$};

\node at (2.35,3) {$V^*$};
\draw[thick,<-] (2.25,2.5) -- (2.25,0.5);
\node at (3,1.5) {$m_a^*$};
\node at (2.35,0) {$V^*$};
\end{scope}

\begin{scope}[shift={(13,0)}]
\node at (0.5, 3.35) {$+$};
\node at (0.5,-0.35) {$+$};
\draw[thick,dashed] (0,3) -- (1,3);
\draw[thick,dashed] (0,0) -- (1,0);
\draw[thick,->] (0.5,0) -- (0.5,3); 
\draw[thick,fill] (0.65,1.5) arc (0:360:1.5mm);
\node at (1.0,1.5) {$\omega$};

\node at (2,1.5) {$=$};

\node at (3.25, 3.35) {$+$};
\node at (3.25,-0.35) {$+$};
\draw[thick,dashed] (2.75,3) -- (3.75,3);
\draw[thick,dashed] (2.75,0) -- (3.75,0);
\draw[thick,->] (3.25,0) -- (3.25,3);

\draw[thick,fill] (3.40,2.35) arc (0:360:1.5mm);
\node at (3.85,2.35) {$a_1$};
\draw[thick,fill] (3.40,1.65) arc (0:360:1.5mm);
\node at (3.85,1.65) {$a_2$};
\node at (3.85,1.0) {$\vdots$};
\draw[thick,fill] (3.40,0.40) arc (0:360:1.5mm);
\node at (3.85,0.40) {$a_n$};

\node at (5.25,3) {$V$};
\draw[thick,<-] (5.25,2.5) -- (5.25,0.5);
\node at (7.80,1.45) {$m_{\omega} = m_{a_1}\cdots \:m_{a_n}$};
\node at (5.25,0) {$V$};
\end{scope}

\end{tikzpicture}
    \caption{Left: to a labelled dot on an upward-oriented interval, functor $\mcF$ associates an endomorphism $m_a$ of $V$. The dual endomorphism of $V^{\ast}$ (given by the transposed matrix of that of $m_a$) is associated by $\mcF$ to the downward-oriented interval with an $a$-dot. Right: a sequence of labelled dots on an upward-oriented interval defines a word $\omega=a_1\cdots a_n$ and the induced endomorphism $m_{\omega}=m_{a_1}\cdots m_{a_n}$ of $V$.}
    \label{figure-3}
\end{figure}
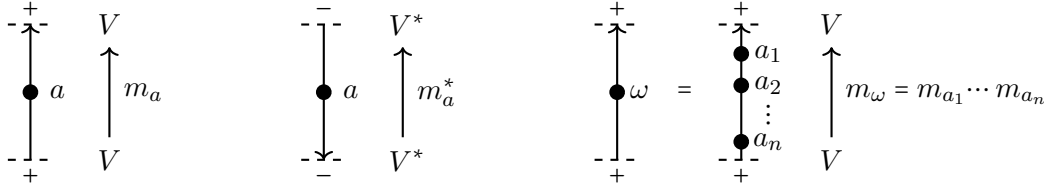

%% file: figure-2.tex
\begin{figure}
    \centering
\begin{tikzpicture}[scale=0.6]

\begin{scope}[shift={(0,0)}]
\node at (0.5,3.35) {$+$};
\draw[thick,dashed] (0,3) -- (1,3);
\draw[thick,dashed] (0,0) -- (1,0);
\draw[thick,<-] (0.5,3) -- (0.5,1.5); 

\node at (2,3) {$V$};
\draw[thick,<-] (2,2.5) -- (2,0.5);
\node at (2,0) {$\kk$};

\node at (2.8,3) {$\ni$};
\node at (3.5,3) {$v_0$};
\draw[thick,<-|] (3.5,2.5) -- (3.5,0.5);
\node at (2.8,0) {$\ni$};
\node at (3.5,0) {$1$};
\end{scope}

\begin{scope}[shift={(8,0)}]
\node at (0.5,-0.35) {$+$};
\draw[thick,dashed] (0,3) -- (1,3);
\draw[thick,dashed] (0,0) -- (1,0);
\draw[thick,->] (0.5,0) -- (0.5,1.5); 

\node at (2,3) {$\kk$};
\draw[thick,<-] (2,2.5) -- (2,0.5);
\node at (2.60,1.5) {$v^{\ast}$};
\node at (2,0) {$V$};

\node at (2.8,3) {$\ni$};
\node at (3.0,0) {$\ni$};

\node at (4,3) {$v^{\ast}(v)$};
\draw[thick,<-|] (4,2.5) -- (4,0.5);
\node at (4,0) {$v$};
\end{scope}

\begin{scope}[shift={(17,0)}]
\draw[thick,dashed] (0,3) -- (1,3);
\draw[thick,dashed] (0,0) -- (1,0);
\draw[thick,<-] (0.5,2.35) -- (0.5,0.65);
\node at (1.5,1.5) {$=$};
\node at (3.2,1.5) {$v^{\ast}(v_0)\in\kk$};
\end{scope}

\end{tikzpicture}
    \caption{Left: two half-intervals and associated maps. Right: evaluation of a floating interval (without dot decorations).}
    \label{figure-2}
\end{figure}
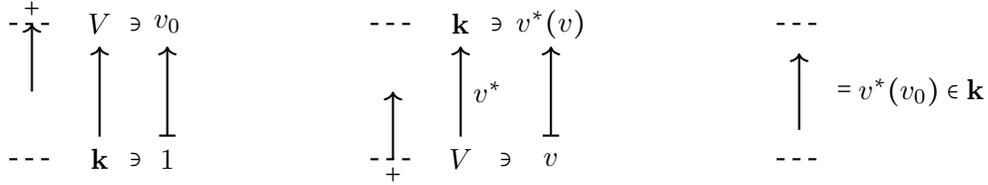

%% file: figure-4.tex
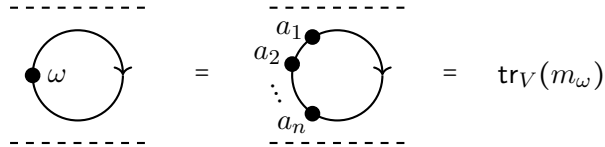
\begin{figure}
    \centering
\begin{tikzpicture}[scale=0.6]

\begin{scope}[shift={(12,0)}]
\draw[thick,dashed] (0,3) -- (3,3);
\draw[thick,dashed] (0,0) -- (3,0);
\draw[thick,<-] (2.5,1.5) arc (0:360:1);
\draw[thick,fill] (0.65,1.5) arc (0:360:1.5mm);
\node at (1.05,1.5) {$\omega$};

\node at (4.25,1.5) {$=$};

\begin{scope}[shift={(5.75,0)}]
\draw[thick,dashed] (0,3) -- (3,3);
\draw[thick,dashed] (0,0) -- (3,0);
\draw[thick,<-] (2.5,1.5) arc (0:360:1);

\draw[thick,fill] (1.10,2.35) arc (0:360:1.5mm);
\node at (0.45,2.50) {$a_1$};
\draw[thick,fill] (0.65,1.75) arc (0:360:1.5mm);
\node at (-0.05,1.95) {$a_2$};
\node at (0.15,1.1) {\rotatebox[origin=c]{-20}{$\ddots$}};
\draw[thick,fill] (1.10,0.65) arc (0:360:1.5mm);
\node at (0.50,0.40) {$a_n$};

\node at (4.0,1.5) {$=$};
\node at (6.25,1.5) {$\tr_V(m_{\omega})$};
\end{scope}

\end{scope}

\end{tikzpicture}
    \caption{ An $\omega$-decorated circle evaluates to the trace of $m_{\omega}$ on $V$.}
    \label{figure-4}
\end{figure}

%% file: figure-5.tex
\begin{figure}
    \centering
\begin{tikzpicture}[scale=0.6]
\begin{scope}[shift={(-1,0)}]
\draw[thick,dashed] (-1,3) -- (4,3);
\draw[thick,dashed] (-1,0) -- (4,0);

\draw[thick,<-] (-0.5,1.5) .. controls (0.2,3) and (2.8,0) .. (3.5,1.5);
\draw[thick,fill] (1.65,1.5) arc (0:360:1.5mm);
\node at (1.5,0.9) {$\omega$};

\node at (5.5,1.5) {$=$};

\begin{scope}[shift={(8,0)}]
\draw[thick,dashed] (-1,3) -- (4,3);
\draw[thick,dashed] (-1,0) -- (4,0);

\draw[thick,<-] (-0.5,1.5) .. controls (0.2,3) and (2.8,0) .. (3.5,1.5);

\draw[thick,fill] (0.70,1.90) arc (0:360:1.5mm);
\node at (0.55,1.35) {$a_1$};
\node at (1.50,0.9) {\rotatebox[origin=c]{-25}{$\cdots$}};
\draw[thick,fill] (2.65,1.10) arc (0:360:1.5mm);
\node at (2.50,0.55) {$a_n$};

\node at (5.5,1.5) {$=$};

\node at (8.00,1.5) {${v}^{\ast}(m_{\omega}v)$};

\end{scope}

\end{scope}

\end{tikzpicture}
    \caption{Evaluation of a floating $\omega$-decorated interval (both endpoints are inner). }
    \label{figure-5}
\end{figure}
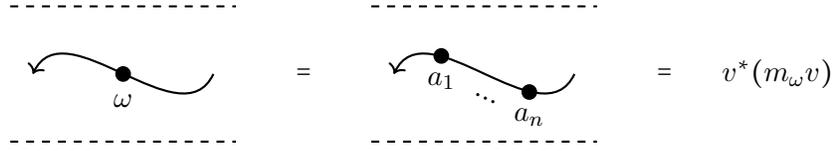

%% file: fig_A1.tex
\begin{figure}
    \centering
\begin{tikzpicture}[scale=0.6,decoration={
    markings,
    mark=at position 0.5 with {\arrow{>}}}]

\draw[thick,fill] (0.15,1) arc (0:360:1.5mm);
\node at (0.15,0.25) {$q_{\init}$};
\node at (-0.5,0.25) {$\ni$};
\node at (-1.25,0.25) {$Q_{\init}$};

\draw[thick,postaction={decorate}] (0,1) -- (2,3);

\draw[thick,fill] (2.15,3) arc (0:360:1.5mm);

\node at (0.4,2.4)  {$a_1$};

\draw[thick,postaction={decorate}] (2,3) -- (4,1);
\node at (5.65,1.75) {$a_3$};

\draw[thick,fill] (4.15,1) arc (0:360:1.5mm);

\node at (3.5,2.4) {$a_2$};

\draw[thick,postaction={decorate}] (4,1) -- (6,3);
\draw[thick,fill] (6.15,3) arc (0:360:1.5mm);

\node at (7,2.2) {\rotatebox[origin=c]{-35}{$\cdots$}};
\node at (8,1.5) {\rotatebox[origin=c]{-35}{$\cdots$}};

\draw[thick,fill] (9.15,1) arc (0:360:1.5mm);

\draw[thick, postaction={decorate}] (9,1) -- (11,3); 

\draw[thick,fill] (11.15,3) arc (0:360:1.5mm);

\node at (11.60,3) {$q_{\t}$};
\node at (12.25,3) {$\in$};
\node at (13,3) {$Q_{\t}$};

\node at (10.65,1.75) {$a_n$};

\end{tikzpicture}
    \caption{Oriented path $a_1a_2\cdots a_n$ in the graph of an automaton, starting at an initial vertex and terminating at an accepting vertex.}
    \label{fig_A1}
\end{figure}
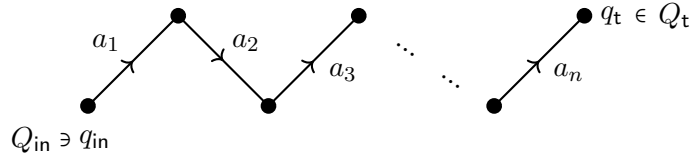

%% file: figure-X1.tex
\begin{figure}
    \centering
\begin{tikzpicture}[>=stealth',shorten >=1pt,auto,node distance=2cm]

  \node[thick,initial,state, initial text=] (x) {$q_1$};
  \node[thick,state] (y) [right of=x]  {$q_2$};
  \node[thick,state,accepting] (z) [right of=y] {$q_3$};

  \path[thick,->] (x)  edge [loop above] node {$a$} (x)
             edge  node {$b\:$} (y)
        (y) edge [loop above] node {$b$} (y)
             edge [bend right=-15] node {$a,b$} (z)
        (z) edge [bend left=40] node {$a\:\:$} (x)
             edge [bend left=15] node {$b$} (y);
\end{tikzpicture} 
\qquad 

    \caption{A nondeterministic automaton on 3 states $q_1$, $q_2$, and $q_3$ that accepts the language $L=(a+b)^{\ast}b(a+b)$ in Example~\ref{ex_language}. It has a single initial state $q_1$, indicated by a short arrow into the state, and a single accepting state $q_3$, indicated by the double border.   }
    \label{figure-X1}
\end{figure}
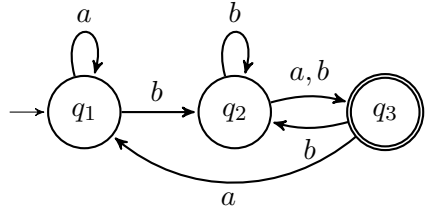

%% file: figure-X2.tex
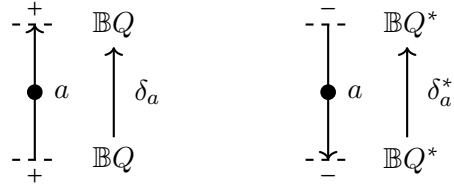
\begin{figure}
    \centering
\begin{tikzpicture}[scale=0.6]
\begin{scope}[shift={(0,0)}]
\node at (0.5, 3.35) {$+$};
\node at (0.5,-0.35) {$+$};
\draw[thick,dashed] (0,3) -- (1,3);
\draw[thick,dashed] (0,0) -- (1,0);
\draw[thick,->] (0.5,0) -- (0.5,3); 
\draw[thick,fill] (0.65,1.5) arc (0:360:1.5mm);
\node at (1.1,1.5) {$a$};

\node at (2.25,3) {$\Bool Q$};
\draw[thick,<-] (2.25,2.5) -- (2.25,0.5);
\node at (3,1.5) {$\delta_a$};
\node at (2.25,0) {$\Bool Q$};
\end{scope}

\begin{scope}[shift={(6.5,0)}]
\node at (0.5, 3.35) {$-$};
\node at (0.5,-0.35) {$-$};
\draw[thick,dashed] (0,3) -- (1,3);
\draw[thick,dashed] (0,0) -- (1,0);
\draw[thick,<-] (0.5,0) -- (0.5,3); 
\draw[thick,fill] (0.65,1.5) arc (0:360:1.5mm);
\node at (1.1,1.5) {$a$};

\node at (2.35,3) {$\Bool Q^*$};
\draw[thick,<-] (2.25,2.5) -- (2.25,0.5);
\node at (3,1.5) {$\delta_a^*$};
\node at (2.35,0) {$\Bool Q^*$};
\end{scope}

\end{tikzpicture}
    \caption{To a labelled dot on an upward-oriented interval functor $\mcF$ associates endomorphism $\delta_a$ of $\Bool Q$. To a dot on a downward-oriented interval functor $\mcF$  associates the dual operator $\delta_a^{\ast}$ on $\Bool Q^{\ast}$.}
    \label{figure-X2}
\end{figure}

%% file: figure-1-1.tex
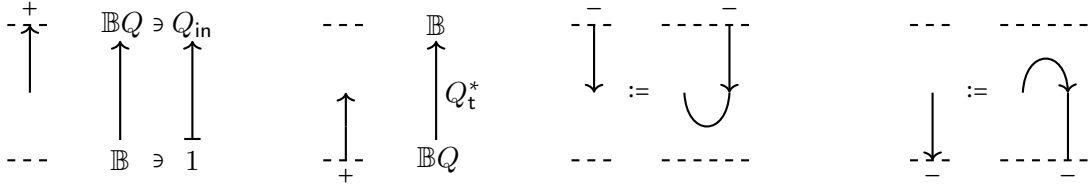
\begin{figure}
    \centering
\begin{tikzpicture}[scale=0.6]

\begin{scope}[shift={(0,0)}]
\draw[thick,dashed] (0,3) -- (1,3);
\draw[thick,dashed] (0,0) -- (1,0);
\node at (0.5,3.35) {$+$};
\draw[thick,->] (0.5,1.5) -- (0.5,3);
\end{scope}

\begin{scope}[shift={(0.8,0)}]
\node at (1.75,3) {$\Bool Q$};
\draw[thick,->] (1.7,0.45) -- (1.7,2.65);
\node at (1.7,0) {$\Bool$};

\node at (3.30,3) {$Q_{\init}$};
\node at (2.55,3) {$\ni$};
\draw[thick,|->] (3.30,0.45) -- (3.30,2.65);
\node at (2.55,0) {$\ni$};
\node at (3.30,0) {$1$};
\end{scope}

\begin{scope}[shift={(7,0)}]
\draw[thick,dashed] (0,3) -- (1,3);
\draw[thick,dashed] (0,0) -- (1,0);
\node at (0.5,-0.35) {$+$};
\draw[thick] (0.5,0) -- (0.5,0.75);
\draw[thick,->] (0.5,0.75) -- (0.5,1.5);
\end{scope}

\begin{scope}[shift={(7.8,0)}]
\node at (1.7,3) {$\Bool$};
\draw[thick,->] (1.7,0.45) -- (1.7,2.65);
\node at (1.75,0) {$\Bool Q$};
\node at (2.30,1.45) {$Q_{\t}^{\ast}$};

\end{scope}

\begin{scope}[shift={(12.5,0)}]
\draw[thick,dashed] (0,3) -- (1,3);
\draw[thick,dashed] (0,0) -- (1,0);
\node at (0.5,3.35) {$-$};
\draw[thick,<-] (0.5,1.5) -- (0.5,3);
\node at (1.5,1.5) {$:=$};

\begin{scope}[shift={(2,0)}]
\draw[thick,dashed] (0,3) -- (2,3);
\draw[thick,dashed] (0,0) -- (2,0);
\node at (1.5,3.35) {$-$};
\draw[thick,<-] (1.5,1.5) -- (1.5,3);

\draw[thick] (0.5,1.5) .. controls (0.55,0.5) and (1.45,0.5) .. (1.5,1.5);

\end{scope}
\end{scope}

\begin{scope}[shift={(20,0)}]
\draw[thick,dashed] (0,3) -- (1,3);
\draw[thick,dashed] (0,0) -- (1,0);
\node at (0.5,-0.35) {$-$};
\draw[thick,->] (0.5,1.5) -- (0.5,0);
\node at (1.5,1.5) {$:=$};

\begin{scope}[shift={(2,0)}]
\draw[thick,dashed] (0,3) -- (2,3);
\draw[thick,dashed] (0,0) -- (2,0);
\node at (1.5,-0.35) {$-$};
\draw[thick] (1.5,1.5) -- (1.5,0);

\draw[thick,->] (0.5,1.5) .. controls (0.55,2.5) and (1.45,2.5) .. (1.5,1.5);

\end{scope}
\end{scope}

\end{tikzpicture}
    \caption{Left: maps assigned to the half-intervals with a $+$ boundary points. Right: defining maps for half-intervals with a $-$ boundary point.
    }
    \label{figure-1-1}
\end{figure}

%% file: figure-5a.tex
\begin{figure}
    \centering
\begin{tikzpicture}[scale=0.6]
\begin{scope}[shift={(-1,0)}]
\draw[thick,dashed] (-1,3) -- (4,3);
\draw[thick,dashed] (-1,0) -- (4,0);

\draw[thick,<-] (-0.5,1.5) .. controls (0.2,3) and (2.8,0) .. (3.5,1.5);
\draw[thick,fill] (1.65,1.5) arc (0:360:1.5mm);
\node at (1.5,0.9) {$\omega$};

\node at (5.5,1.5) {$=$};

\begin{scope}[shift={(8,0)}]
\draw[thick,dashed] (-1,3) -- (4,3);
\draw[thick,dashed] (-1,0) -- (4,0);

\draw[thick,<-] (-0.5,1.5) .. controls (0.2,3) and (2.8,0) .. (3.5,1.5);

\draw[thick,fill] (0.70,1.90) arc (0:360:1.5mm);
\node at (0.55,1.35) {$a_n$};
\node at (1.50,0.9) {\rotatebox[origin=c]{-25}{$\cdots$}};
\draw[thick,fill] (2.65,1.10) arc (0:360:1.5mm);
\node at (2.50,0.55) {$a_1$};

\node at (5.5,1.5) {$=$};

\node at (9.0,1.5) {$\begin{cases} 
1 & \mathsf{if }\ \omega\in L_{(Q)}, \\
0 & \mathsf{otherwise}\end{cases}$};

\end{scope}

\end{scope}

\end{tikzpicture}
    \caption{Evaluation of a floating $\omega$-decorated interval (both endpoints are inner). Note that word $\omega$ is written in the opposite direction from that in the linear case (over a field) in Figure~\ref{figure-5}. This is due to the opposite conventions, where in linear algebra the actions are usually on the left, while in automata theory the actions are on the right.}
    \label{figure-5a}
\end{figure}

%% file: figure-4a.tex
\begin{figure}
    \centering
\begin{tikzpicture}[scale=0.6]

\begin{scope}[shift={(12,0)}]
\draw[thick,dashed] (-0.5,3) -- (3.5,3);
\draw[thick,dashed] (-0.5,0) -- (3.5,0);
\draw[thick,<-] (2.5,1.5) arc (0:360:1);
\draw[thick,fill] (0.65,1.5) arc (0:360:1.5mm);
\node at (1.05,1.5) {$\omega$};

\node at (5.00,1.5) {$=$};

\begin{scope}[shift={(7.25,0)}]
\draw[thick,dashed] (-1,3) -- (3.5,3);
\draw[thick,dashed] (-1,0) -- (3.5,0);
\draw[thick,<-] (2.5,1.5) arc (0:360:1);

\draw[thick,fill] (1.10,2.35) arc (0:360:1.5mm);
\node at (0.35,2.50) {$a_n$};
\draw[thick,fill] (0.65,1.75) arc (0:360:1.5mm);
\node at (-0.35,1.95) {$a_{n-1}$};
\node at (0.15,1.1) {\rotatebox[origin=c]{-20}{$\ddots$}};
\draw[thick,fill] (1.10,0.65) arc (0:360:1.5mm);
\node at (0.50,0.40) {$a_1$};

\node at (4.75,1.5) {$=$};
\node at (7.25,1.5) {$\tr_{\Bool Q}(\delta_{\omega})$};
\end{scope}

\end{scope}

\end{tikzpicture}
    \caption{Evaluation of an $\omega$-decorated circle is $1$ if and only if there is a closed path $\omega$ in the graph of $(Q)$.
    }
    \label{figure-4a}
\end{figure}

%% file: mfig_022.tex
\begin{figure}
    \centering
\begin{tikzpicture}[scale=0.6,decoration={
    markings,
    mark=at position 0.6 with {\arrow{>}}}]

\begin{scope}[shift={(0,0)}]

\draw[thick] (0,0) -- (25,0) -- (25,6) -- (0,6) -- (0,0);
\draw[thick] (0,3) -- (25,3);

\draw[thick] (4,0) -- (4,6);
\draw[thick] (12,0) -- (12,6);
\draw[thick] (16,0) -- (16,6);
\draw[thick] (20,0) -- (20,6);

\node at (2,5.2) {$\mathcal{F}(+)$};
\node at (2,4.2) {$+$};
\draw[thick, dashed] (1,3.8) -- (3,3.8); 

\node at (2,2.5) {free};
\node at (2,1.5) {$\Bool$-module};
\node at (2,0.5) {$\Bool$Q};

\node at (7.50,2.5) {coevaluation,};
\node at (8,1.5) {evaluation maps};

\node at (5,0.5) {$\Bool$};
\draw[thick,->] (5.65,0.65) -- (8.15,0.65);
\draw[thick,<-] (5.65,0.35) -- (8.15,0.35);
\node at (10,0.5) {$\Bool Q\otimes \Bool Q^*$};

\node at (14,2) {transition};
\node at (14,1) {map $\delta_a$};

\node at (18,2) {set of initial};
\node at (18,1) {states $Q_{\init}$};

\node at (22.5,2) {set of accepting};
\node at (22.5,1) {states $Q_{\t}$};

\end{scope}

\begin{scope}[shift={(4,3)}]

\draw[thick, dashed] (0.5,2.30) -- (3.5,2.30); 
\draw[thick, dashed] (0.5,0.70) -- (3.5,0.70); 
\draw[thick,postaction={decorate}] (2.75,2.3) .. controls (2.55,0.65) and (1.45,0.65) .. (1.25,2.3);
\node at (1.25,2.6) {$+$};
\node at (2.75,2.6) {$-$};
\node at (4,1.0) {$,$};
\end{scope}

\begin{scope}[shift={(8,3)}]

\draw[thick, dashed] (0.5,2.30) -- (3.5,2.30); 
\draw[thick, dashed] (0.5,0.70) -- (3.5,0.70); 
\draw[thick,postaction={decorate}] (1.25,0.7) .. controls (1.45,2.35) and (2.55,2.35) .. (2.75,0.7);
\node at (1.25,0.4) {$+$};
\node at (2.75,0.4) {$-$};
\end{scope}

\begin{scope}[shift={(12,3)}]
\draw[thick, dashed] (1,2.30) -- (3,2.30); 
\draw[thick, dashed] (1,0.70) -- (3,0.70); 

\node at (2,2.6) {$+$};
\node at (2,0.4) {$+$};

\draw[thick,->] (2,0.7) -- (2,2.3);
\draw[thick,fill] (2.15,1.5) arc (0:360:1.5mm);
\node at (2.5,1.5) {$a$};
\end{scope}

\begin{scope}[shift={(16,3)}]
\draw[thick, dashed] (1,2.30) -- (3,2.30); 
\draw[thick, dashed] (1,0.70) -- (3,0.70); 

\node at (2,2.6) {$+$};

\draw[thick,->] (2,1.20) -- (2,2.3);

\end{scope}

\begin{scope}[shift={(20.50,3)}]
\draw[thick, dashed] (1,2.30) -- (3,2.30); 
\draw[thick, dashed] (1,0.70) -- (3,0.70); 

\node at (2,0.4) {$+$};

\draw[thick,->] (2,0.70) -- (2,1.8);

\end{scope}

\end{tikzpicture}
    \caption{Summary table for the automata $\longleftrightarrow$ Boolean TQFT correspondence.}
    \label{mfig_022}
\end{figure}

%% file: mfig_001.tex
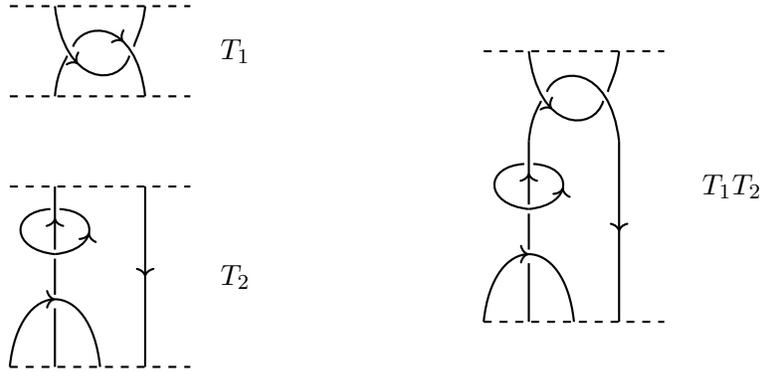
\begin{figure}
    \centering
\begin{tikzpicture}[scale=0.6,decoration={
    markings,
    mark=at position 0.5 with {\arrow{>}}}]

\begin{scope}[shift={(0,4)}]

\draw[thick,dashed] (0,4) -- (4,4);
\draw[thick,dashed] (0,2) -- (4,2);

\node at (5,3) {$T_1$};

\draw[thick] (1,2) .. controls (1.05,2.6) and (1.25,2.8) .. (1.27,2.9);
\draw[thick,postaction={decorate}] (1.4,3.1) .. controls (1.6,3.65) and (2.8,3.75) .. (3,2);

\draw[thick,postaction={decorate}] (1,4) .. controls (1.2,2.2) and (2.45,2.2) .. (2.65,2.9);
\draw[thick] (2.75,3.1) .. controls (2.80,3.3) and (2.95,3.5) .. (3,4);
\end{scope}

\begin{scope}[shift={(0,0)}]

\draw[thick,dashed] (0,4) -- (4,4);
\draw[thick,dashed] (0,0) -- (4,0);

\node at (5,2) {$T_2$};

\draw[thick,postaction={decorate}] (3,4) -- (3,0);

\draw[thick,postaction={decorate}] (0,0) .. controls (0.2,2) and (1.8,2) .. (2,0);
 
\draw[thick] (0.9,3.5) .. controls (0,3.45) and (0,2.6) .. (1,2.5);

\draw[thick,postaction={decorate}] (1,2.5) .. controls (2,2.6) and (2,3.45) .. (1.1,3.5);

\draw[thick,postaction={decorate}] (1,2.6) -- (1,4);
\draw[thick] (1,1.6) -- (1,2.4);
\draw[thick] (1,0) -- (1,1.3);

\end{scope}

\begin{scope}[shift={(10.5,1)}]

\draw[thick,dashed] (0,6) -- (4,6);
\draw[thick,dashed] (0,0) -- (4,0);

\draw[thick,postaction={decorate}] (3,4) -- (3,0);

\draw[thick] (1,4) .. controls (1.05,4.6) and (1.25,4.8) .. (1.27,4.9);
\draw[thick] (1.4,5.1) .. controls (1.6,5.65) and (2.8,5.75) .. (3,4);

\draw[thick,postaction={decorate}] (1,6) .. controls (1.2,4.2) and (2.45,4.2) .. (2.65,4.9);
\draw[thick] (2.75,5.1) .. controls (2.80,5.3) and (2.95,5.5) .. (3,6);

\draw[thick,postaction={decorate}] (0,0) .. controls (0.2,2) and (1.8,2) .. (2,0);

\draw[thick] (0.9,3.5) .. controls (0,3.45) and (0,2.6) .. (1,2.5);

\draw[thick,postaction={decorate}] (1,2.5) .. controls (2,2.6) and (2,3.45) .. (1.1,3.5);

\draw[thick,postaction={decorate}] (1,2.6) -- (1,4);
\draw[thick] (1,1.6) -- (1,2.4);
\draw[thick] (1,0) -- (1,1.3);

\node at (5.5,3) {$T_1T_2$};
\end{scope}

\end{tikzpicture}
    \caption{Two tangles and their composition.}
    \label{mfig_001}
\end{figure}

%% file: mfig_002.tex
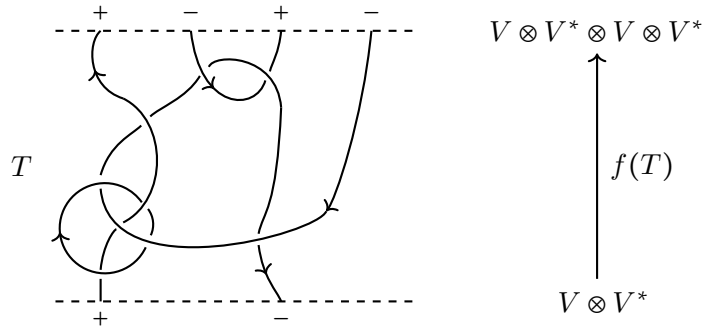
\begin{figure}
    \centering
\begin{tikzpicture}[scale=0.6,decoration={
    markings,
    mark=at position 0.5 with {\arrow{>}}}]
\begin{scope}[shift={(0,0)}]

\draw[thick,dashed] (0,6) -- (8,6);
\draw[thick,dashed] (0,0) -- (8,0);

\node at (12,6) {$V\otimes V^*\otimes V\otimes V^*$};
\draw[thick,->] (12,0.5) -- (12,5.5);
\node at (13,3) {$f(T)$};
\node at (12.15,0) {$V\otimes V^*$};

\node at (1,6.4) {$+$};
\node at (3,6.4) {$-$};
\node at (5,6.4) {$+$};   
\node at (7,6.4) {$-$};

\node at (1,-0.4) {$+$};
\node at (5,-0.4) {$-$};

\draw[thick,postaction={decorate}] (3,6) .. controls (3.2,4.2) and (4.45,4.2) .. (4.65,4.9);
\draw[thick] (4.75,5.1) .. controls (4.80,5.3) and (4.95,5.5) .. (5,6);

\draw[thick] (1,2.5) .. controls (1.2,0.25) and (5.8,1.50) .. (6,2);

\draw[thick,<-] (6,2) .. controls (6.2,2.25) and (6.8,4) .. (7,6);

\draw[thick] (3.4,5.2) .. controls (3.5,5.5) and (4.9,5.5) .. (5,4.3);

\draw[thick] (4.5,1.5) .. controls (4.6,2) and (4.9,2) .. (5,4.3);

\draw[thick,postaction={decorate}] (4.5,1.25) .. controls (4.6,0.50) and (4.9,0.25) .. (5,0);

\draw[thick,postaction={decorate}] (1.5,4.5) .. controls (0.7,4.75) and (0.7,5.75) .. (1,6);

\draw[thick] (1.5,4.5) .. controls (2.5,4) and (2.5,2.25) .. (1.5,1.75);

\draw[thick] (2.1,4.1) .. controls (2.3,4.3) and (2.9,4.5) .. (3.2,5);

\draw[thick,postaction={decorate}] (2,1.2) arc (335:35:1);

\draw[thick] (2.05,2.05) arc (37:-25:0.57);

\draw[thick] (1,2.8) .. controls (1.2,3.5) and (1.7,3.7) .. (1.9,3.9);

\draw[thick] (1,0.7) .. controls (1.05,1.25) and (1.30,1.6) .. (1.35,1.6);

\draw[thick] (1,0) -- (1,0.5);

\node at (-0.75,3) {$T$};

\end{scope}

\end{tikzpicture}
    \caption{This tangle $T$ is a morphism from $(+-)$ to $(+-+-)$. Tensor products of representations for the interwiner $f(T)$ are shown on the right.}
    \label{mfig_002}
\end{figure}

%% file: mfig_003.tex
\begin{figure}
    \centering
\begin{tikzpicture}[scale=0.6]

\draw[thick,dashed] (0,4) -- (4,4);
\draw[thick,dashed] (0,0) -- (4,0);

\draw[thick,->] (2.8,1.6) arc (-60:-200:1);
\draw[thick] (3.08,1.88) arc (-40:170:1);

\draw[thick] (2.5,0.5) arc (-60:120:1);
\draw[thick,->] (1.2,2.0) arc (140:305:1);

\node at (3.75,2) {$L$};

\node at (6,4) {$\mathbb{C}$};
\draw[thick,->] (6,0.5) -- (6,3.5);
\node at (7,2) {$f(L)$};
\node at (6,0) {$\mathbb{C}$};

\end{tikzpicture}

    \caption{ The invariant of a link $L$ is a scalar $f(L)\in \Z[q^{1/D},q^{-1/D}]$.}
    \label{mfig_003}
\end{figure}
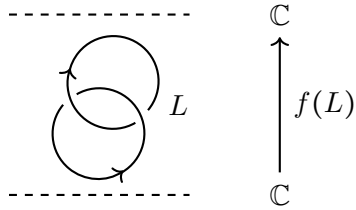

%% file: mfig_004.tex
\begin{tikzpicture}[scale=0.6]

\node at (-6.5,2) {$(III)$};

\node at (-3,2) {$q^NP_N\Big( $};
\node at (0,2) {$\Big)$};

\node at (1.60,2) {$-q^{-N}P_N\Big($};
\node at (5,2) {$\Big)$};
\node at (7.5,2) {$=(q-q^{-1})P_N\Big($}; 

\node at (11.63,2) {$\Big).$};

\draw[thick,<-] (-1.5,2.5) -- (-0.5,1.5);
\draw[thick] (-1.5,1.5) -- (-1.1,1.9);
\draw[thick,->] (-0.9,2.1) -- (-0.5,2.5);

\draw[thick,->] (3.5,1.5) -- (4.5,2.5);
\draw[thick] (4.1,1.9) -- (4.5,1.5);
\draw[thick,<-] (3.5,2.5) -- (3.9,2.1);

\draw[thick,<-] (10,2.5) .. controls (10.35,2.25) and (10.35,1.75) .. (10,1.5);

\draw[thick,<-] (11,2.5) .. controls (10.65,2.25) and (10.65,1.75) .. (11,1.5);

\end{tikzpicture}

%% file: mfig_005.tex
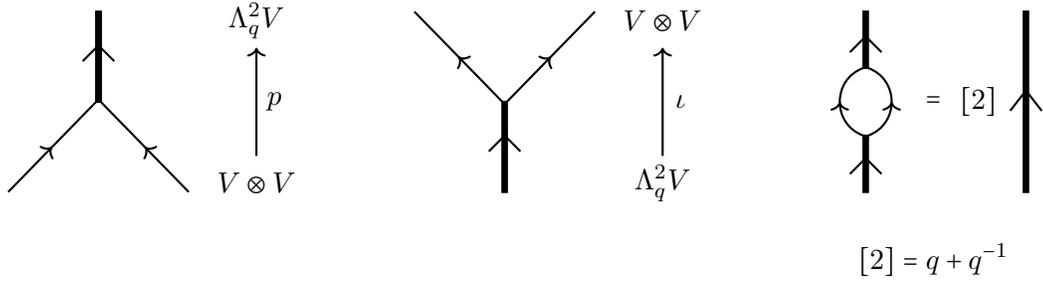
\begin{figure}
    \centering
\begin{tikzpicture}[scale=0.6,decoration={
    markings,
    mark=at position 0.5 with {\arrow{>}}}]
\begin{scope}[shift={(0,0)}]  

\draw[thick] (1.95,2) -- (1.95,4);
\draw[thick, postaction={decorate}] (0,0) -- (1.95,2);

\draw[thick] (2.05,2) -- (2.05,4);
\draw[thick,postaction={decorate}] (4,0) -- (2.05,2);

\draw[thick] (1.65,2.9) -- (2,3.25);
\draw[thick] (2,3.25) -- (2.35,2.9);

\node at (5.5,3.8) {$\Lambda_q^2 V$};
\draw[thick,->] (5.5,0.8) -- (5.5,3.2);
\node at (5.9,2) {$p$};

\node at (5.5,0.2) {$V\otimes V$};

\draw[thick,fill] (1.95,2) rectangle (2.05,4); 
\end{scope}

\begin{scope}[shift={(9,0)}]

\draw[thick,fill] (1.95,0) rectangle (2.05,2); 

\draw[thick] (1.65,0.9) -- (2,1.25);
\draw[thick] (2,1.25) -- (2.35,0.9);

\draw[thick, postaction={decorate}] (1.95,2) -- (0,4);

\draw[thick,postaction={decorate}] (2.05,2) -- (4,4);

\node at (5.5,3.8) {$V\otimes V$};
\draw[thick,->] (5.5,0.8) -- (5.5,3.2);
\node at (5.9,2) {$\iota$};
\node at (5.5,0.2) {$\Lambda_q^2 V$};

\end{scope}


\begin{scope}[shift={(17,0)}]

\draw[thick,fill] (1.95,2.75) rectangle (2.05,4); 
\draw[thick,fill] (1.95,0) rectangle (2.05,1.25); 

\draw[thick] (1.65,3.1) -- (2,3.45);
\draw[thick] (2,3.45) -- (2.35,3.1);

\draw[thick] (1.65,0.4) -- (2,0.75);
\draw[thick] (2,0.75) -- (2.35,0.4);

\draw[thick,postaction={decorate}] (1.95,1.25) .. controls (1.25,1.5) and (1.25,2.5) .. (1.95,2.75);

\draw[thick,postaction={decorate}] (2.05,1.25) .. controls (2.75,1.5) and (2.75,2.5) .. (2.05,2.75);

\node at (3.5,2) {$=$};

\node at (4.5,2) {$[2]$};

\draw[thick, fill] (5.5,0) rectangle (5.6,4);

\draw[thick] (5.2,1.8) -- (5.55,2.25);
\draw[thick] (5.55,2.25) -- (5.9,1.8);

\node at (3.5,-1.5) {$[2]=q+q^{-1}$};
\end{scope}

\end{tikzpicture}
    \caption{Left: diagram of projection $p$ onto $\Lambda^2_q V$. Middle: diagram of inclusion $\iota$ back into $V\otimes V$. Right: these maps are scaled so that the composition $p\circ\iota$ is the identity times $q+q^{-1}$.  }
    \label{mfig_005}
\end{figure}

%% file: mfig_006.tex
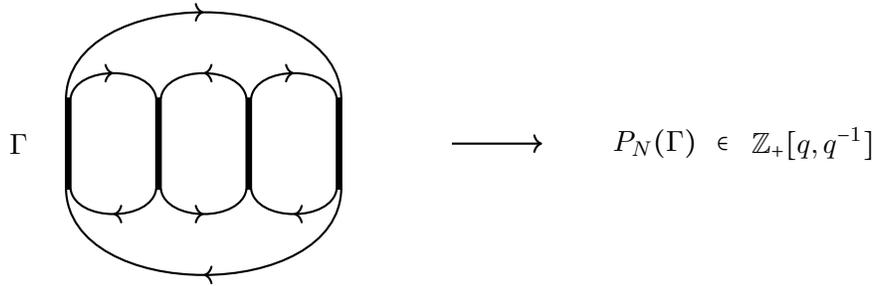
\begin{figure}
    \centering
\begin{tikzpicture}[scale=0.6,decoration={
    markings,
    mark=at position 0.5 with {\arrow{>}}}]
\begin{scope}[shift={(0,0)}]


\node at (-2.10,2) {$\Gamma$};

\draw[thick,fill] (-1.05,1) rectangle (-0.95,3);

\draw[thick,fill] (0.95,1) rectangle (1.05,3);

\draw[thick,fill] (2.95,1) rectangle (3.05,3);

\draw[thick,fill] (4.95,1) rectangle (5.05,3);

\draw[thick,postaction={decorate}] (-0.95,3) .. controls (-0.85,3.75) and (0.85,3.75) .. (0.95,3);

\draw[thick,postaction={decorate}] (0.95,1) .. controls (0.85,0.25) and (-0.85,0.25) .. (-0.95,1);

\draw[thick,postaction={decorate}] (2.95,3) .. controls (2.85,3.75) and (1.15,3.75) .. (1.05,3);

\draw[thick,postaction={decorate}] (1.05,1) .. controls (1.15,0.25) and (2.85,0.25) .. (2.95,1);

\draw[thick,postaction={decorate}] (3.05,3) .. controls (3.15,3.75) and (4.85,3.75) .. (4.95,3);

\draw[thick,postaction={decorate}] (4.95,1) .. controls (4.85,0.25) and (3.15,0.25) .. (3.05,1);

\draw[thick,postaction={decorate}] (-1.05,3) .. controls (-1,5.55) and (5,5.55) .. (5.05,3);

\draw[thick,postaction={decorate}] (5.05,1) .. controls (5,-1.55) and (-1,-1.55) .. (-1.05,1);

\end{scope}

\draw[thick,->] (7.5,2) -- (9.5,2);

\node at (12,2) {$P_N(\Gamma)$};

\node at (13.5,2) {$\in$};

\node at (15.5,2) {$\mathbb{Z}_+[q,q^{-1}]$};

\end{tikzpicture}
   
    \caption{A MOY graph $\Gamma$.}
    \label{mfig_006}
\end{figure}

%% file: mfig_008.tex
\begin{tikzpicture}[scale=0.6]
\begin{scope}[shift={(0,0)}]

\node at (-6.5,2) {$(I)$};

\node at (-2.5,2) {$P_N\Bigg( $};
\node at (0,2) {$\Bigg)$};

\draw[thick,->] (-1.5,1.25) -- (-0.5,2.75);
\draw[thick] (-0.9,1.9) -- (-0.5,1.25);
\draw[thick,<-] (-1.5,2.75) -- (-1.1,2.1);

\node at (2,2) {$=q^{N-1}P_N\Bigg($}; 
\node at (5.6,2) {$\Bigg)$};

\draw[thick,<-] (4,2.75) .. controls (4.40,2.50) and (4.40,1.50) .. (4,1.25);

\draw[thick,<-] (5.1,2.75) .. controls (4.70,2.50) and (4.70,1.50) .. (5.1,1.25);

\begin{scope}[shift={(0.3,0)}]
\node at (7,2) {$-\: q^N P_N\Bigg($};
\node at (10.13,2) {$\Bigg),$};

\draw[thick,->] (8.5,1.25) -- (9.05,1.50);
\draw[thick] (9.05,1.50) -- (9.05,2.50);
\draw[thick,<-]  (8.5,2.75) -- (9.05,2.50);

\draw[thick,->] (9.7,1.25) -- (9.15,1.50);
\draw[thick] (9.15,1.50) -- (9.15,2.50);
\draw[thick,->]  (9.15,2.50) -- (9.7,2.75);

\draw[thick,fill] (9.05,1.5) -- (9.05,2.5) -- (9.15,2.5) -- (9.15,1.5)-- (9.05,1.5);
\end{scope}

\end{scope}
\begin{scope}[shift={(0,-4)}]

\node at (-6.5,2) {$(II)$};

\node at (-2.5,2) {$P_N\Bigg( $};
\node at (0,2) {$\Bigg)$};

\draw[thick,<-] (-1.5,2.75) -- (-0.5,1.25);
\draw[thick] (-1.5,1.25) -- (-1.1,1.9);
\draw[thick,->] (-0.9,2.1) -- (-0.5,2.75);

\node at (2,2) {$=q^{1-N}P_N\Bigg($}; 
\node at (5.6,2) {$\Bigg)$};

\draw[thick,<-] (4,2.75) .. controls (4.40,2.50) and (4.40,1.50) .. (4,1.25);

\draw[thick,<-] (5.1,2.75) .. controls (4.70,2.50) and (4.70,1.50) .. (5.1,1.25);

\begin{scope}[shift={(0.5,0)}]
\node at (7,2) {$-\: q^{-N} P_N\Bigg($};
\node at (10.3,2) {$\Bigg)$.};

\begin{scope}[shift={(0.2,0)}]
  
\draw[thick,->] (8.5,1.25) -- (9.05,1.50);
\draw[thick] (9.05,1.50) -- (9.05,2.50);
\draw[thick,<-]  (8.5,2.75) -- (9.05,2.50);

\draw[thick,->] (9.7,1.25) -- (9.15,1.50);
\draw[thick] (9.15,1.50) -- (9.15,2.50);
\draw[thick,->]  (9.15,2.50) -- (9.7,2.75);

\draw[thick,fill] (9.05,1.5) -- (9.05,2.5) -- (9.15,2.5) -- (9.15,1.5)-- (9.05,1.5);
\end{scope}
\end{scope}
 
\end{scope}

\end{tikzpicture}

%% file: mfig_007.tex
\begin{figure}
    \centering

\begin{tikzpicture}[scale=0.6]
\begin{scope}[shift={(0,0)}]

\draw[thick,->] (0,0) -- (1,1);
\draw[thick] (1,1) -- (2,2);

\draw[thick,<-] (3,1) -- (4,0);
\draw[thick] (2,2) -- (3,1);

\draw[thick,->] (2,2) -- (2,3);
\draw[thick] (2,3) -- (2,4);

\node at (3.35,3) {$a+b$};
\node at (-0.5,0) {$a$};
\node at (4.5,0) {$b$};
\end{scope}

\begin{scope}[shift={(8,0)}]

\draw[thick] (2,1) -- (2,2);
\draw[thick,->] (2,0) -- (2,1);
\node at (3.5,1) {$a+b$};

\draw[thick] (0,4) -- (1,3);
\draw[thick,<-] (1,3) -- (2,2);

\draw[thick] (3,3) -- (4,4);
\draw[thick,->] (2,2) -- (3,3);

\node at (-0.5,3.75) {$a$};
\node at (4.5,3.75) {$b$};
\end{scope}


\begin{scope}[shift={(1,-4.5)}]

\draw[thick,<-] (0,2) arc (180:-180:1);
\node at (2,3) {$a$};

\node at (3,2) {$=$};

\node at (4.75,2) {$\Bigg[\hspace{0.65cm} \Bigg]$};
\node at (6.5,2) {$=$};

\node at (4.75,2.5) {$N$};
\node at (4.75,1.5) {$a$};

\node at (9,2) {$\dfrac{[N]!}{[a]![N-a]!}$};

\end{scope}

\end{tikzpicture}
    \caption{Top: $(a,b)$-vertices of general MOY graphs. Bottom: Circle of thickess $a$ evaluates to a quantum binomial coefficient.}
    \label{mfig_007}
\end{figure}
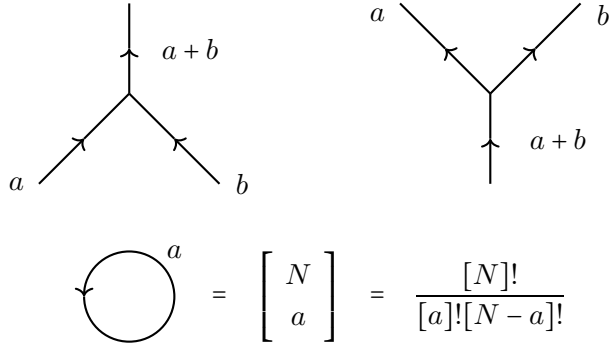

%% file: mfig_023.tex
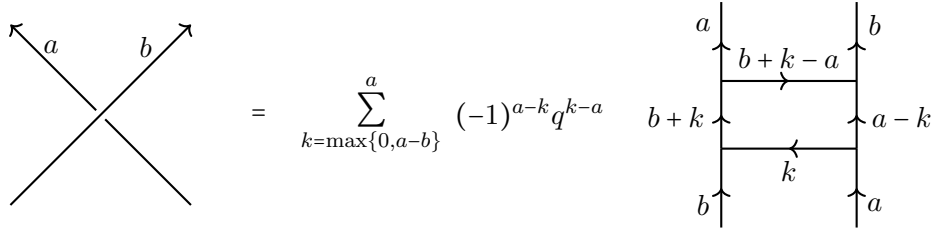
\begin{figure}
    \centering
\begin{tikzpicture}[scale=0.6,decoration={markings,mark=at position 0.5 with {\arrow{>}}}]
\begin{scope}[shift={(0,0)}]

\draw[thick,->] (0,0) -- (4,4);

\draw[thick] (2.1,1.9) -- (4,0);
\draw[thick,<-] (0,4) -- (1.9,2.1); 

\node at (0.9,3.5) {$a$};
\node at (3,3.5) {$b$};

\node at (5.5,2) {$=$};

\node at (8,2) {$\displaystyle{\sum_{k=\max\{0,a-b\}}^a}$};
\node at (11.5,2) {$(-1)^{a-k}q^{k-a}$};

\end{scope}

\begin{scope}[shift={(13.75,-0.25)}]

\draw[thick,postaction={decorate}] (2,3) -- (2,4.75);
\draw[thick,postaction={decorate}] (5,3) -- (5,4.75);

\draw[thick,postaction={decorate}] (2,1.5) -- (2,3);
\draw[thick,postaction={decorate}] (5,1.5) -- (5,3);

\draw[thick,postaction={decorate}] (2,-0.25) -- (2,1.5);
\draw[thick,postaction={decorate}] (5,-0.25) -- (5,1.5);

\draw[thick,postaction={decorate}] (2,3) -- (5,3);
\draw[thick,postaction={decorate}] (5,1.5) -- (2,1.5);

\node at (1.6,4.25) {$a$};
\node at (5.4,4.25) {$b$};

\node at (6,2.2) {$a-k$};
\node at (1,2.2) {$b+k$};

\node at (3.5,3.5) {$b+k-a$};
\node at (3.5,1) {$k$};

\node at (1.6,0.25) {$b$};
\node at (5.4,0.25) {$a$};

\end{scope}

\end{tikzpicture}
    \caption{Decomposing a negative crossing for two strands of arbitrary thickness $a,b\le N$. }
    \label{mfig_023}
\end{figure}

%% file: mfig_024.tex
\begin{figure}
    \centering
\begin{tikzpicture}[scale=0.6,decoration={markings,mark=at position 0.5 with {\arrow{>}}}]
\begin{scope}[shift={(0,0)}]

\draw[thick,->] (4,0) -- (0,4);

\draw[thick,->] (2.1,2.1) -- (4,4);
\draw[thick] (0,0) -- (1.9,1.9); 

\node at (0.9,3.5) {$a$};
\node at (3,3.5) {$b$};

\node at (5.5,2) {$=$};

\node at (8,2) {$\displaystyle{\sum_{k=\max\{0,a-b\}}^a}$};
\node at (11.5,2) {$(-1)^{k-a}q^{a-k}$};

\end{scope}

\begin{scope}[shift={(13.75,-0.25)}]

\draw[thick,postaction={decorate}] (2,3) -- (2,4.75);
\draw[thick,postaction={decorate}] (5,3) -- (5,4.75);

\draw[thick,postaction={decorate}] (2,1.5) -- (2,3);
\draw[thick,postaction={decorate}] (5,1.5) -- (5,3);

\draw[thick,postaction={decorate}] (2,-0.25) -- (2,1.5);
\draw[thick,postaction={decorate}] (5,-0.25) -- (5,1.5);

\draw[thick,postaction={decorate}] (2,3) -- (5,3);
\draw[thick,postaction={decorate}] (5,1.5) -- (2,1.5);

\node at (1.6,4.25) {$a$};
\node at (5.4,4.25) {$b$};

\node at (6,2.2) {$a-k$};
\node at (1,2.2) {$b+k$};

\node at (3.5,3.5) {$b+k-a$};
\node at (3.5,1) {$k$};

\node at (1.6,0.25) {$b$};
\node at (5.4,0.25) {$a$};

\end{scope}

\end{tikzpicture}
    \caption{Decomposition of a positive crossing with thicknesses $a,b$.}
    \label{mfig_024}
\end{figure}
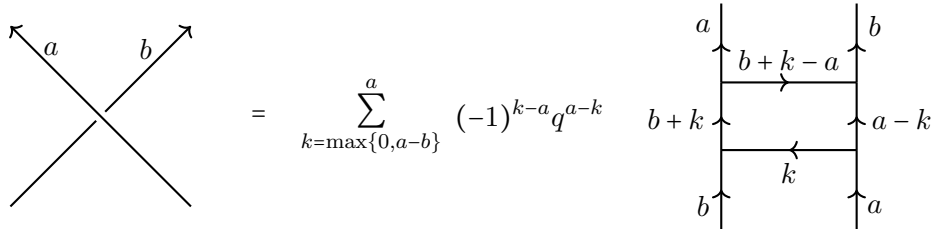

%% file: mfig_009.tex
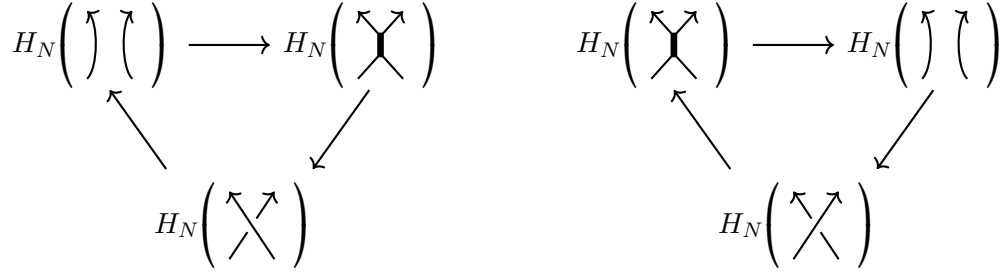
\begin{figure}
    \centering
\begin{tikzpicture}[scale=0.6]
\begin{scope}[shift={(0,0)}]

\node at (-0.2,4) {$H_N\Bigg(\hspace{1cm} \Bigg)$};

\draw[thick,->] (2,4) -- (3.8,4);

\node at (5.8,4) {$H_N\Bigg(\hspace{1cm} \Bigg)$};

\node at (2.95,0) {$H_N\Bigg(\hspace{1cm} \Bigg)$};

\draw[thick,<-] (0.25,3) -- (1.50,1.25);
\draw[thick,<-] (4.75,1.25) -- (6.00,3);

\draw[thick,<-] (-0.25,4.75) .. controls (0,4.50) and (0,3.50) .. (-0.25,3.25);

\draw[thick,<-] (0.75,4.75) .. controls (0.50,4.50) and (0.50,3.50) .. (0.75,3.25);

\begin{scope}[shift={(-0.1,-2)}]
\draw[thick,<-] (3,2.75) -- (4,1.25);
\draw[thick] (3,1.25) -- (3.4,1.85);
\draw[thick,->] (3.6,2.15) -- (4,2.75);
\end{scope}

\begin{scope}[shift={(0.25,0)}]
\draw[thick] (5.5,3.25) -- (5.95,3.75);
\draw[thick] (6.05,3.75) -- (6.5,3.25);

\draw[thick,fill] (5.95,3.75) rectangle (6.05,4.25);

\draw[thick,<-] (5.5,4.75) -- (5.95,4.25);
\draw[thick,->] (6.05,4.25) -- (6.5,4.75);
\end{scope}
\end{scope}

\begin{scope}[shift={(12.5,0)}]

\node at (-0.2,4) {$H_N\Bigg(\hspace{1cm} \Bigg)$};

\draw[thick,->] (2,4) -- (3.8,4);

\node at (5.8,4) {$H_N\Bigg(\hspace{1cm} \Bigg)$};

\node at (2.95,0) {$H_N\Bigg(\hspace{1cm} \Bigg)$};

\draw[thick,<-] (0.25,3) -- (1.50,1.25);
\draw[thick,<-] (4.75,1.25) -- (6.00,3);

\begin{scope}[shift={(6,0)}]
\draw[thick,<-] (-0.25,4.75) .. controls (0,4.50) and (0,3.50) .. (-0.25,3.25);

\draw[thick,<-] (0.75,4.75) .. controls (0.50,4.50) and (0.50,3.50) .. (0.75,3.25);
\end{scope}

\begin{scope}[shift={(-0.1,-2)}]

\draw[thick,->] (3,1.25) -- (4,2.75);

\draw[thick] (4,1.25) -- (3.6,1.85);
\draw[thick,->] (3.4,2.15) -- (3,2.75);
\end{scope}

\begin{scope}[shift={(-5.75,0)}]
\draw[thick] (5.5,3.25) -- (5.95,3.75);
\draw[thick] (6.05,3.75) -- (6.5,3.25);

\draw[thick,fill] (5.95,3.75) rectangle (6.05,4.25);

\draw[thick,<-] (5.5,4.75) -- (5.95,4.25);
\draw[thick,->] (6.05,4.25) -- (6.5,4.75);
\end{scope}
\end{scope}

\end{tikzpicture}
    \caption{Long exact sequences to compute homology of a crossing given those of its planar resolutions.}
    \label{mfig_009}
\end{figure}

%% file: mfig_011.tex
\begin{tikzpicture}[scale=0.6]
\node at (0,2) {$M^0$};
\node at (4,2) {$M^1$};

\draw[thick,->] (0.5,2.5) .. controls (0.75,3.1) and (3.25,3.1) .. (3.5,2.5);
\draw[thick,<-] (0.5,1.5) .. controls (0.75,0.9) and (3.25,0.9) .. (3.5,1.5);

\node at (2,3.4) {$d_0$};
\node at (2,0.6) {$d_1$};
\end{tikzpicture}

%% file: complex_001.tex
\begin{center}
\begin{tikzpicture}[scale=0.6]
\node at (0,2) {$S$};
\node at (3,2) {$S$};
\node at (6,2) {$S$};

\draw[thick,->] (0.5,2) -- (2.5,2);
\draw[thick,->] (3.5,2) -- (5.5,2);
\node at (1.5,2.5) {$\pi_{12}$};
\node at (4.5,2.5) {$x_1-x_2$};

\node at (9,2) {where};
\node at (12,2) {$\pi_{12}$};
\node at (13,2) {$=$};
\node at (15.5,2) {$\dfrac{x_1^{N+1}-x_2^{N+1}}{x_1-x_2}$};
\node at (18,2) {$=$};
\node at (20.25,2) {$\displaystyle{\sum_{i=0}^N} \:\: x_1^i x_2^{N-i}.$};
\end{tikzpicture}
\end{center}

%% file: mfig_012.tex
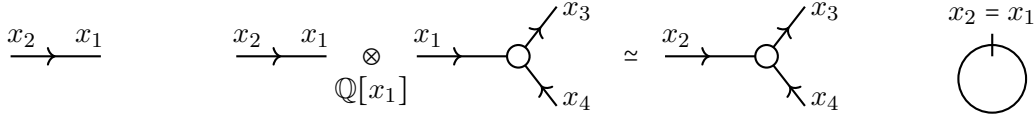
\begin{figure}
    \centering
\begin{tikzpicture}[scale=0.6,decoration={
    markings,
    mark=at position 0.5 with {\arrow{>}}}]
\begin{scope}[shift={(0,0)}]

\draw[thick,postaction={decorate}] (0,2) -- (2,2);

\node at (0.25,2.4) {$x_2$};
\node at (1.75,2.4) {$x_1$};
\end{scope}

\begin{scope}[shift={(5,0)}]

\draw[thick,postaction={decorate}] (0,2) -- (2,2);

\node at (0.25,2.4) {$x_2$};
\node at (1.75,2.4) {$x_1$};
\node at (3,2) {$\otimes$};

\node at (3,1.25) {$\Q[x_1]$};

\draw[thick,postaction={decorate}] (4,2) -- (6,2);

\node at (4.25,2.4) {$x_1$};
\draw[thick] (6.5,2) arc (0:360:0.25);

\draw[thick,postaction={decorate}] (6.40,2.2) -- (7.1,3.1);
\node at (7.55,3) {$x_3$};

\draw[thick,postaction={decorate}] (7.1,0.90) -- (6.40,1.80);
\node at (7.55,1) {$x_4$};

\node at (8.75,2) {$\simeq$};
\end{scope}

\begin{scope}[shift={(10.5,0)}]
\draw[thick,postaction={decorate}] (4,2) -- (6,2);

\node at (4.25,2.4) {$x_2$};
\draw[thick] (6.5,2) arc (0:360:0.25);

\draw[thick,postaction={decorate}] (6.40,2.2) -- (7.1,3.1);
\node at (7.55,3) {$x_3$};

\draw[thick,postaction={decorate}] (7.1,0.90) -- (6.40,1.80);
\node at (7.55,1) {$x_4$};
\end{scope}

\begin{scope}[shift={(20.25,-0.5)}]
\draw[thick] (2.25,2) arc (0:360:0.75);

\node at (1.5,3.4) {$x_2=x_1$};

\draw[thick] (1.5,3) -- (1.5,2.5);

\end{scope}

\end{tikzpicture}
    \caption{Left: an arc with variables as labels at the end point. Middle: tensoring with a factorization $M$, forgetting variable $x_1$, and relabeling $x_2$ back to $x_1$ gives a functor isomorphic to the identity functor in the homotopy category $MF_\omega$. Right: closing up an arc into a circle and equating $x_1=x_2$.}
    \label{mfig_012}
\end{figure}

%% file: mfig_013.tex
\begin{figure}
    \centering
\begin{tikzpicture}[scale=0.6,decoration={
    markings,
    mark=at position 0.5 with {\arrow{>}}}]
\begin{scope}[shift={(0,0)}]

\draw[thick,postaction={decorate}] (0,0) -- (1.935,1.25);
\draw[thick,postaction={decorate}] (4,0) -- (2.065,1.25);

\draw[thick,fill] (1.935,1.25) rectangle (2.065,2.75);

\draw[thick,postaction={decorate}] (1.935,2.75) -- (0,4);
\draw[thick,postaction={decorate}] (2.065,2.75) -- (4,4);

\node at (-0.25,4.4) {$x_1$};
\node at (4.25,4.4) {$x_2$};

\node at (-0.25,-0.4) {$x_3$};
\node at (4.25,-0.4) {$x_4$};

\node at (6.5,2) {$\simeq$};
\end{scope}

\begin{scope}[shift={(8,0)}]

\draw[thick] (1.95,0) -- (1.95,4);
\draw[thick] (2.05,0) -- (2.05,4);

\draw[thick,fill] (1.95,0) rectangle (2.05,4);

\draw[thick] (1.65,1.75) -- (2,2.25);
\draw[thick] (2,2.25) -- (2.35,1.75);

\node at (2, 4.4) {$x_1+x_2,x_1x_2$};
\node at (2,-0.4) {$x_3+x_4,x_3x_4$};

\node at (5.00,2) {$\simeq$};
\end{scope}

\begin{scope}[shift={(14.25,0)}]

\draw[thick,postaction={decorate}] (1,0) -- (1,4);
\draw[thick,postaction={decorate}] (2.75,0) -- (2.75,4);

\node at (1,4.4) {$x_1$};
\node at (2.75,4.4) {$x_2$};

\node at (1,-0.4) {$x_3$};
\node at (2.75,-0.4) {$x_4$};

\end{scope}

\end{tikzpicture}
    \caption{Double edge factorization can be thought of as the identity factorization for the subring of symmetric functions in two variables.}
    \label{mfig_013}
\end{figure}
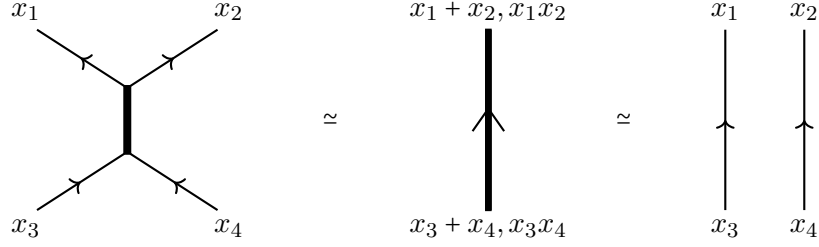

%% file: mfig_014.tex
\begin{figure}
    \centering
\begin{tikzpicture}[scale=0.6,decoration={
    markings,
    mark=at position 0.85 with {\arrow{>}}}]
\begin{scope}[shift={(0,0)}]

\draw[thick,postaction={decorate}] (0,3) .. controls (0.2,3) and (1.5,3) .. (2,2.05);
\draw[thick,postaction={decorate}] (0,1) .. controls (0.2,1) and (1.5,1) .. (2,1.95);

\draw[thick,fill] (2,1.95) rectangle (4,2.05);

\draw[thick,->] (4,2.05) .. controls (4.5,3.5) and (6,3.5) .. (7.5,3.5);

\draw[thick,->] (4,1.95) .. controls (5,1) and (5.5,1) .. (6,1);

\node at (1.25,3.4) {$i_3$};
\draw[thick] (0.9,2.6) -- (1.1,3.1);

\node at (1.25,0.5) {$i_4$};
\draw[thick] (0.9,1.4) -- (1.1,0.9);

\node at (4.6,3.7) {$i_1$};
\draw[thick] (4.7,3.3) -- (5,2.8);

\node at (4.6,0.75) {$i_2$};
\draw[thick] (4.8,1.0) -- (5.0,1.5);

\node at (6.25,4.15) {$i_5$};
\draw[thick] (6.25,3.75) -- (6.25,3.25);

\end{scope}

\begin{scope}[shift={(11,0)}]

\draw[thick] (0,2.5) -- (2,1.05);
\draw[thick] (0,-0.5) -- (2,0.95);

\draw[thick,fill] (2,0.95) rectangle (4,1.05);

\draw[thick,postaction={decorate}] (4,1.05) -- (6,2.45);
\draw[thick] (4,0.95) -- (6,-0.5);

\draw[thick,fill] (6,2.55) rectangle (8,2.45);

\draw[thick] (4,4) -- (6,2.55);

\draw[thick] (8,2.55) -- (10,4);
\draw[thick] (8,2.45) -- (10,0.5);

\draw[thick] (4.80,2.00) -- (5.30,1.45);
\node at (5.6,1.1){$i$}; 
\end{scope}

\end{tikzpicture}
    \caption{Left: markings around a thick edge and an additional mark on an adjacent thin edge. Right: potentials $\omega$ add up to $0$ for a closed diagram, due to cancellation of terms. For a mark labelled $i$, the contribution of potentials is $x_i^{N+1}-x_i^{N+1}=0$, with the signs opposite due to opposite orientations.}
    \label{mfig_014}
\end{figure}
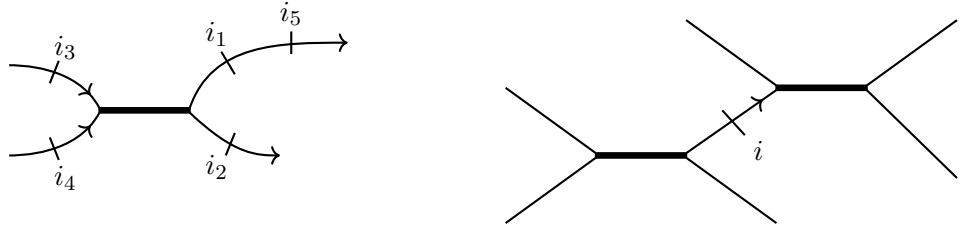

%% file: mfig_015.tex
\begin{figure}
    \centering

\begin{tikzpicture}[scale=0.6,decoration={
    markings,
    mark=at position 0.7 with {\arrow{>}}}]
\begin{scope}[shift={(0,0)}]

\draw[thick,->] (0,0) .. controls (0.5,0.5) and (0.5,2.5) .. (0,3);

\draw[thick,->] (2,0) .. controls (1.5,0.5) and (1.5,2.5) .. (2,3);

\draw[thick,->] (3,1.75) -- (4,1.75);
\draw[thick,<-] (3,1.25) -- (4,1.25);

\node at (1,-0.75) {$\Gamma_0$};

\end{scope}

\begin{scope}[shift={(4.75,0)}]

\draw[thick,postaction={decorate}] (0.9,2) -- (0,3);
\draw[thick,postaction={decorate}] (1.1,2) -- (2,3);

\draw[line width = 0.05in] (1,1) -- (1,2);

\draw[thick,postaction={decorate}] (0,0) -- (0.9,1);
\draw[thick,postaction={decorate}] (2,0) -- (1.1,1);

\node at (1,-0.75) {$\Gamma_1$};

\end{scope}

\begin{scope}[shift={(8.75,0)}]

\node at (0,1.5) {$M_{\Gamma_0}$};

\draw[thick,->] (1,1.75) -- (2,1.75);
\draw[thick,<-] (1,1.25) -- (2,1.25);
\node at (1.5,2.2) {$\chi_0$};
\node at (1.5,0.75) {$\chi_1$};

\node at (3,1.5) {$M_{\Gamma_1}$};
\end{scope}

\begin{scope}[shift={(13.75,-0.25)}]

\node at (0,3) {$0$}; 
\draw[thick,->] (0.5,3) -- (1,3);

\draw[thick,<-] (1.25,3.75) .. controls (1.5,3.5) and (1.5,2.5) .. (1.25,2.25);

\draw[thick,->] (2,2.25) .. controls (1.75,2.5) and (1.75,3.5) .. (2,3.75);

\draw[thick,->] (2.5,3) -- (3,3);

\begin{scope}[shift={(0.1,0)}]
\draw[line width=0.03in] (3.5,2.75) -- (3.5,3.25);
\draw[thick,postaction={decorate}] (3.05,2.25) -- (3.45,2.75);
\draw[thick,postaction={decorate}] (3.95,2.25) -- (3.55,2.75);

\draw[thick,postaction={decorate}] (3.45,3.25) -- (3.05,3.75);
\draw[thick,postaction={decorate}] (3.55,3.25) -- (3.95,3.75);
\end{scope}

\draw[thick,->] (4.25,3) -- (4.75,3);

\node at (5.25,3) {$0$};
\end{scope}

\begin{scope}[shift={(13.75,-2.75)}]

\node at (0,3) {$0$}; 
\draw[thick,->] (0.5,3) -- (1,3);

\begin{scope}[shift={(2,0)}]
\draw[thick,<-] (1.25,3.75) .. controls (1.5,3.5) and (1.5,2.5) .. (1.25,2.25);

\draw[thick,->] (2,2.25) .. controls (1.75,2.5) and (1.75,3.5) .. (2,3.75);
\end{scope}

\draw[thick,->] (2.5,3) -- (3,3);

\begin{scope}[shift={(-2,0)}]
\draw[line width=0.03in] (3.5,2.75) -- (3.5,3.25);
\draw[thick,postaction={decorate}] (3.05,2.25) -- (3.45,2.75);
\draw[thick,postaction={decorate}] (3.95,2.25) -- (3.55,2.75);

\draw[thick,postaction={decorate}] (3.45,3.25) -- (3.05,3.75);
\draw[thick,postaction={decorate}] (3.55,3.25) -- (3.95,3.75);
\end{scope}

\draw[thick,->] (4.25,3) -- (4.75,3);

\node at (5.25,3) {$0$};  
\end{scope}

\begin{scope}[shift={(20.5,0.375)}]

\draw[thick,->] (1.25,1.75) -- (0,3);
 
\draw[thick] (0,1.75) -- (0.55,2.30);
\draw[thick,->] (0.7,2.46) -- (1.25,3);
\end{scope}

\begin{scope}[shift={(20.5,-2.125)}]
\draw[thick,->] (0,1.75) -- (1.25,3);

\draw[thick,<-] (0,3) -- (0.55,2.45);
\draw[thick] (0.7,2.3) -- (1.25,1.75);
\end{scope}

\end{tikzpicture}

    \caption{Left: MOY graphs $\Gamma_0,\Gamma_1$ and corresponding factorizations. Right: forming complexes for positive and negative crossings out of factorizations $M_{\Gamma_0}$ and $M_{\Gamma_1}$.}
    \label{mfig_015}
\end{figure}
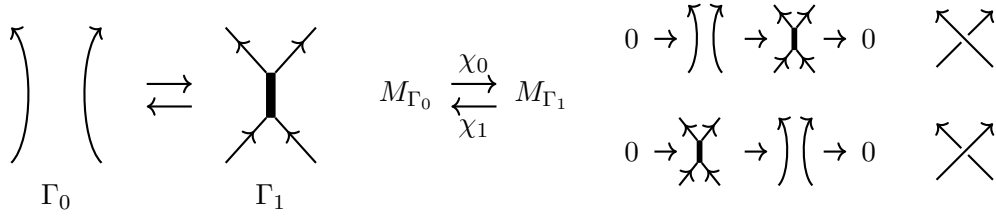

%% file: mfig_020.tex
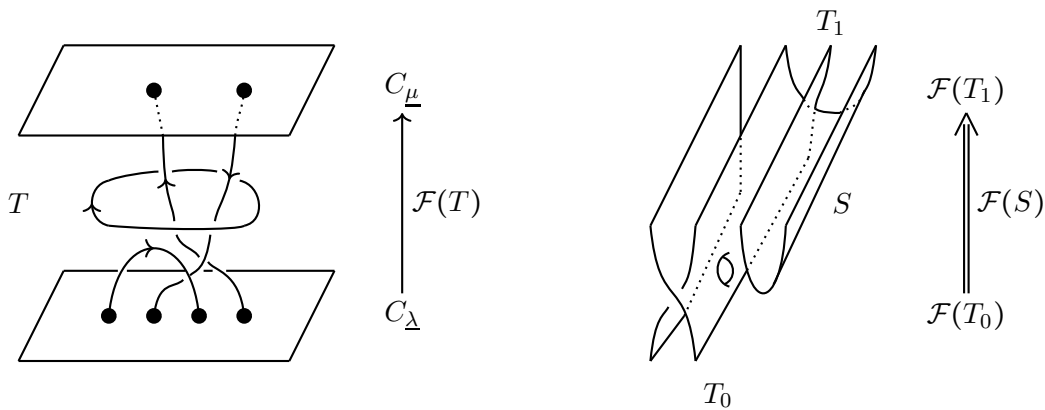
\begin{figure}
    \centering
\begin{tikzpicture}[scale=0.6,decoration={
    markings,
    mark=at position 0.5 with {\arrow{>}}}]

\begin{scope}[shift={(0,0)}]

\draw[thick] (-1,0) -- (0,2);
\draw[thick] (5,0) -- (6,2);

\draw[thick] (-1,0) -- (5,0);

\draw[thick,postaction={decorate}] (1,1) .. controls (1.2,3) and (2.8,3) .. (3,1);

\draw[thick,fill] (2.15,6) arc (0:360:1.5mm);
\draw[thick,fill] (4.15,6) arc (0:360:1.5mm);

\draw[thick,fill] (1.15,1) arc (0:360:1.5mm);
\draw[thick,fill] (2.15,1) arc (0:360:1.5mm);
\draw[thick,fill] (3.15,1) arc (0:360:1.5mm);
\draw[thick,fill] (4.15,1) arc (0:360:1.5mm);

\draw[thick,postaction={decorate}] (1,3) .. controls (0.5,3.1) and (0.5,3.9) .. (1,4);

\node at (-1,3.5) {$T$};

\draw[thick] (4,3) .. controls (4.5,3.1) and (4.4,3.9) .. (3.9,4);

\draw[thick] (1,3) .. controls (2.1,2.9) and (3.9,2.9) .. (4,3);

\draw[thick] (2,1) .. controls (2.1,1.8) and (2.45,1.6) .. (2.65,1.8);

\draw[thick] (2.8,1.9) .. controls (2.9,2) and (3.15,2) .. (3.25,2.8);

\draw[thick,dotted] (3.8,5) .. controls (3.85,5.5) and (3.95,5.5) .. (4,6);
\draw[thick,postaction={decorate}] (3.8,5) .. controls (3.7,3.6) and (3.45,3.6) .. (3.35,3.1); 

\draw[thick,dotted] (2,6) .. controls (2.05,5.75) and (2.15,5.25) .. (2.2,5);
\draw[thick,postaction={decorate}] (2.45,3.1) .. controls (2.43,3.5) and (2.22,3.5) .. (2.2,5);

\draw[thick] (2.5,2.8) .. controls (2.6,2.5) and (2.9,2.5) .. (3,2.3);

\draw[thick] (4,1) .. controls (3.9,2) and (3.3,2) .. (3.2,2.2);

\draw[thick] (0,7) -- (6,7);
\draw[thick] (-1,5) -- (5,5);
\draw[thick] (-1,5) -- (0,7);
\draw[thick] (5,5) -- (6,7);

\node at (7.5,6) {$C_{\underline{\mu}}$};
\draw[thick,->] (7.5,1.5) -- (7.5,5.5);
\node at (8.5,3.5) {$\mcF(T)$};
\node at (7.5,1) {$C_{\underline{\lambda}}$};

\draw[thick] (0,2) -- (1.2,2);
\draw[thick] (1.4,2) -- (2.7,2);
\draw[thick] (3.1,2) -- (3.3,2);
\draw[thick] (3.6,2) -- (6,2);

\draw[thick] (1,4) .. controls (1.2,4.1) and (1.9,4.2) .. (2.1,4.2);
\draw[thick] (2.4,4.2) .. controls (2.6,4.25) and (3.2,4.25) .. (3.4,4.2);

\end{scope}

\begin{scope}[shift={(13,0)}]

\draw[thick] (0,3) .. controls (0.2,1) and (0.8,2) .. (1,0);
\draw[thick] (0,0) .. controls (0.1,1) and (0.3,1) .. (0.4,1.2);
\draw[thick] (0.7,1.6) .. controls (0.8,1.8) and (0.9,1.8) .. (1,3);

\draw[thick] (2,3) .. controls (2.1,1) and (2.9,1) .. (3,3);

\draw[thick] (3,7) .. controls (3.1,6) and (3.3,6) .. (3.4,5.8);

\draw[thick,dotted] (3.4,5.8) .. controls (3.5,5.6) and (3.55,5.6) .. (3.65,5.6);

\draw[thick] (3.65,5.6) .. controls (3.75,5.5) and (4.1,5.5) .. (4.2,5.5);

\draw[thick,dotted] (4.2,5.5) .. controls (4.3,5.6) and (4.6,5.6) .. (4.7,5.9);
\draw[thick] (4.7,5.9) .. controls (4.8,6.1) and (4.9,6.1) .. (5,7);

\draw[thick] (2.73,1.7) -- (4.75,6.0);

\draw[thick,dotted] (2,3.75) -- (2,5);
\draw[thick] (2,5) -- (2,7);

\draw[thick] (0,0) -- (0.75,1);
\draw[thick,dotted] (0.75,1) -- (2,3.75);

\draw[thick] (3.65,5.6) .. controls (3.7,6) and (3.9,6) .. (4,7);
\draw[thick,dotted] (3.5,4.5) .. controls (3.6,5) and (3.55,5) .. (3.65,5.6); 

\draw[thick] (1,0) -- (2.14,2.01);
\draw[thick,dotted] (2.14,2.01) -- (3.5,4.5);

\begin{scope}[shift={(-0.1,-0.1)}]
\draw[thick] (1.85,2.5) .. controls (1.5,2.4) and (1.5,1.85) .. (1.85,1.75);
\draw[thick] (1.71,2.4) .. controls (2,2.3) and (2,1.9) .. (1.73,1.8);
\end{scope}

\draw[thick] (0,3) -- (2,7);
\draw[thick] (1,3) -- (3,7);
\draw[thick] (2,3) -- (4,7);
\draw[thick] (3,3) -- (5,7);

\begin{scope}[shift={(-0.5,0)}]
\node at (7.5,6) {$\mcF(T_1)$};
\draw[thick] (7.45,1.5) -- (7.45,5.25);
\draw[thick] (7.55,1.5) -- (7.55,5.25);

\draw[thick] (7.25,5.0) -- (7.5,5.5);
\draw[thick] (7.75,5.0) -- (7.5,5.5);

\node at (8.5,3.5) {$\mcF(S)$};
\node at (7.5,1) {$\mcF(T_0)$};
\end{scope}

\node at (4,7.5) {$T_1$};
\node at (1.5,-0.75) {$T_0$};
\node at (4.25,3.5) {$S$};
\end{scope}

\end{tikzpicture}
    \caption{Left: To a tangle $T$ there is assigned an exact functor $\mcF(T)$. Middle: A tangle cobordism $S$ between tangles $T_0$ and $T_1$. Right: A natural transformation $\mcF(S)$ is associated to $S$.}
    \label{mfig_020}
\end{figure}

%% file: mfig_010.tex
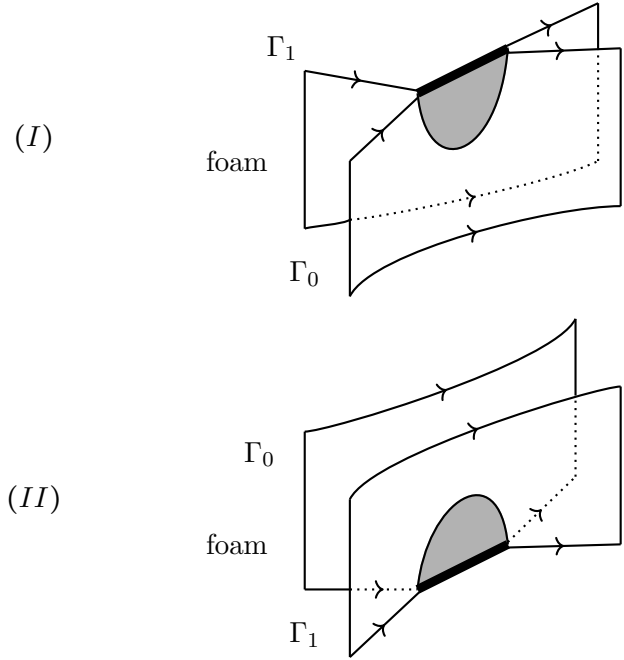
\begin{figure}
    \centering
\begin{tikzpicture}[scale=0.6,decoration={
    markings,
    mark=at position 0.5 with {\arrow{>}}}]
\begin{scope}[shift={(0,0)}]

\node at (-4.5,2) {$(I)$};

\draw[thick,fill,black!30] (4,3) .. controls (4.2,1) and (5.8,1.5) .. (6,4) -- (4,3);

\draw[line width=0.12cm] (4,3) -- (6,4);
\draw[thick] (4,3) .. controls (4.2,1) and (5.8,1.5) .. (6,4);

\draw[thick,postaction={decorate}] (2.5,1.5) -- (4,2.9);

\draw[thick,postaction={decorate}] (1.5,3.5) -- (4,3.06);

\draw[thick] (2.5,-1.5) -- (2.5,1.5);
\draw[thick] (1.5,0) -- (1.5,3.5);

\draw[thick,postaction={decorate}] (5.95,4.04) -- (8,5);
\draw[thick,postaction={decorate}] (6.0,3.9) -- (8.5,4);

\draw[thick] (8,4) -- (8,5);
\draw[thick,dotted] (8,1.5) -- (8,4);
\draw[thick] (8.5,0.5) -- (8.5,4);

\draw[thick,postaction={decorate}] (2.5,-1.5) .. controls (3,-0.5) and (7,0.5) .. (8.5,0.5); 

\node at (0,1.5) {foam};
\node at (1,4) {$\Gamma_1$};
\node at (1.5,-1) {$\Gamma_0$};

\draw[thick,dotted,postaction={decorate}] (2.5,0.2) .. controls (3,0.2) and (7,1) .. (8,1.5);

\draw[thick] (1.5,0) .. controls (1.6,0.05) and (2.4,0.1) .. (2.5,0.2);
\end{scope}

\begin{scope}[shift={(0,-8)}]

\node at (-4.5,2) {$(II)$};

\node at (0,1) {foam};
\node at (0.5,3) {$\Gamma_0$};
\node at (1.5,-1) {$\Gamma_1$};

\draw[thick,fill,black!30] (4,0) .. controls (4.2,2) and (5.8,3) .. (6,1) -- (4,0);

\draw[line width=0.12cm] (4,0) -- (6,1);
\draw[thick] (4,0) .. controls (4.2,2) and (5.8,3) .. (6,1);

\draw[thick] (1.5,0) -- (2.5,0);
\draw[thick,dotted,postaction={decorate}] (2.5,0) -- (4,0);

\draw[thick,postaction={decorate}] (2.5,-1.5) -- (4.05,-0.05);
 
\draw[thick] (1.5,0) -- (1.5,3.5);
\draw[thick] (2.5,-1.5) -- (2.5,2);

\draw[thick,postaction={decorate},dotted] (6,1.05) -- (7.5,2.5);
\draw[thick,postaction={decorate}] (6,0.93) -- (8.5,1);

\draw[thick] (8.5,1) -- (8.5,4.5);

\draw[thick] (7.5,4.3) -- (7.5,6);
\draw[thick,dotted] (7.5,2.5) -- (7.5,4.3);

\draw[thick,postaction={decorate}] (1.5,3.5) .. controls (2,3.5) and (7,5) .. (7.5,6);
\draw[thick,postaction={decorate}] (2.5,2) .. controls (3,3) and (8,4.5) .. (8.5,4.5);

\end{scope}

\end{tikzpicture}
    \caption{Foam cobordisms between two resolutions $\Gamma_0,\Gamma_1$ of a crossing, also see Figure~\ref{mfig_015}.}
    \label{mfig_010}
\end{figure}

%% file: mfig_016.tex
\begin{figure}
    \centering
\begin{tikzpicture}[scale=0.6,decoration={
    markings,
    mark=at position 0.5 with {\arrow{>}}}]

\begin{scope}[shift={(0,0)}]

\draw[thick,fill,black!30] (4,3) .. controls (4.2,1) and (5.8,1.5) .. (6,4) -- (4,3);

\draw[line width=0.12cm] (4,3) -- (6,4);
\draw[thick] (4,3) .. controls (4.2,1) and (5.8,1.5) .. (6,4);

\draw[thick,postaction={decorate}] (2.5,1.5) -- (4,2.9);

\draw[thick,postaction={decorate}] (1.5,3.5) -- (4,3.06);

\draw[thick] (2.5,-1.5) -- (2.5,1.5);
\draw[thick] (1.5,0) -- (1.5,3.5);

\draw[thick,postaction={decorate}] (5.95,4.04) -- (8,5);
\draw[thick,postaction={decorate}] (6.0,3.9) -- (8.5,4);

\draw[thick] (8,4) -- (8,5);
\draw[thick,dotted] (8,1.5) -- (8,4);
\draw[thick] (8.5,0.5) -- (8.5,4);

\draw[thick,postaction={decorate}] (2.5,-1.5) .. controls (3,-0.5) and (7,0.5) .. (8.5,0.5); 

\node at (0.75,1.5) {$F$};
\node at (1,4) {$\Gamma_1$};
\node at (1.5,-1) {$\Gamma_0$};

\draw[thick,dotted,postaction={decorate}] (2.5,0.2) .. controls (3,0.2) and (7,1) .. (8,1.5);

\draw[thick] (1.5,0) .. controls (1.6,0.05) and (2.4,0.1) .. (2.5,0.2);
\end{scope}

\begin{scope}[shift={(10,0)}]

\node at (0,4.5) {$H(\Gamma_1)$};
\draw[thick,->] (0,1.5) -- (0,4);
\node at (0.75,2.75) {$[F]$};
\node at (0,1) {$H(\Gamma_0)$};
\end{scope}

\begin{scope}[shift={(16,0)}]

\draw[thick,->] (0.25,1.25) -- (0.25,2.75);
\draw[thick,->] (1.75,4) -- (4,3.25);
\draw[thick,<-] (1.75,0) -- (4,0.75);

\node at (0,4) {$H\Bigg( \hspace{1cm}\Bigg)$};

\node at (0,0) {$H\Bigg( \hspace{1cm}\Bigg)$};

\node at (4,2) {$H\Bigg( \hspace{1cm}\Bigg)$};

\begin{scope}[shift={(-1.5,-3)}]
\draw[thick,<-] (1.25,3.75) .. controls (1.5,3.5) and (1.5,2.5) .. (1.25,2.25);
\draw[thick,->] (2.25,2.25) .. controls (2,2.5) and (2,3.5) .. (2.25,3.75);
\end{scope}

\begin{scope}[shift={(-3.25,1)}]
\draw[line width=0.03in] (3.5,2.75) -- (3.5,3.25);
\draw[thick,postaction={decorate}] (3.05,2.25) -- (3.45,2.75);
\draw[thick,postaction={decorate}] (3.95,2.25) -- (3.55,2.75);

\draw[thick,->] (3.45,3.25) -- (3.05,3.75);
\draw[thick,->] (3.55,3.25) -- (3.95,3.75);
\end{scope}

\begin{scope}[shift={(3.65,-0.35)}]

\draw[thick,->] (1.25,1.6) -- (0,3.15);

\draw[thick] (0,1.6) -- (0.55,2.27);
\draw[thick,->] (0.7,2.46) -- (1.25,3.15);
\end{scope}

\end{scope}

\end{tikzpicture}
    \caption{Foam $F$, a sort of ``singular saddle'' cobordism between webs $\Gamma_0$ and $\Gamma_1$, induces a map of state spaces for these webs which enters the inductive construction of link homology groups via cones or SES shown on the right. }
    \label{mfig_016}
\end{figure}
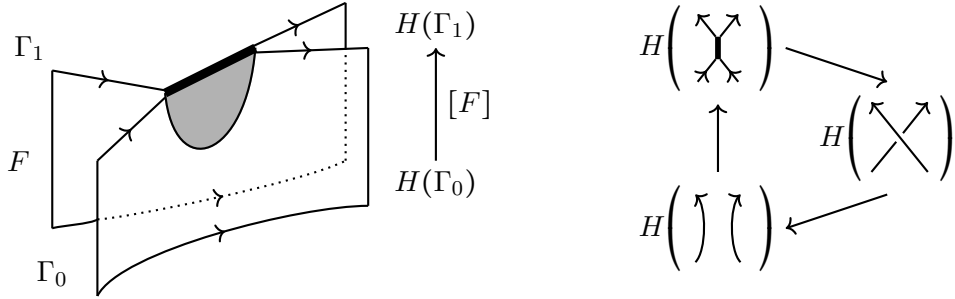

%% file: mfig_021.tex
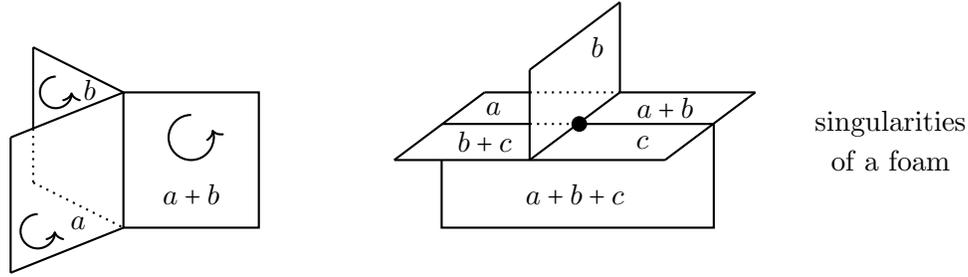
\begin{figure}
    \centering
\begin{tikzpicture}[scale=0.6]
\begin{scope}[shift={(0,0)}]

\draw[thick] (2,0) rectangle (5,3);
\node at (3.5,0.75) {$a+b$};

\draw[thick] (0,4) -- (2,3);
\draw[thick,dotted] (0,1) -- (2,0);
\node at (1.25,3.05) {$b$};

\draw[thick] (0,2.2) -- (0,4);
\draw[thick,dotted] (0,2.2) -- (0,1);

\draw[thick] (-0.5,2) -- (2,3);
\draw[thick] (-0.5,-1) -- (2,0);
\draw[thick] (-0.5,2) -- (-0.5,-1);
\node at (1,0.1) {$a$};

\draw[thick,->] (0.5,3.35) arc (90:360:0.35);

\draw[thick,->] (0.10,0.30) arc (90:360:0.38);

\draw[thick,->] (3.5,2.50) arc (90:380:0.5);

\end{scope}

\begin{scope}[shift={(9,0)}]

\draw[thick] (1,3) -- (2,3);
\draw[thick] (-1,1.5) -- (5,1.5);
\draw[thick] (-1,1.5) -- (1,3);

\draw[thick,dotted] (2,3) -- (4,3);
\draw[thick] (5,1.5) -- (7,3);

\draw[thick] (0.05,2.3) -- (2,2.3);
\draw[thick,dotted] (2,2.3) -- (3,2.3);

\draw[thick] (3.08,2.3) -- (6.08,2.3);

\draw[thick] (4,3) -- (4,5);
\draw[thick] (2,1.5) -- (2,3.5);

\draw[thick] (2,3.5) -- (4,5);
\draw[thick] (2,1.5) -- (4,3);

\draw[thick,fill] (3.25,2.3) arc (0:360:1.5mm);

\draw[thick] (4,3) -- (7,3);

\node at (1.2,2.65) {$a$};
\node at (3.5,4) {$b$};
\node at (1,1.9) {$b+c$};
\node at (5,2.65) {$a+b$};

\node at (4.5,1.9) {$c$};

\draw[thick] (0.05,1.5) -- (0.05,0) -- (6.08,0) -- (6.08,2.3);

\node at (3,0.75) {$a+b+c$};

\node at (10,2.3) {singularities};
\node at (10,1.5) {of a foam};
\end{scope}

\end{tikzpicture}
    \caption{Left: an $(a,b)$ seam of a $\mathsf{GL}_N$ foam and adjacent three facets. Facets are oriented in a compatible way, with orientation preserved when moving from an $a$ or $b$ facet to the adjacent $a+b$ facet. Orientation is reversed between adjacent $a$ and $b$ facets. Right: a singular vertex where four seams meet. }
    \label{mfig_021}
\end{figure}

%% file: mfig_000.tex
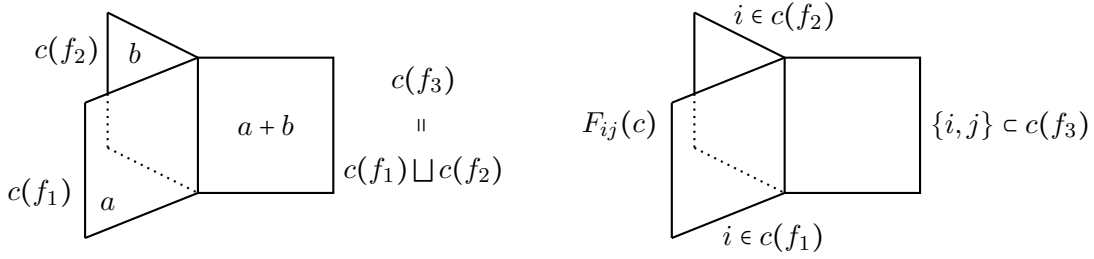
\begin{figure}
    \centering
\begin{tikzpicture}[scale=0.6]
\begin{scope}[shift={(0,0)}]

\draw[thick] (2,0) rectangle (5,3);
\node at (3.5,1.5) {$a+b$};

\draw[thick] (0,4) -- (2,3);
\draw[thick,dotted] (0,1) -- (2,0);
\node at (0.6,3.10) {$b$};

\draw[thick] (0,2.2) -- (0,4);
\draw[thick,dotted] (0,2.2) -- (0,1);

\draw[thick] (-0.5,2) -- (2,3);
\draw[thick] (-0.5,-1) -- (2,0);
\draw[thick] (-0.5,2) -- (-0.5,-1);
\node at (0,-0.25) {$a$};

\node at (-0.9,3.15) {$c(f_2)$};
\node at (-1.5,0) {$c(f_1)$};

\node at (7,2.50) {$c(f_3)$};
\node at (7,1.50) {\rotatebox[origin=c]{90}{$=$}};
\node at (7,0.50) {$c(f_1)\bigsqcup c(f_2)$};

\end{scope}

\begin{scope}[shift={(13,0)}]
\draw[thick] (2,0) rectangle (5,3);
\node at (7.0,1.5) {$\{i,j\}\subset c(f_3)$};

\node at (-1.65,1.5) {$F_{ij}(c)$};

\draw[thick] (0,4) -- (2,3);
\draw[thick,dotted] (0,1) -- (2,0);
\node at (2.00,3.9) {$i\in c(f_2)$};

\draw[thick] (0,2.2) -- (0,4);
\draw[thick,dotted] (0,2.2) -- (0,1);

\draw[thick] (-0.5,2) -- (2,3);
\draw[thick] (-0.5,-1) -- (2,0);
\draw[thick] (-0.5,2) -- (-0.5,-1);
\node at (1.75,-1) {$i\in c(f_1)$};

\end{scope}

\end{tikzpicture}
    \caption{Left: Facets $f_1$, $f_2$ of thickness $a$ and $b$, respectively, merging along a seam are colored by disjoint subsets $c(f_1)$ and $c(f_2)$ of $\{1,\dots, N\}$ of cardinality $a$, respectively $b$. Right: if $i\in c(f_1),j\in c(f_2)$, the two thin facets belong to $F_{ij}(c)$ but the thick facet does not.}
    \label{mfig_000}
\end{figure}

%% file: mfig_017.tex
\begin{figure}
    \centering
\begin{tikzpicture}[scale=0.6,decoration={
    markings,
    mark=at position 0.5 with {\arrow{>}}}]
\begin{scope}[shift={(0,0)}]

\draw[thick] (2,0) rectangle (4,3);
\node at (4.5,1.00) {$f_3$};

\draw[thick] (0,4) -- (2,3);
\draw[thick,dotted] (0,1) -- (2,0);
\node at (-0.5,3) {$f_2$};

\draw[thick] (0,2.2) -- (0,4);
\draw[thick,dotted] (0,2.2) -- (0,1);

\draw[thick] (-0.5,2) -- (2,3);
\draw[thick] (-0.5,-1) -- (2,0);
\draw[thick] (-0.5,2) -- (-0.5,-1);
\node at (-1,0.5) {$f_1$};

\node at (1.25,4) {$c(f_2)$};
\node at (1,-1) {$c(f_1)$};

\node at (5,4) {$c(f_1)\bigsqcup c(f_2)$};
\node at (5,3.25) {\rotatebox[origin=c]{90}{$=$}};
\node at (5,2.50) {$c(f_3)$};

\end{scope}

\begin{scope}[shift={(8.75,0)}]

\draw[thick] (2,0) rectangle (4,3);
\node at (3.5,2.5) {$3$};

\draw[thick] (0,4) -- (2,3);
\draw[thick,dotted] (0,1) -- (2,0);
\node at (0.5,3.2) {$1$};

\draw[thick] (0,2.2) -- (0,4);
\draw[thick,dotted] (0,2.2) -- (0,1);

\draw[thick] (-0.5,2) -- (2,3);
\draw[thick] (-0.5,-1) -- (2,0);
\draw[thick] (-0.5,2) -- (-0.5,-1);
\node at (0,-0.25) {$2$};

\node at (1.25,4) {$\{ 5 \}$};
\node at (1.25,-1.1) {$\{ 1,4\}$};
\node at (5.3,1.5) {$\{1,4,5\}$};
\end{scope}

\begin{scope}[shift={(18.00,0)}]

\draw[thick,fill,black!25] (0,1) -- (2,0) -- (2,3) -- (0,4) -- (0,1);

\draw[thick,fill,black!10] (-0.5,-1) -- (2,0) -- (2,3) -- (-0.5,2) -- (-0.5,-1);

\draw[thick] (2,0) rectangle (4,3);

\draw[thick] (0,4) -- (2,3);
\draw[thick,dotted] (0,1) -- (2,0);
 
\draw[thick] (0,2.2) -- (0,4);
\draw[thick,dotted] (0,2.2) -- (0,1);

\draw[thick] (-0.5,2) -- (2,3);
\draw[thick] (-0.5,-1) -- (2,0);
\draw[thick] (-0.5,2) -- (-0.5,-1);

\node at (3.25,3.5) {$F_{45}(c)$};

\end{scope}

\end{tikzpicture}
    \caption{Left: Thin facets $f_1,f_2$ are colored by disjoint subsets $c(f_1),c(f_2)$. Thick facet $f_3$ is colored by the union of these subsets. Middle: an example of a coloring when $f_1,f_2$ have thickness $2,1$ and are colored by disjoint subsets $\{1,4\}$ and $\{5\}$, respectively. Facet $f_3$ is then colored by their union $\{1,4,5\}$. Right: part of the bicolored surface $F_{45}(c)$ for this coloring, containing both thin but not the thick facet along that seam.}
    \label{mfig_017}
\end{figure}
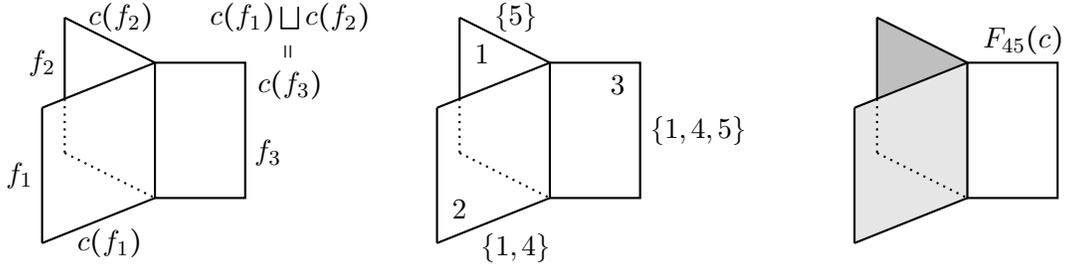

%% file: mfig_X1.tex
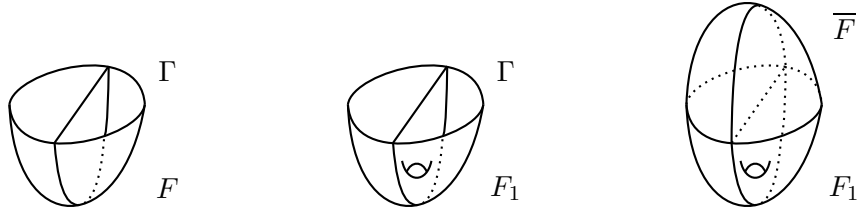
\begin{figure}
    \centering

\begin{tikzpicture}[scale=0.6]
\begin{scope}[shift={(0,0)}]
\draw[thick] (0,2) .. controls (0.1,2.75) and (2.9,3.5) .. (3,2);
\draw[thick] (0,2) .. controls (0.1,0.5) and (2.9,1.25) .. (3,2);
\draw[thick] (1,1.15) -- (2.2,2.85);

\draw[thick] (0,2) .. controls (0.1,-1) and (2.5,-1) .. (3,2);

\draw[thick] (1,1.15) .. controls (1.1,-0.2) and (1.4,-0.2) .. (1.5,-0.23);

\draw[thick] (2.2,2.85) .. controls (2.18,1.4) and (2.12,1.3) .. (2.1,1.3);
\draw[thick,dotted] (1.5,-0.23) .. controls (1.6,-0.23) and (2,-0.23) .. (2.1,1.3);

\node at (3.5,2.75) {$\Gamma$};
\node at (3.5,0.15) {$F$};

\end{scope}

\begin{scope}[shift={(7.5,0)}]
\draw[thick] (0,2) .. controls (0.1,2.75) and (2.9,3.5) .. (3,2);
\draw[thick] (0,2) .. controls (0.1,0.5) and (2.9,1.25) .. (3,2);
\draw[thick] (1,1.15) -- (2.2,2.85);

\draw[thick] (0,2) .. controls (0.1,-1) and (2.5,-1) .. (3,2);

\draw[thick] (1,1.15) .. controls (1.1,-0.2) and (1.4,-0.2) .. (1.5,-0.23);

\draw[thick] (2.2,2.85) .. controls (2.18,1.4) and (2.12,1.3) .. (2.1,1.3);
\draw[thick,dotted] (1.5,-0.23) .. controls (1.6,-0.23) and (2,-0.23) .. (2.1,1.3);

\node at (3.5,2.75) {$\Gamma$};
\node at (3.5,0.15) {$F_1$};

\draw[thick] (1.2,0.75) .. controls (1.3,0.25) and (1.75,0.25) .. (1.85,0.75);

\draw[thick] (1.3,0.5) .. controls (1.4,0.75) and (1.65,0.75) .. (1.75,0.5);
\end{scope}

\begin{scope}[shift={(15,0)}]
\draw[thick,dotted] (0,2) .. controls (0.1,2.75) and (2.9,3.5) .. (3,2);
\draw[thick] (0,2) .. controls (0.1,0.5) and (2.9,1.25) .. (3,2);
\draw[thick,dotted] (1,1.15) -- (2.2,2.85);

\draw[thick] (0,2) .. controls (0.1,-1) and (2.5,-1) .. (3,2);

\draw[thick] (1,1.15) .. controls (1.1,-0.2) and (1.4,-0.2) .. (1.5,-0.23);

\draw[thick,dotted] (2.2,2.85) .. controls (2.18,1.4) and (2.12,1.3) .. (2.1,1.3);
\draw[thick,dotted] (1.5,-0.23) .. controls (1.6,-0.23) and (2,-0.23) .. (2.1,1.3);

\node at (3.5,3.75) {$\overline{F}$};
\node at (3.5,0.15) {$F_1$};

\draw[thick] (1.2,0.75) .. controls (1.3,0.25) and (1.75,0.25) .. (1.85,0.75);

\draw[thick] (1.3,0.5) .. controls (1.4,0.75) and (1.65,0.75) .. (1.75,0.5);

\draw[thick] (0,2) .. controls (0.1,5) and (2.5,5) .. (3,2);

\draw[thick] (1,1.15) .. controls (1.01,4.2) and (1.5,4.2) .. (1.6,4.2);
\draw[thick,dotted] (1.6,4.2) .. controls (1.7,4.1) and (2.1,4.1) .. (2.2,2.85);

\end{scope}

\end{tikzpicture}
    \caption{Left: A foam $F$ with boundary $\Gamma$. Middle: another foam $F_1$ with $\partial F_1 = \Gamma$. Right: gluing $F$ and $F_1$ along $\Gamma$ into a closed foam. Here $\overline{F}$ is $F$ reflected in the $z$-plane.}
    \label{mfig_X1}
\end{figure}